\documentclass[10pt,a4,twoside,rm]{article}
\usepackage[left=4.5cm,top=4cm,right=4.5cm,bottom=2cm]{geometry}
\usepackage{amsmath, amsthm, amssymb}
\usepackage{graphicx}
\usepackage[all, knot]{xy}
\usepackage{amssymb}
\usepackage{fancyhdr}
\usepackage{titlesec}
\usepackage{multirow}
\usepackage{youngtab}
\usepackage{pstricks}
\usepackage{psfrag}

\def\linew(#1){%
/xywidth{#1 setlinewidth}def}

\def\union{\xy 0;/r.07pc/:
{\ar@{<->} (0,0)*{};(50,0)*{}};
(25,10)*{\Huge ? };
\endxy}

\def\flechahaciaabajo{\xy 0;/r.03pc/:
{\ar@{->} (0,0)*{};(0,-30)*{}};
\endxy}

\def\action{\xy 0;/r1pc/:
(-12,0); (12,0) **\dir{-};
(-7,1); (-7,-1); **\dir{.};
(-2,1); (-2,-1) **\dir{.};
(3,1); (3,-1) **\dir{.};
(8,1); (8,-1) **\dir{.};
(0,.5); (0,-.5) **\dir{-};
(6,.5); (6,-.5) **\dir{-};
(10,.5); (10,-.5) **\dir{-};
(-4,.5); (-4,-.5) **\dir{-};
(-10,.5); (-10,-.5) **\dir{-};
(13,0)*{\cdots};
(-13,0)*{\cdots};
{\ar@/^1pc/(0.2,1)*{};(5.8,1)*{}};
{\ar@/^1pc/(6.2,1)*{};(9.8,1)*{}};
{\ar@/^1pc/(-3.8,1)*{};(-0.2,1)*{}};
{\ar@/^1pc/(-9.8,1)*{};(-4.2,1)*{}};
(0,-1.5)*{\scriptstyle \bf 0};
(-4,-1.5)*{\scriptstyle \bf -4};
(-10,-1.5)*{\scriptstyle \bf-10};
(6,-1.5)*{ \scriptstyle \bf 6};
(10,-1.5)*{\scriptstyle \bf 10};
(-7,-1.5)*{ \scriptstyle -7};
(-2,-1.5)*{\scriptstyle -2};
(3,-1.5)*{ \scriptstyle 3};
(8,-1.5)*{\scriptstyle  8};
\endxy}

\def\ejemploprimero{\xy 0;/r.8pc/:
(0,-8)*{ w(\sub{t}_0)};
(22,-8)*{ w(\sub{t}_1)};
(0,0)*{ \bullet};
(-1,-1)*{ \bullet};
(1,-1)*{ \bullet};
(-2,-2)*{ \bullet};
(2,-2)*{ \bullet};
(0,-2)*{ \bullet};
(-3,-3)*{ \bullet};
(3,-3)*{ \bullet};
(-1,-3)*{ \bullet};
(1,-3)*{ \bullet};
(4,-4)*{ \bullet};
(-4,-4)*{ \bullet};
(2,-4)*{ \bullet};
(-2,-4)*{ \bullet};
(0,-4)*{ \bullet};
(-3,-5)*{ \bullet};
(3,-5)*{ \bullet};
(-1,-5)*{ \bullet};
(1,-5)*{ \bullet};
(-5,-5)*{ \bullet};
(5,-5)*{ \bullet};
(4,-6)*{ \bullet};
(-4,-6)*{ \bullet};
(2,-6)*{ \bullet};
(-2,-6)*{ \bullet};
(0,-6)*{ \bullet};
(6,-6)*{ \bullet};
(-6,-6)*{ \bullet};
(0,0); (-1,-1) **\dir{-};
(0,-2); (-1,-1) **\dir{-};
(0,-2); (-1,-3) **\dir{-};
(0,-4); (-1,-3) **\dir{-};
(0,-4); (2,-6) **\dir{-};
(22,0)*{ \bullet};
(21,-1)*{ \bullet};
(23,-1)*{ \bullet};
(20,-2)*{ \bullet};
(24,-2)*{ \bullet};
(22,-2)*{ \bullet};
(21,-3)*{ \bullet};
(23,-3)*{ \bullet};
(19,-3)*{ \bullet};
(25,-3)*{ \bullet};
(26,-4)*{ \bullet};
(18,-4)*{ \bullet};
(24,-4)*{ \bullet};
(20,-4)*{ \bullet};
(22,-4)*{ \bullet};
(19,-5)*{ \bullet};
(25,-5)*{ \bullet};
(21,-5)*{ \bullet};
(23,-5)*{ \bullet};
(17,-5)*{ \bullet};
(27,-5)*{ \bullet};
(26,-6)*{ \bullet};
(18,-6)*{ \bullet};
(24,-6)*{ \bullet};
(20,-6)*{ \bullet};
(22,-6)*{ \bullet};
(28,-6)*{ \bullet};
(16,-6)*{ \bullet};
(22,0); (23,-1) **\dir{-};
(22,-2); (23,-1) **\dir{-};
(22,-2); (21,-3) **\dir{-};
(22,-4); (21,-3) **\dir{-};
(22,-4); (24,-6) **\dir{-};
\endxy}

\def\ejemplosegundo{\xy 0;/r.8pc/:
(0,-8)*{ w(\sub{t}_2)};
(22,-8)*{ w(\sub{t}_3)};
(0,0)*{ \bullet};
(-1,-1)*{ \bullet};
(1,-1)*{ \bullet};
(-2,-2)*{ \bullet};
(2,-2)*{ \bullet};
(0,-2)*{ \bullet};
(-3,-3)*{ \bullet};
(3,-3)*{ \bullet};
(-1,-3)*{ \bullet};
(1,-3)*{ \bullet};
(4,-4)*{ \bullet};
(-4,-4)*{ \bullet};
(2,-4)*{ \bullet};
(-2,-4)*{ \bullet};
(0,-4)*{ \bullet};
(-3,-5)*{ \bullet};
(3,-5)*{ \bullet};
(-1,-5)*{ \bullet};
(1,-5)*{ \bullet};
(-5,-5)*{ \bullet};
(5,-5)*{ \bullet};
(4,-6)*{ \bullet};
(-4,-6)*{ \bullet};
(2,-6)*{ \bullet};
(-2,-6)*{ \bullet};
(0,-6)*{ \bullet};
(6,-6)*{ \bullet};
(-6,-6)*{ \bullet};
(0,0); (1,-1) **\dir{-};
(0,-2); (1,-1) **\dir{-};
(0,-2); (1,-3) **\dir{-};
(0,-4); (1,-3) **\dir{-};
(0,-4); (2,-6) **\dir{-};
(22,0)*{ \bullet};
(21,-1)*{ \bullet};
(23,-1)*{ \bullet};
(20,-2)*{ \bullet};
(24,-2)*{ \bullet};
(22,-2)*{ \bullet};
(21,-3)*{ \bullet};
(23,-3)*{ \bullet};
(19,-3)*{ \bullet};
(25,-3)*{ \bullet};
(26,-4)*{ \bullet};
(18,-4)*{ \bullet};
(24,-4)*{ \bullet};
(20,-4)*{ \bullet};
(22,-4)*{ \bullet};
(19,-5)*{ \bullet};
(25,-5)*{ \bullet};
(21,-5)*{ \bullet};
(23,-5)*{ \bullet};
(17,-5)*{ \bullet};
(27,-5)*{ \bullet};
(26,-6)*{ \bullet};
(18,-6)*{ \bullet};
(24,-6)*{ \bullet};
(20,-6)*{ \bullet};
(22,-6)*{ \bullet};
(28,-6)*{ \bullet};
(16,-6)*{ \bullet};
(22,0); (24,-2) **\dir{-};
(22,-4); (24,-2) **\dir{-};
(22,-4); (24,-6) **\dir{-};
\endxy}

\def\ejemplotercero{\xy 0;/r.8pc/:
(0,-8)*{ w(\sub{t}_4)};
(22,-8)*{ w(\sub{t}_5)};
(0,0)*{ \bullet};
(-1,-1)*{ \bullet};
(1,-1)*{ \bullet};
(-2,-2)*{ \bullet};
(2,-2)*{ \bullet};
(0,-2)*{ \bullet};
(-3,-3)*{ \bullet};
(3,-3)*{ \bullet};
(-1,-3)*{ \bullet};
(1,-3)*{ \bullet};
(4,-4)*{ \bullet};
(-4,-4)*{ \bullet};
(2,-4)*{ \bullet};
(-2,-4)*{ \bullet};
(0,-4)*{ \bullet};
(-3,-5)*{ \bullet};
(3,-5)*{ \bullet};
(-1,-5)*{ \bullet};
(1,-5)*{ \bullet};
(-5,-5)*{ \bullet};
(5,-5)*{ \bullet};
(4,-6)*{ \bullet};
(-4,-6)*{ \bullet};
(2,-6)*{ \bullet};
(-2,-6)*{ \bullet};
(0,-6)*{ \bullet};
(6,-6)*{ \bullet};
(-6,-6)*{ \bullet};
(0,0); (2,-2) **\dir{-};
(1,-3); (2,-2) **\dir{-};
(1,-3); (2,-4) **\dir{-};
(1,-5); (2,-4) **\dir{-};
(1,-5); (2,-6) **\dir{-};
(22,0)*{ \bullet};
(21,-1)*{ \bullet};
(23,-1)*{ \bullet};
(20,-2)*{ \bullet};
(24,-2)*{ \bullet};
(22,-2)*{ \bullet};
(21,-3)*{ \bullet};
(23,-3)*{ \bullet};
(19,-3)*{ \bullet};
(25,-3)*{ \bullet};
(26,-4)*{ \bullet};
(18,-4)*{ \bullet};
(24,-4)*{ \bullet};
(20,-4)*{ \bullet};
(22,-4)*{ \bullet};
(19,-5)*{ \bullet};
(25,-5)*{ \bullet};
(21,-5)*{ \bullet};
(23,-5)*{ \bullet};
(17,-5)*{ \bullet};
(27,-5)*{ \bullet};
(26,-6)*{ \bullet};
(18,-6)*{ \bullet};
(24,-6)*{ \bullet};
(20,-6)*{ \bullet};
(22,-6)*{ \bullet};
(28,-6)*{ \bullet};
(16,-6)*{ \bullet};
(22,0); (25,-3) **\dir{-};
(23,-5); (25,-3) **\dir{-};
(23,-5); (24,-6) **\dir{-};
\endxy}

\def\dibujos{\xy 0;/r.5pc/:
\xy 0;/r.6pc/:
(-3,-3); (0,0) **\dir{-};
(-3,-3); (3,-9) **\dir{-};
(-1,-13); (3,-9) **\dir{-};
(0,0)*{ \bullet};
(-1,-1)*{ \bullet};
(1,-1)*{ \bullet};
(-2,-2)*{ \bullet};
(2,-2)*{ \bullet};
(0,-2)*{ \bullet};
(-3,-3)*{ \bullet};
(3,-3)*{ \bullet};
(-1,-3)*{ \bullet};
(1,-3)*{ \bullet};
(4,-4)*{ \bullet};
(-4,-4)*{ \bullet};
(2,-4)*{ \bullet};
(-2,-4)*{ \bullet};
(0,-4)*{ \bullet};
(-3,-5)*{ \bullet};
(3,-5)*{ \bullet};
(-1,-5)*{ \bullet};
(1,-5)*{ \bullet};
(-5,-5)*{ \bullet};
(5,-5)*{ \bullet};
(4,-6)*{ \bullet};
(-4,-6)*{ \bullet};
(2,-6)*{ \bullet};
(-2,-6)*{ \bullet};
(0,-6)*{ \bullet};
(6,-6)*{ \bullet};
(-6,-6)*{ \bullet};
(-3,-7)*{ \bullet};
(3,-7)*{ \bullet};
(-1,-7)*{ \bullet};
(1,-7)*{ \bullet};
(-5,-7)*{ \bullet};
(5,-7)*{ \bullet};
(-7,-7)*{ \bullet};
(7,-7)*{ \bullet};
(4,-8)*{ \bullet};
(-4,-8)*{ \bullet};
(2,-8)*{ \bullet};
(-2,-8)*{ \bullet};
(0,-8)*{ \bullet};
(6,-8)*{ \bullet};
(-6,-8)*{ \bullet};
(8,-8)*{ \bullet};
(-8,-8)*{ \bullet};
(-3,-9)*{ \bullet};
(3,-9)*{ \bullet};
(-1,-9)*{ \bullet};
(1,-9)*{ \bullet};
(-5,-9)*{ \bullet};
(5,-9)*{ \bullet};
(-7,-9)*{ \bullet};
(7,-9)*{ \bullet};
(-9,-9)*{ \bullet};
(9,-9)*{ \bullet};
(4,-10)*{ \bullet};
(-4,-10)*{ \bullet};
(2,-10)*{ \bullet};
(-2,-10)*{ \bullet};
(0,-10)*{ \bullet};
(6,-10)*{ \bullet};
(-6,-10)*{ \bullet};
(8,-10)*{ \bullet};
(-8,-10)*{ \bullet};
(10,-10)*{ \bullet};
(-10,-10)*{ \bullet};
(-3,-11)*{ \bullet};
(3,-11)*{ \bullet};
(-1,-11)*{ \bullet};
(1,-11)*{ \bullet};
(-5,-11)*{ \bullet};
(5,-11)*{ \bullet};
(-7,-11)*{ \bullet};
(7,-11)*{ \bullet};
(-9,-11)*{ \bullet};
(9,-11)*{ \bullet};
(-11,-11)*{ \bullet};
(11,-11)*{ \bullet};
(4,-12)*{ \bullet};
(-4,-12)*{ \bullet};
(2,-12)*{ \bullet};
(-2,-12)*{ \bullet};
(0,-12)*{ \bullet};
(6,-12)*{ \bullet};
(-6,-12)*{ \bullet};
(8,-12)*{ \bullet};
(-8,-12)*{ \bullet};
(10,-12)*{ \bullet};
(-10,-12)*{ \bullet};
(12,-12)*{ \bullet};
(-12,-12)*{ \bullet};
(-3,-13)*{ \bullet};
(3,-13)*{ \bullet};
(-1,-13)*{ \bullet};
(1,-13)*{ \bullet};
(-5,-13)*{ \bullet};
(5,-13)*{ \bullet};
(-7,-13)*{ \bullet};
(7,-13)*{ \bullet};
(-9,-13)*{ \bullet};
(9,-13)*{ \bullet};
(-11,-13)*{ \bullet};
(11,-13)*{ \bullet};
(-13,-13)*{ \bullet};
(13,-13)*{ \bullet};
(-12,0); (-12,-14.5) **\dir{.};
(-7,0); (-7,-14.5) **\dir{.};
(-2,0); (-2,-14.5) **\dir{.};
(3,0); (3,-14.5) **\dir{.};
(8,0); (8,-14.5) **\dir{.};
(13,0); (13,-14.5) **\dir{.};
{\ar@{->} (-14,-15)*{};(14,-15)*{}};
(-12,-15.5)*{\scriptscriptstyle -12};
(-7,-15.5)*{\scriptscriptstyle -7};
(-2,-15.5)*{\scriptscriptstyle -2};
(3,-15.5)*{\scriptscriptstyle 3};
(8,-15.5)*{\scriptscriptstyle 8};
(13,-15.5)*{\scriptscriptstyle 13};
(0,-14.5)*{\scriptscriptstyle  \text{weight}};
\endxy}

\def\dibujoss{\xy 0;/r.6pc/:
(0,0); (-13,-13) **\dir{-};
(-11,-13); (-12,-12) **\dir{-};
(-7,-7); (-2,-12) **\dir{-};
(-3,-13); (-2,-12) **\dir{-};
(-1,-13); (-2,-12) **\dir{-};
(-2,-2); (9,-13) **\dir{-};
(3,-7); (-2,-12) **\dir{-};
(7,-13); (8,-12) **\dir{-};
(0,0)*{ \bullet};
(-1,-1)*{ \bullet};
(1,-1)*{ \bullet};
(-2,-2)*{ \bullet};
(2,-2)*{ \bullet};
(0,-2)*{ \bullet};
(-3,-3)*{ \bullet};
(3,-3)*{ \bullet};
(-1,-3)*{ \bullet};
(1,-3)*{ \bullet};
(4,-4)*{ \bullet};
(-4,-4)*{ \bullet};
(2,-4)*{ \bullet};
(-2,-4)*{ \bullet};
(0,-4)*{ \bullet};
(-3,-5)*{ \bullet};
(3,-5)*{ \bullet};
(-1,-5)*{ \bullet};
(1,-5)*{ \bullet};
(-5,-5)*{ \bullet};
(5,-5)*{ \bullet};
(4,-6)*{ \bullet};
(-4,-6)*{ \bullet};
(2,-6)*{ \bullet};
(-2,-6)*{ \bullet};
(0,-6)*{ \bullet};
(6,-6)*{ \bullet};
(-6,-6)*{ \bullet};
(-3,-7)*{ \bullet};
(3,-7)*{ \bullet};
(-1,-7)*{ \bullet};
(1,-7)*{ \bullet};
(-5,-7)*{ \bullet};
(5,-7)*{ \bullet};
(-7,-7)*{ \bullet};
(7,-7)*{ \bullet};
(4,-8)*{ \bullet};
(-4,-8)*{ \bullet};
(2,-8)*{ \bullet};
(-2,-8)*{ \bullet};
(0,-8)*{ \bullet};
(6,-8)*{ \bullet};
(-6,-8)*{ \bullet};
(8,-8)*{ \bullet};
(-8,-8)*{ \bullet};
(-3,-9)*{ \bullet};
(3,-9)*{ \bullet};
(-1,-9)*{ \bullet};
(1,-9)*{ \bullet};
(-5,-9)*{ \bullet};
(5,-9)*{ \bullet};
(-7,-9)*{ \bullet};
(7,-9)*{ \bullet};
(-9,-9)*{ \bullet};
(9,-9)*{ \bullet};
(4,-10)*{ \bullet};
(-4,-10)*{ \bullet};
(2,-10)*{ \bullet};
(-2,-10)*{ \bullet};
(0,-10)*{ \bullet};
(6,-10)*{ \bullet};
(-6,-10)*{ \bullet};
(8,-10)*{ \bullet};
(-8,-10)*{ \bullet};
(10,-10)*{ \bullet};
(-10,-10)*{ \bullet};
(-3,-11)*{ \bullet};
(3,-11)*{ \bullet};
(-1,-11)*{ \bullet};
(1,-11)*{ \bullet};
(-5,-11)*{ \bullet};
(5,-11)*{ \bullet};
(-7,-11)*{ \bullet};
(7,-11)*{ \bullet};
(-9,-11)*{ \bullet};
(9,-11)*{ \bullet};
(-11,-11)*{ \bullet};
(11,-11)*{ \bullet};
(4,-12)*{ \bullet};
(-4,-12)*{ \bullet};
(2,-12)*{ \bullet};
(-2,-12)*{ \bullet};
(0,-12)*{ \bullet};
(6,-12)*{ \bullet};
(-6,-12)*{ \bullet};
(8,-12)*{ \bullet};
(-8,-12)*{ \bullet};
(10,-12)*{ \bullet};
(-10,-12)*{ \bullet};
(12,-12)*{ \bullet};
(-12,-12)*{ \bullet};
(-3,-13)*{ \bullet};
(3,-13)*{ \bullet};
(-1,-13)*{ \bullet};
(1,-13)*{ \bullet};
(-5,-13)*{ \bullet};
(5,-13)*{ \bullet};
(-7,-13)*{ \bullet};
(7,-13)*{ \bullet};
(-9,-13)*{ \bullet};
(9,-13)*{ \bullet};
(-11,-13)*{ \bullet};
(11,-13)*{ \bullet};
(-13,-13)*{ \bullet};
(13,-13)*{ \bullet};
(-12,0); (-12,-14.5) **\dir{.};
(-7,0); (-7,-14.5) **\dir{.};
(-2,0); (-2,-14.5) **\dir{.};
(3,0); (3,-14.5) **\dir{.};
(8,0); (8,-14.5) **\dir{.};
(13,0); (13,-14.5) **\dir{.};
{\ar@{->} (-14,-15)*{};(14,-15)*{}};
(-12,-15.5)*{\scriptscriptstyle -12};
(-7,-15.5)*{\scriptscriptstyle -7};
(-2,-15.5)*{\scriptscriptstyle -2};
(3,-15.5)*{\scriptscriptstyle 3};
(8,-15.5)*{\scriptscriptstyle 8};
(13,-15.5)*{\scriptscriptstyle 13};
(0,-14.5)*{\scriptscriptstyle  \text{weight}};
\endxy}

\def\pascal{\xy 0;/r2pc/:
(0,1); (0,-3) **\dir{.};
{\ar@{->} (0,.5)*{};(2,.5)*{}};
(4.3,0.5)*{ \scriptstyle \text{central vertical axis}};
(0,0)*{ \bullet};
(-1,-1)*{ \bullet};
(1,-1)*{ \bullet};
(-2,-2)*{ \bullet};
(2,-2)*{ \bullet};
(0,-2)*{ \bullet};
(-3,-3)*{ \bullet};
(3,-3)*{ \bullet};
(-1,-3)*{ \bullet};
(1,-3)*{ \bullet};
{\ar@{->} (-4,-5)*{};(4,-5)*{}};
(-3,-5.5)*{\scriptstyle -3};
(-2,-5.5)*{\scriptstyle -2};
(-1,-5.5)*{\scriptstyle -1};
(0,-5.5)*{\scriptstyle 0};
(1,-5.5)*{\scriptstyle1};
(2,-5.5)*{ \scriptstyle 2};
(3,-5.5)*{\scriptstyle 3};
(0,-4.5)*{\scriptstyle  \text{weight}};
(4,-5.5)*{\ldots};
(-4,-5.5)*{\ldots};
{\ar@{->} (-6,0)*{};(-6,-5)*{}};
(-5.5,0)*{\scriptstyle 0};
(-5.5,-1)*{\scriptstyle 1};
(-5.5,-2)*{\scriptstyle 2};
(-5.5,-3)*{\scriptstyle 3};
(-7,-1.5)*{\scriptstyle  \text{level}};
\endxy
}

\def\mayor{\xy 0;/r1.3pc/:
(0,0)*{ \bullet};
(-1,-1)*{ \bullet};
(1,-1)*{ \bullet};
(-2,-2)*{ \bullet};
(2,-2)*{ \bullet};
(0,-2)*{ \bullet};
(-3,-3)*{ \bullet};
(3,-3)*{ \bullet};
(-1,-3)*{ \bullet};
(1,-3)*{ \bullet};
(4,-4)*{ \bullet};
(-4,-4)*{ \bullet};
(2,-4)*{ \bullet};
(-2,-4)*{ \bullet};
(0,-4)*{ \bullet};
(-3,-5)*{ \bullet};
(3,-5)*{ \bullet};
(-1,-5)*{ \bullet};
(1,-5)*{ \bullet};
(-5,-5)*{ \bullet};
(5,-5)*{ \bullet};
(4,-6)*{ \bullet};
(-4,-6)*{ \bullet};
(2,-6)*{ \bullet};
(-2,-6)*{ \bullet};
(0,-6)*{ \bullet};
(6,-6)*{ \bullet};
(-6,-6)*{ \bullet};
(0,1); (0,-6) **\dir{.};
(0,0); (-1,-1) **\dir{-};
(0,-2); (-1,-1) **\dir{-};
(0,-2); (-1,-3) **\dir{-};
(0,-4); (-1,-3) **\dir{-};
(0,-4); (1,-5) **\dir{-};
(2,-6); (1,-5) **\dir{-};
{\ar@{->} (-7,-8)*{};(7,-8)*{}};
(-6,-8.5)*{\scriptstyle -6};
(-5,-8.5)*{\scriptstyle -5};
(-4,-8.5)*{\scriptstyle -4};
(-3,-8.5)*{\scriptstyle -3};
(-2,-8.5)*{\scriptstyle -2};
(-1,-8.5)*{\scriptstyle -1};
(0,-8.5)*{\scriptstyle 0};
(1,-8.5)*{\scriptstyle1};
(2,-8.5)*{ \scriptstyle 2};
(3,-8.5)*{\scriptstyle 3};
(6,-8.5)*{\scriptstyle 6};
(5,-8.5)*{\scriptstyle 5};
(4,-8.5)*{\scriptstyle 4};
(0,-7.5)*{\scriptstyle  \text{weight}};
{\ar@{->} (-8,0)*{};(-8,-7)*{}};
(-7.5,0)*{\scriptstyle 0};
(-7.5,-1)*{\scriptstyle 1};
(-7.5,-2)*{\scriptstyle 2};
(-7.5,-3)*{\scriptstyle 3};
(-7.5,-4)*{\scriptstyle 4};
(-7.5,-5)*{\scriptstyle 5};
(-7.5,-6)*{\scriptstyle 6};
(-9,-3)*{\scriptstyle  \text{level}};
\endxy }

\def\walls{\xy 0;/r.8pc/:
(-3,-3); (0,0) **\dir{-};
(-3,-3); (3,-9) **\dir{-};
(-1,-13); (3,-9) **\dir{-};
(0,0)*{ \bullet};
(-1,-1)*{ \bullet};
(1,-1)*{ \bullet};
(-2,-2)*{ \bullet};
(2,-2)*{ \bullet};
(0,-2)*{ \bullet};
(-3,-3)*{ \bullet};
(3,-3)*{ \bullet};
(-1,-3)*{ \bullet};
(1,-3)*{ \bullet};
(4,-4)*{ \bullet};
(-4,-4)*{ \bullet};
(2,-4)*{ \bullet};
(-2,-4)*{ \bullet};
(0,-4)*{ \bullet};
(-3,-5)*{ \bullet};
(3,-5)*{ \bullet};
(-1,-5)*{ \bullet};
(1,-5)*{ \bullet};
(-5,-5)*{ \bullet};
(5,-5)*{ \bullet};
(4,-6)*{ \bullet};
(-4,-6)*{ \bullet};
(2,-6)*{ \bullet};
(-2,-6)*{ \bullet};
(0,-6)*{ \bullet};
(6,-6)*{ \bullet};
(-6,-6)*{ \bullet};
(-3,-7)*{ \bullet};
(3,-7)*{ \bullet};
(-1,-7)*{ \bullet};
(1,-7)*{ \bullet};
(-5,-7)*{ \bullet};
(5,-7)*{ \bullet};
(-7,-7)*{ \bullet};
(7,-7)*{ \bullet};
(4,-8)*{ \bullet};
(-4,-8)*{ \bullet};
(2,-8)*{ \bullet};
(-2,-8)*{ \bullet};
(0,-8)*{ \bullet};
(6,-8)*{ \bullet};
(-6,-8)*{ \bullet};
(8,-8)*{ \bullet};
(-8,-8)*{ \bullet};
(-3,-9)*{ \bullet};
(3,-9)*{ \bullet};
(-1,-9)*{ \bullet};
(1,-9)*{ \bullet};
(-5,-9)*{ \bullet};
(5,-9)*{ \bullet};
(-7,-9)*{ \bullet};
(7,-9)*{ \bullet};
(-9,-9)*{ \bullet};
(9,-9)*{ \bullet};
(4,-10)*{ \bullet};
(-4,-10)*{ \bullet};
(2,-10)*{ \bullet};
(-2,-10)*{ \bullet};
(0,-10)*{ \bullet};
(6,-10)*{ \bullet};
(-6,-10)*{ \bullet};
(8,-10)*{ \bullet};
(-8,-10)*{ \bullet};
(10,-10)*{ \bullet};
(-10,-10)*{ \bullet};
(-3,-11)*{ \bullet};
(3,-11)*{ \bullet};
(-1,-11)*{ \bullet};
(1,-11)*{ \bullet};
(-5,-11)*{ \bullet};
(5,-11)*{ \bullet};
(-7,-11)*{ \bullet};
(7,-11)*{ \bullet};
(-9,-11)*{ \bullet};
(9,-11)*{ \bullet};
(-11,-11)*{ \bullet};
(11,-11)*{ \bullet};
(4,-12)*{ \bullet};
(-4,-12)*{ \bullet};
(2,-12)*{ \bullet};
(-2,-12)*{ \bullet};
(0,-12)*{ \bullet};
(6,-12)*{ \bullet};
(-6,-12)*{ \bullet};
(8,-12)*{ \bullet};
(-8,-12)*{ \bullet};
(10,-12)*{ \bullet};
(-10,-12)*{ \bullet};
(12,-12)*{ \bullet};
(-12,-12)*{ \bullet};
(-3,-13)*{ \bullet};
(3,-13)*{ \bullet};
(-1,-13)*{ \bullet};
(1,-13)*{ \bullet};
(-5,-13)*{ \bullet};
(5,-13)*{ \bullet};
(-7,-13)*{ \bullet};
(7,-13)*{ \bullet};
(-9,-13)*{ \bullet};
(9,-13)*{ \bullet};
(-11,-13)*{ \bullet};
(11,-13)*{ \bullet};
(-13,-13)*{ \bullet};
(13,-13)*{ \bullet};
(-12,0); (-12,-14.5) **\dir{.};
(-7,0); (-7,-14.5) **\dir{.};
(-2,0); (-2,-14.5) **\dir{.};
(3,0); (3,-14.5) **\dir{.};
(8,0); (8,-14.5) **\dir{.};
(13,0); (13,-14.5) **\dir{.};
{\ar@{->} (-14,-15)*{};(14,-15)*{}};
(-12,-15.5)*{\scriptscriptstyle -12};
(-7,-15.5)*{\scriptscriptstyle -7};
(-2,-15.5)*{\scriptscriptstyle -2};
(3,-15.5)*{\scriptscriptstyle 3};
(8,-15.5)*{\scriptscriptstyle 8};
(13,-15.5)*{\scriptscriptstyle 13};
(0,-14.5)*{\scriptscriptstyle  \text{weight}};
\endxy
}

\def\maximales{\xy 0;/r.8pc/:
(0,0); (-13,-13) **\dir{-};
(-11,-13); (-12,-12) **\dir{-};
(-7,-7); (-2,-12) **\dir{-};
(-3,-13); (-2,-12) **\dir{-};
(-1,-13); (-2,-12) **\dir{-};
(-2,-2); (9,-13) **\dir{-};
(3,-7); (-2,-12) **\dir{-};
(7,-13); (8,-12) **\dir{-};
(0,0)*{ \bullet};
(-1,-1)*{ \bullet};
(1,-1)*{ \bullet};
(-2,-2)*{ \bullet};
(2,-2)*{ \bullet};
(0,-2)*{ \bullet};
(-3,-3)*{ \bullet};
(3,-3)*{ \bullet};
(-1,-3)*{ \bullet};
(1,-3)*{ \bullet};
(4,-4)*{ \bullet};
(-4,-4)*{ \bullet};
(2,-4)*{ \bullet};
(-2,-4)*{ \bullet};
(0,-4)*{ \bullet};
(-3,-5)*{ \bullet};
(3,-5)*{ \bullet};
(-1,-5)*{ \bullet};
(1,-5)*{ \bullet};
(-5,-5)*{ \bullet};
(5,-5)*{ \bullet};
(4,-6)*{ \bullet};
(-4,-6)*{ \bullet};
(2,-6)*{ \bullet};
(-2,-6)*{ \bullet};
(0,-6)*{ \bullet};
(6,-6)*{ \bullet};
(-6,-6)*{ \bullet};
(-3,-7)*{ \bullet};
(3,-7)*{ \bullet};
(-1,-7)*{ \bullet};
(1,-7)*{ \bullet};
(-5,-7)*{ \bullet};
(5,-7)*{ \bullet};
(-7,-7)*{ \bullet};
(7,-7)*{ \bullet};
(4,-8)*{ \bullet};
(-4,-8)*{ \bullet};
(2,-8)*{ \bullet};
(-2,-8)*{ \bullet};
(0,-8)*{ \bullet};
(6,-8)*{ \bullet};
(-6,-8)*{ \bullet};
(8,-8)*{ \bullet};
(-8,-8)*{ \bullet};
(-3,-9)*{ \bullet};
(3,-9)*{ \bullet};
(-1,-9)*{ \bullet};
(1,-9)*{ \bullet};
(-5,-9)*{ \bullet};
(5,-9)*{ \bullet};
(-7,-9)*{ \bullet};
(7,-9)*{ \bullet};
(-9,-9)*{ \bullet};
(9,-9)*{ \bullet};
(4,-10)*{ \bullet};
(-4,-10)*{ \bullet};
(2,-10)*{ \bullet};
(-2,-10)*{ \bullet};
(0,-10)*{ \bullet};
(6,-10)*{ \bullet};
(-6,-10)*{ \bullet};
(8,-10)*{ \bullet};
(-8,-10)*{ \bullet};
(10,-10)*{ \bullet};
(-10,-10)*{ \bullet};
(-3,-11)*{ \bullet};
(3,-11)*{ \bullet};
(-1,-11)*{ \bullet};
(1,-11)*{ \bullet};
(-5,-11)*{ \bullet};
(5,-11)*{ \bullet};
(-7,-11)*{ \bullet};
(7,-11)*{ \bullet};
(-9,-11)*{ \bullet};
(9,-11)*{ \bullet};
(-11,-11)*{ \bullet};
(11,-11)*{ \bullet};
(4,-12)*{ \bullet};
(-4,-12)*{ \bullet};
(2,-12)*{ \bullet};
(-2,-12)*{ \bullet};
(0,-12)*{ \bullet};
(6,-12)*{ \bullet};
(-6,-12)*{ \bullet};
(8,-12)*{ \bullet};
(-8,-12)*{ \bullet};
(10,-12)*{ \bullet};
(-10,-12)*{ \bullet};
(12,-12)*{ \bullet};
(-12,-12)*{ \bullet};
(-3,-13)*{ \bullet};
(3,-13)*{ \bullet};
(-1,-13)*{ \bullet};
(1,-13)*{ \bullet};
(-5,-13)*{ \bullet};
(5,-13)*{ \bullet};
(-7,-13)*{ \bullet};
(7,-13)*{ \bullet};
(-9,-13)*{ \bullet};
(9,-13)*{ \bullet};
(-11,-13)*{ \bullet};
(11,-13)*{ \bullet};
(-13,-13)*{ \bullet};
(13,-13)*{ \bullet};
(-12,0); (-12,-14.5) **\dir{.};
(-7,0); (-7,-14.5) **\dir{.};
(-2,0); (-2,-14.5) **\dir{.};
(3,0); (3,-14.5) **\dir{.};
(8,0); (8,-14.5) **\dir{.};
(13,0); (13,-14.5) **\dir{.};
{\ar@{->} (-14,-15)*{};(14,-15)*{}};
(-12,-15.5)*{\scriptscriptstyle -12};
(-7,-15.5)*{\scriptscriptstyle -7};
(-2,-15.5)*{\scriptscriptstyle -2};
(3,-15.5)*{\scriptscriptstyle 3};
(8,-15.5)*{\scriptscriptstyle 8};
(13,-15.5)*{\scriptscriptstyle 13};
(0,-14.5)*{\scriptscriptstyle  \text{weight}};
\endxy
}

\def\pascalbeameruno{\xy 0;/r.7pc/:
(0,0)*{ \bullet};
(-1,-1)*{ \bullet};
(1,-1)*{ \bullet};
(-2,-2)*{ \bullet};
(2,-2)*{ \bullet};
(0,-2)*{ \bullet};
(-3,-3)*{ \bullet};
(3,-3)*{ \bullet};
(-1,-3)*{ \bullet};
(1,-3)*{ \bullet};
(4,-4)*{ \bullet};
(-4,-4)*{ \bullet};
(2,-4)*{ \bullet};
(-2,-4)*{ \bullet};
(0,-4)*{ \bullet};
(-3,-5)*{ \bullet};
(3,-5)*{ \bullet};
(-1,-5)*{ \bullet};
(1,-5)*{ \bullet};
(-5,-5)*{ \bullet};
(5,-5)*{ \bullet};
(4,-6)*{ \bullet};
(-4,-6)*{ \bullet};
(2,-6)*{ \bullet};
(-2,-6)*{ \bullet};
(0,-6)*{ \bullet};
(6,-6)*{ \bullet};
(-6,-6)*{ \bullet};
(-3,-7)*{ \bullet};
(3,-7)*{ \bullet};
(-1,-7)*{ \bullet};
(1,-7)*{ \bullet};
(-5,-7)*{ \bullet};
(5,-7)*{ \bullet};
(-7,-7)*{ \bullet};
(7,-7)*{ \bullet};
(4,-8)*{ \bullet};
(-4,-8)*{ \bullet};
(2,-8)*{ \bullet};
(-2,-8)*{ \bullet};
(0,-8)*{ \bullet};
(6,-8)*{ \bullet};
(-6,-8)*{ \bullet};
(8,-8)*{ \bullet};
(-8,-8)*{ \bullet};
(-3,-9)*{ \bullet};
(3,-9)*{ \bullet};
(-1,-9)*{ \bullet};
(1,-9)*{ \bullet};
(-5,-9)*{ \bullet};
(5,-9)*{ \bullet};
(-7,-9)*{ \bullet};
(7,-9)*{ \bullet};
(-9,-9)*{ \bullet};
(9,-9)*{ \bullet};
(4,-10)*{ \bullet};
(-4,-10)*{ \bullet};
(2,-10)*{ \bullet};
(-2,-10)*{ \bullet};
(0,-10)*{ \bullet};
(6,-10)*{ \bullet};
(-6,-10)*{ \bullet};
(8,-10)*{ \bullet};
(-8,-10)*{ \bullet};
(10,-10)*{ \bullet};
(-10,-10)*{ \bullet};
(-3,-11)*{ \bullet};
(3,-11)*{ \bullet};
(-1,-11)*{ \bullet};
(1,-11)*{ \bullet};
(-5,-11)*{ \bullet};
(5,-11)*{ \bullet};
(-7,-11)*{ \bullet};
(7,-11)*{ \bullet};
(-9,-11)*{ \bullet};
(9,-11)*{ \bullet};
(-11,-11)*{ \bullet};
(11,-11)*{ \bullet};
(4,-12)*{ \bullet};
(-4,-12)*{ \bullet};
(2,-12)*{ \bullet};
(-2,-12)*{ \bullet};
(0,-12)*{ \bullet};
(6,-12)*{ \bullet};
(-6,-12)*{ \bullet};
(8,-12)*{ \bullet};
(-8,-12)*{ \bullet};
(10,-12)*{ \bullet};
(-10,-12)*{ \bullet};
(12,-12)*{ \bullet};
(-12,-12)*{ \bullet};
(-3,-13)*{ \bullet};
(3,-13)*{ \bullet};
(-1,-13)*{ \bullet};
(1,-13)*{ \bullet};
(-5,-13)*{ \bullet};
(5,-13)*{ \bullet};
(-7,-13)*{ \bullet};
(7,-13)*{ \bullet};
(-9,-13)*{ \bullet};
(9,-13)*{ \bullet};
(-11,-13)*{ \bullet};
(11,-13)*{ \bullet};
(-13,-13)*{ \bullet};
(13,-13)*{ \bullet};
(14,-14)*{ \bullet};
(-14,-14)*{ \bullet};
(4,-14)*{ \bullet};
(-4,-14)*{ \bullet};
(2,-14)*{ \bullet};
(-2,-14)*{ \bullet};
(0,-14)*{ \bullet};
(6,-14)*{ \bullet};
(-6,-14)*{ \bullet};
(8,-14)*{ \bullet};
(-8,-14)*{ \bullet};
(10,-14)*{ \bullet};
(-10,-14)*{ \bullet};
(12,-14)*{ \bullet};
(-12,-14)*{ \bullet};
(-15,-15)*{ \bullet};
(15,-15)*{ \bullet};
(-3,-15)*{ \bullet};
(3,-15)*{ \bullet};
(-1,-15)*{ \bullet};
(1,-15)*{ \bullet};
(-5,-15)*{ \bullet};
(5,-15)*{ \bullet};
(-7,-15)*{ \bullet};
(7,-15)*{ \bullet};
(-9,-15)*{ \bullet};
(9,-15)*{ \bullet};
(-11,-15)*{ \bullet};
(11,-15)*{ \bullet};
(-13,-15)*{ \bullet};
(13,-15)*{ \bullet};
{\ar@{<->} (-15,-17)*{};(15,-17)*{}};
{\ar@{->} (-18,0)*{};(-18,-15)*{}};
(0,-18)*{\scriptscriptstyle WEIGHT};
(-18,1)*{\scriptscriptstyle LEVEL};
\endxy
}

\def\pascalbeamerdos{\xy 0;/r.7pc/:
(0,0)*{ \bullet};
(-1,-1)*{ \bullet};
(1,-1)*{ \bullet};
(-2,-2)*{ \bullet};
(2,-2)*{ \bullet};
(0,-2)*{ \bullet};
(-3,-3)*{ \bullet};
(3,-3)*{ \bullet};
(-1,-3)*{ \bullet};
(1,-3)*{ \bullet};
(4,-4)*{ \bullet};
(-4,-4)*{ \bullet};
(2,-4)*{ \bullet};
(-2,-4)*{ \bullet};
(0,-4)*{ \bullet};
(-3,-5)*{ \bullet};
(3,-5)*{ \bullet};
(-1,-5)*{ \bullet};
(1,-5)*{ \bullet};
(-5,-5)*{ \bullet};
(5,-5)*{ \bullet};
(4,-6)*{ \bullet};
(-4,-6)*{ \bullet};
(2,-6)*{ \bullet};
(-2,-6)*{ \bullet};
(0,-6)*{ \bullet};
(6,-6)*{ \bullet};
(-6,-6)*{ \bullet};
(-3,-7)*{ \bullet};
(3,-7)*{ \bullet};
(-1,-7)*{ \bullet};
(1,-7)*{ \bullet};
(-5,-7)*{ \bullet};
(5,-7)*{ \bullet};
(-7,-7)*{ \bullet};
(7,-7)*{ \bullet};
(4,-8)*{ \bullet};
(-4,-8)*{ \bullet};
(2,-8)*{ \bullet};
(-2,-8)*{ \bullet};
(0,-8)*{ \bullet};
(6,-8)*{ \bullet};
(-6,-8)*{ \bullet};
(8,-8)*{ \bullet};
(-8,-8)*{ \bullet};
(-3,-9)*{ \bullet};
(3,-9)*{ \bullet};
(-1,-9)*{ \bullet};
(1,-9)*{ \bullet};
(-5,-9)*{ \bullet};
(5,-9)*{ \bullet};
(-7,-9)*{ \bullet};
(7,-9)*{ \bullet};
(-9,-9)*{ \bullet};
(9,-9)*{ \bullet};
(4,-10)*{ \bullet};
(-4,-10)*{ \bullet};
(2,-10)*{ \bullet};
(-2,-10)*{ \bullet};
(0,-10)*{ \bullet};
(6,-10)*{ \bullet};
(-6,-10)*{ \bullet};
(8,-10)*{ \bullet};
(-8,-10)*{ \bullet};
(10,-10)*{ \bullet};
(-10,-10)*{ \bullet};
(-3,-11)*{ \bullet};
(3,-11)*{ \bullet};
(-1,-11)*{ \bullet};
(1,-11)*{ \bullet};
(-5,-11)*{ \bullet};
(5,-11)*{ \bullet};
(-7,-11)*{ \bullet};
(7,-11)*{ \bullet};
(-9,-11)*{ \bullet};
(9,-11)*{ \bullet};
(-11,-11)*{ \bullet};
(11,-11)*{ \bullet};
(4,-12)*{ \bullet};
(-4,-12)*{ \bullet};
(2,-12)*{ \bullet};
(-2,-12)*{ \bullet};
(0,-12)*{ \bullet};
(6,-12)*{ \bullet};
(-6,-12)*{ \bullet};
(8,-12)*{ \bullet};
(-8,-12)*{ \bullet};
(10,-12)*{ \bullet};
(-10,-12)*{ \bullet};
(12,-12)*{ \bullet};
(-12,-12)*{ \bullet};
(-3,-13)*{ \bullet};
(3,-13)*{ \bullet};
(-1,-13)*{ \bullet};
(1,-13)*{ \bullet};
(-5,-13)*{ \bullet};
(5,-13)*{ \bullet};
(-7,-13)*{ \bullet};
(7,-13)*{ \bullet};
(-9,-13)*{ \bullet};
(9,-13)*{ \bullet};
(-11,-13)*{ \bullet};
(11,-13)*{ \bullet};
(-13,-13)*{ \bullet};
(13,-13)*{ \bullet};
(14,-14)*{ \bullet};
(-14,-14)*{ \bullet};
(4,-14)*{ \bullet};
(-4,-14)*{ \bullet};
(2,-14)*{ \bullet};
(-2,-14)*{ \bullet};
(0,-14)*{ \bullet};
(6,-14)*{ \bullet};
(-6,-14)*{ \bullet};
(8,-14)*{ \bullet};
(-8,-14)*{ \bullet};
(10,-14)*{ \bullet};
(-10,-14)*{ \bullet};
(12,-14)*{ \bullet};
(-12,-14)*{ \bullet};
(-15,-15)*{ \bullet};
(15,-15)*{ \bullet};
(-3,-15)*{ \bullet};
(3,-15)*{ \bullet};
(-1,-15)*{ \bullet};
(1,-15)*{ \bullet};
(-5,-15)*{ \circ};
(5,-15)*{ \bullet};
(-7,-15)*{ \bullet};
(7,-15)*{ \bullet};
(-9,-15)*{ \bullet};
(9,-15)*{ \bullet};
(-11,-15)*{ \bullet};
(11,-15)*{ \bullet};
(-13,-15)*{ \bullet};
(13,-15)*{ \bullet};
{\ar@{<->} (-15,-17)*{};(15,-17)*{}};
{\ar@{->} (-18,0)*{};(-18,-15)*{}};
(0,-18)*{\scriptscriptstyle WEIGHT};
(-18,1)*{\scriptscriptstyle LEVEL};
\endxy
}

\def\pascalbeamertres{\xy 0;/r.7pc/:
{(0,0); (1,-1); \connect[!\linew(1.5)\rgbcolor(1 0 0)]  \dir{-}};
{(-1,-3); (1,-1); \connect[!\linew(1.5)\rgbcolor(1 0 0)]  \dir{-}};
{(-1,-3); (1,-5); \connect[!\linew(1.5)\rgbcolor(1 0 0)]  \dir{-}};
{(-1,-7); (1,-5); \connect[!\linew(1.5)\rgbcolor(1 0 0)]  \dir{-}};
{(-1,-7); (0,-8); \connect[!\linew(1.5)\rgbcolor(1 0 0)]  \dir{-}};
{(-2,-10); (0,-8); \connect[!\linew(1.5)\rgbcolor(1 0 0)]  \dir{-}};
{(-2,-10); (-1,-11); \connect[!\linew(1.5)\rgbcolor(1 0 0)]  \dir{-}};
{(-5,-15); (-1,-11); \connect[!\linew(1.5)\rgbcolor(1 0 0)]  \dir{-}};
(0,0)*{ \bullet};
(-1,-1)*{ \bullet};
(1,-1)*{ \bullet};
(-2,-2)*{ \bullet};
(2,-2)*{ \bullet};
(0,-2)*{ \bullet};
(-3,-3)*{ \bullet};
(3,-3)*{ \bullet};
(-1,-3)*{ \bullet};
(1,-3)*{ \bullet};
(4,-4)*{ \bullet};
(-4,-4)*{ \bullet};
(2,-4)*{ \bullet};
(-2,-4)*{ \bullet};
(0,-4)*{ \bullet};
(-3,-5)*{ \bullet};
(3,-5)*{ \bullet};
(-1,-5)*{ \bullet};
(1,-5)*{ \bullet};
(-5,-5)*{ \bullet};
(5,-5)*{ \bullet};
(4,-6)*{ \bullet};
(-4,-6)*{ \bullet};
(2,-6)*{ \bullet};
(-2,-6)*{ \bullet};
(0,-6)*{ \bullet};
(6,-6)*{ \bullet};
(-6,-6)*{ \bullet};
(-3,-7)*{ \bullet};
(3,-7)*{ \bullet};
(-1,-7)*{ \bullet};
(1,-7)*{ \bullet};
(-5,-7)*{ \bullet};
(5,-7)*{ \bullet};
(-7,-7)*{ \bullet};
(7,-7)*{ \bullet};
(4,-8)*{ \bullet};
(-4,-8)*{ \bullet};
(2,-8)*{ \bullet};
(-2,-8)*{ \bullet};
(0,-8)*{ \bullet};
(6,-8)*{ \bullet};
(-6,-8)*{ \bullet};
(8,-8)*{ \bullet};
(-8,-8)*{ \bullet};
(-3,-9)*{ \bullet};
(3,-9)*{ \bullet};
(-1,-9)*{ \bullet};
(1,-9)*{ \bullet};
(-5,-9)*{ \bullet};
(5,-9)*{ \bullet};
(-7,-9)*{ \bullet};
(7,-9)*{ \bullet};
(-9,-9)*{ \bullet};
(9,-9)*{ \bullet};
(4,-10)*{ \bullet};
(-4,-10)*{ \bullet};
(2,-10)*{ \bullet};
(-2,-10)*{ \bullet};
(0,-10)*{ \bullet};
(6,-10)*{ \bullet};
(-6,-10)*{ \bullet};
(8,-10)*{ \bullet};
(-8,-10)*{ \bullet};
(10,-10)*{ \bullet};
(-10,-10)*{ \bullet};
(-3,-11)*{ \bullet};
(3,-11)*{ \bullet};
(-1,-11)*{ \bullet};
(1,-11)*{ \bullet};
(-5,-11)*{ \bullet};
(5,-11)*{ \bullet};
(-7,-11)*{ \bullet};
(7,-11)*{ \bullet};
(-9,-11)*{ \bullet};
(9,-11)*{ \bullet};
(-11,-11)*{ \bullet};
(11,-11)*{ \bullet};
(4,-12)*{ \bullet};
(-4,-12)*{ \bullet};
(2,-12)*{ \bullet};
(-2,-12)*{ \bullet};
(0,-12)*{ \bullet};
(6,-12)*{ \bullet};
(-6,-12)*{ \bullet};
(8,-12)*{ \bullet};
(-8,-12)*{ \bullet};
(10,-12)*{ \bullet};
(-10,-12)*{ \bullet};
(12,-12)*{ \bullet};
(-12,-12)*{ \bullet};
(-3,-13)*{ \bullet};
(3,-13)*{ \bullet};
(-1,-13)*{ \bullet};
(1,-13)*{ \bullet};
(-5,-13)*{ \bullet};
(5,-13)*{ \bullet};
(-7,-13)*{ \bullet};
(7,-13)*{ \bullet};
(-9,-13)*{ \bullet};
(9,-13)*{ \bullet};
(-11,-13)*{ \bullet};
(11,-13)*{ \bullet};
(-13,-13)*{ \bullet};
(13,-13)*{ \bullet};
(14,-14)*{ \bullet};
(-14,-14)*{ \bullet};
(4,-14)*{ \bullet};
(-4,-14)*{ \bullet};
(2,-14)*{ \bullet};
(-2,-14)*{ \bullet};
(0,-14)*{ \bullet};
(6,-14)*{ \bullet};
(-6,-14)*{ \bullet};
(8,-14)*{ \bullet};
(-8,-14)*{ \bullet};
(10,-14)*{ \bullet};
(-10,-14)*{ \bullet};
(12,-14)*{ \bullet};
(-12,-14)*{ \bullet};
(-15,-15)*{ \bullet};
(15,-15)*{ \bullet};
(-3,-15)*{ \bullet};
(3,-15)*{ \bullet};
(-1,-15)*{ \bullet};
(1,-15)*{ \bullet};
(-5,-15)*{ \bullet};
(5,-15)*{ \bullet};
(-7,-15)*{ \bullet};
(7,-15)*{ \bullet};
(-9,-15)*{ \bullet};
(9,-15)*{ \bullet};
(-11,-15)*{ \bullet};
(11,-15)*{ \bullet};
(-13,-15)*{ \bullet};
(13,-15)*{ \bullet};
{\ar@{<->} (-15,-17)*{};(15,-17)*{}};
{\ar@{->} (-18,0)*{};(-18,-15)*{}};
(0,-18)*{\scriptscriptstyle WEIGHT};
(-18,1)*{\scriptscriptstyle LEVEL};
\endxy
}

\def\pascalbeamertresmayor{\xy 0;/r.7pc/:
{(0,0); (-1,-1); \connect[!\linew(1.5)\rgbcolor(1 0 0)]  \dir{-}};
{(0,-2); (-1,-1); \connect[!\linew(1.5)\rgbcolor(1 0 0)]  \dir{-}};
{(0,-2); (-1,-3); \connect[!\linew(1.5)\rgbcolor(1 0 0)]  \dir{-}};
{(0,-4); (-1,-3); \connect[!\linew(1.5)\rgbcolor(1 0 0)]  \dir{-}};
{(0,-4); (-1,-5); \connect[!\linew(1.5)\rgbcolor(1 0 0)]  \dir{-}};
{(0,-6); (-1,-5); \connect[!\linew(1.5)\rgbcolor(1 0 0)]  \dir{-}};
{(0,-6); (-9,-15); \connect[!\linew(1.5)\rgbcolor(1 0 0)]  \dir{-}};
(0,0)*{ \bullet};
(-1,-1)*{ \bullet};
(1,-1)*{ \bullet};
(-2,-2)*{ \bullet};
(2,-2)*{ \bullet};
(0,-2)*{ \bullet};
(-3,-3)*{ \bullet};
(3,-3)*{ \bullet};
(-1,-3)*{ \bullet};
(1,-3)*{ \bullet};
(4,-4)*{ \bullet};
(-4,-4)*{ \bullet};
(2,-4)*{ \bullet};
(-2,-4)*{ \bullet};
(0,-4)*{ \bullet};
(-3,-5)*{ \bullet};
(3,-5)*{ \bullet};
(-1,-5)*{ \bullet};
(1,-5)*{ \bullet};
(-5,-5)*{ \bullet};
(5,-5)*{ \bullet};
(4,-6)*{ \bullet};
(-4,-6)*{ \bullet};
(2,-6)*{ \bullet};
(-2,-6)*{ \bullet};
(0,-6)*{ \bullet};
(6,-6)*{ \bullet};
(-6,-6)*{ \bullet};
(-3,-7)*{ \bullet};
(3,-7)*{ \bullet};
(-1,-7)*{ \bullet};
(1,-7)*{ \bullet};
(-5,-7)*{ \bullet};
(5,-7)*{ \bullet};
(-7,-7)*{ \bullet};
(7,-7)*{ \bullet};
(4,-8)*{ \bullet};
(-4,-8)*{ \bullet};
(2,-8)*{ \bullet};
(-2,-8)*{ \bullet};
(0,-8)*{ \bullet};
(6,-8)*{ \bullet};
(-6,-8)*{ \bullet};
(8,-8)*{ \bullet};
(-8,-8)*{ \bullet};
(-3,-9)*{ \bullet};
(3,-9)*{ \bullet};
(-1,-9)*{ \bullet};
(1,-9)*{ \bullet};
(-5,-9)*{ \bullet};
(5,-9)*{ \bullet};
(-7,-9)*{ \bullet};
(7,-9)*{ \bullet};
(-9,-9)*{ \bullet};
(9,-9)*{ \bullet};
(4,-10)*{ \bullet};
(-4,-10)*{ \bullet};
(2,-10)*{ \bullet};
(-2,-10)*{ \bullet};
(0,-10)*{ \bullet};
(6,-10)*{ \bullet};
(-6,-10)*{ \bullet};
(8,-10)*{ \bullet};
(-8,-10)*{ \bullet};
(10,-10)*{ \bullet};
(-10,-10)*{ \bullet};
(-3,-11)*{ \bullet};
(3,-11)*{ \bullet};
(-1,-11)*{ \bullet};
(1,-11)*{ \bullet};
(-5,-11)*{ \bullet};
(5,-11)*{ \bullet};
(-7,-11)*{ \bullet};
(7,-11)*{ \bullet};
(-9,-11)*{ \bullet};
(9,-11)*{ \bullet};
(-11,-11)*{ \bullet};
(11,-11)*{ \bullet};
(4,-12)*{ \bullet};
(-4,-12)*{ \bullet};
(2,-12)*{ \bullet};
(-2,-12)*{ \bullet};
(0,-12)*{ \bullet};
(6,-12)*{ \bullet};
(-6,-12)*{ \bullet};
(8,-12)*{ \bullet};
(-8,-12)*{ \bullet};
(10,-12)*{ \bullet};
(-10,-12)*{ \bullet};
(12,-12)*{ \bullet};
(-12,-12)*{ \bullet};
(-3,-13)*{ \bullet};
(3,-13)*{ \bullet};
(-1,-13)*{ \bullet};
(1,-13)*{ \bullet};
(-5,-13)*{ \bullet};
(5,-13)*{ \bullet};
(-7,-13)*{ \bullet};
(7,-13)*{ \bullet};
(-9,-13)*{ \bullet};
(9,-13)*{ \bullet};
(-11,-13)*{ \bullet};
(11,-13)*{ \bullet};
(-13,-13)*{ \bullet};
(13,-13)*{ \bullet};
(14,-14)*{ \bullet};
(-14,-14)*{ \bullet};
(4,-14)*{ \bullet};
(-4,-14)*{ \bullet};
(2,-14)*{ \bullet};
(-2,-14)*{ \bullet};
(0,-14)*{ \bullet};
(6,-14)*{ \bullet};
(-6,-14)*{ \bullet};
(8,-14)*{ \bullet};
(-8,-14)*{ \bullet};
(10,-14)*{ \bullet};
(-10,-14)*{ \bullet};
(12,-14)*{ \bullet};
(-12,-14)*{ \bullet};
(-15,-15)*{ \bullet};
(15,-15)*{ \bullet};
(-3,-15)*{ \bullet};
(3,-15)*{ \bullet};
(-1,-15)*{ \bullet};
(1,-15)*{ \bullet};
(-5,-15)*{ \bullet};
(5,-15)*{ \bullet};
(-7,-15)*{ \bullet};
(7,-15)*{ \bullet};
(-9,-15)*{ \bullet};
(9,-15)*{ \bullet};
(-11,-15)*{ \bullet};
(11,-15)*{ \bullet};
(-13,-15)*{ \bullet};
(13,-15)*{ \bullet};
{\ar@{<->} (-15,-17)*{};(15,-17)*{}};
{\ar@{->} (-18,0)*{};(-18,-15)*{}};
(0,-18)*{\scriptscriptstyle WEIGHT};
(-18,1)*{\scriptscriptstyle LEVEL};
\endxy
}

\def\pascalbeamer{\xy 0;/r.7pc/:
(0,0)*{ \bullet};
(-1,-1)*{ \bullet};
(1,-1)*{ \bullet};
(-2,-2)*{ \bullet};
(2,-2)*{ \bullet};
(0,-2)*{ \bullet};
(-3,-3)*{ \bullet};
(3,-3)*{ \bullet};
(-1,-3)*{ \bullet};
(1,-3)*{ \bullet};
(4,-4)*{ \bullet};
(-4,-4)*{ \bullet};
(2,-4)*{ \bullet};
(-2,-4)*{ \bullet};
(0,-4)*{ \bullet};
(-3,-5)*{ \bullet};
(3,-5)*{ \bullet};
(-1,-5)*{ \bullet};
(1,-5)*{ \bullet};
(-5,-5)*{ \bullet};
(5,-5)*{ \bullet};
(4,-6)*{ \bullet};
(-4,-6)*{ \bullet};
(2,-6)*{ \bullet};
(-2,-6)*{ \bullet};
(0,-6)*{ \bullet};
(6,-6)*{ \bullet};
(-6,-6)*{ \bullet};
(-3,-7)*{ \bullet};
(3,-7)*{ \bullet};
(-1,-7)*{ \bullet};
(1,-7)*{ \bullet};
(-5,-7)*{ \bullet};
(5,-7)*{ \bullet};
(-7,-7)*{ \bullet};
(7,-7)*{ \bullet};
(4,-8)*{ \bullet};
(-4,-8)*{ \bullet};
(2,-8)*{ \bullet};
(-2,-8)*{ \bullet};
(0,-8)*{ \bullet};
(6,-8)*{ \bullet};
(-6,-8)*{ \bullet};
(8,-8)*{ \bullet};
(-8,-8)*{ \bullet};
(-3,-9)*{ \bullet};
(3,-9)*{ \bullet};
(-1,-9)*{ \bullet};
(1,-9)*{ \bullet};
(-5,-9)*{ \bullet};
(5,-9)*{ \bullet};
(-7,-9)*{ \bullet};
(7,-9)*{ \bullet};
(-9,-9)*{ \bullet};
(9,-9)*{ \bullet};
(4,-10)*{ \bullet};
(-4,-10)*{ \bullet};
(2,-10)*{ \bullet};
(-2,-10)*{ \bullet};
(0,-10)*{ \bullet};
(6,-10)*{ \bullet};
(-6,-10)*{ \bullet};
(8,-10)*{ \bullet};
(-8,-10)*{ \bullet};
(10,-10)*{ \bullet};
(-10,-10)*{ \bullet};
(-3,-11)*{ \bullet};
(3,-11)*{ \bullet};
(-1,-11)*{ \bullet};
(1,-11)*{ \bullet};
(-5,-11)*{ \bullet};
(5,-11)*{ \bullet};
(-7,-11)*{ \bullet};
(7,-11)*{ \bullet};
(-9,-11)*{ \bullet};
(9,-11)*{ \bullet};
(-11,-11)*{ \bullet};
(11,-11)*{ \bullet};
(4,-12)*{ \bullet};
(-4,-12)*{ \bullet};
(2,-12)*{ \bullet};
(-2,-12)*{ \bullet};
(0,-12)*{ \bullet};
(6,-12)*{ \bullet};
(-6,-12)*{ \bullet};
(8,-12)*{ \bullet};
(-8,-12)*{ \bullet};
(10,-12)*{ \bullet};
(-10,-12)*{ \bullet};
(12,-12)*{ \bullet};
(-12,-12)*{ \bullet};
(-3,-13)*{ \bullet};
(3,-13)*{ \bullet};
(-1,-13)*{ \bullet};
(1,-13)*{ \bullet};
(-5,-13)*{ \bullet};
(5,-13)*{ \bullet};
(-7,-13)*{ \bullet};
(7,-13)*{ \bullet};
(-9,-13)*{ \bullet};
(9,-13)*{ \bullet};
(-11,-13)*{ \bullet};
(11,-13)*{ \bullet};
(-13,-13)*{ \bullet};
(13,-13)*{ \bullet};
(14,-14)*{ \bullet};
(-14,-14)*{ \bullet};
(4,-14)*{ \bullet};
(-4,-14)*{ \bullet};
(2,-14)*{ \bullet};
(-2,-14)*{ \bullet};
(0,-14)*{ \bullet};
(6,-14)*{ \bullet};
(-6,-14)*{ \bullet};
(8,-14)*{ \bullet};
(-8,-14)*{ \bullet};
(10,-14)*{ \bullet};
(-10,-14)*{ \bullet};
(12,-14)*{ \bullet};
(-12,-14)*{ \bullet};
(-15,-15)*{ \bullet};
(15,-15)*{ \bullet};
(-3,-13)*{ \bullet};
(3,-15)*{ \bullet};
(-1,-15)*{ \bullet};
(1,-15)*{ \bullet};
(-5,-15)*{ \bullet};
(5,-15)*{ \bullet};
(-7,-15)*{ \bullet};
(7,-15)*{ \bullet};
(-9,-15)*{ \bullet};
(9,-15)*{ \bullet};
(-11,-15)*{ \bullet};
(11,-15)*{ \bullet};
(-13,-15)*{ \bullet};
(13,-15)*{ \bullet};
(-12,0); (-12,-14.5) **\dir{-};
(-7,0); (-7,-14.5) **\dir{-};
(-2,0); (-2,-14.5) **\dir{-};
(3,0); (3,-14.5) **\dir{-};
(8,0); (8,-14.5) **\dir{-};
(13,0); (13,-14.5) **\dir{-};
{\ar@{->} (-14,-15)*{};(14,-15)*{}};
(-12,-15.5)*{\scriptscriptstyle -12};
(-7,-15.5)*{\scriptscriptstyle -7};
(-2,-15.5)*{\scriptscriptstyle -2};
(3,-15.5)*{\scriptscriptstyle 3};
(8,-15.5)*{\scriptscriptstyle 8};
(13,-15.5)*{\scriptscriptstyle 13};
(0,-14.5)*{\scriptscriptstyle  \text{weight}};
\endxy
}

\def\maximales{\xy 0;/r.8pc/:
(0,0); (-13,-13) **\dir{-};
(-11,-13); (-12,-12) **\dir{-};
(-7,-7); (-2,-12) **\dir{-};
(-3,-13); (-2,-12) **\dir{-};
(-1,-13); (-2,-12) **\dir{-};
(-2,-2); (9,-13) **\dir{-};
(3,-7); (-2,-12) **\dir{-};
(7,-13); (8,-12) **\dir{-};
(0,0)*{ \bullet};
(-1,-1)*{ \bullet};
(1,-1)*{ \bullet};
(-2,-2)*{ \bullet};
(2,-2)*{ \bullet};
(0,-2)*{ \bullet};
(-3,-3)*{ \bullet};
(3,-3)*{ \bullet};
(-1,-3)*{ \bullet};
(1,-3)*{ \bullet};
(4,-4)*{ \bullet};
(-4,-4)*{ \bullet};
(2,-4)*{ \bullet};
(-2,-4)*{ \bullet};
(0,-4)*{ \bullet};
(-3,-5)*{ \bullet};
(3,-5)*{ \bullet};
(-1,-5)*{ \bullet};
(1,-5)*{ \bullet};
(-5,-5)*{ \bullet};
(5,-5)*{ \bullet};
(4,-6)*{ \bullet};
(-4,-6)*{ \bullet};
(2,-6)*{ \bullet};
(-2,-6)*{ \bullet};
(0,-6)*{ \bullet};
(6,-6)*{ \bullet};
(-6,-6)*{ \bullet};
(-3,-7)*{ \bullet};
(3,-7)*{ \bullet};
(-1,-7)*{ \bullet};
(1,-7)*{ \bullet};
(-5,-7)*{ \bullet};
(5,-7)*{ \bullet};
(-7,-7)*{ \bullet};
(7,-7)*{ \bullet};
(4,-8)*{ \bullet};
(-4,-8)*{ \bullet};
(2,-8)*{ \bullet};
(-2,-8)*{ \bullet};
(0,-8)*{ \bullet};
(6,-8)*{ \bullet};
(-6,-8)*{ \bullet};
(8,-8)*{ \bullet};
(-8,-8)*{ \bullet};
(-3,-9)*{ \bullet};
(3,-9)*{ \bullet};
(-1,-9)*{ \bullet};
(1,-9)*{ \bullet};
(-5,-9)*{ \bullet};
(5,-9)*{ \bullet};
(-7,-9)*{ \bullet};
(7,-9)*{ \bullet};
(-9,-9)*{ \bullet};
(9,-9)*{ \bullet};
(4,-10)*{ \bullet};
(-4,-10)*{ \bullet};
(2,-10)*{ \bullet};
(-2,-10)*{ \bullet};
(0,-10)*{ \bullet};
(6,-10)*{ \bullet};
(-6,-10)*{ \bullet};
(8,-10)*{ \bullet};
(-8,-10)*{ \bullet};
(10,-10)*{ \bullet};
(-10,-10)*{ \bullet};
(-3,-11)*{ \bullet};
(3,-11)*{ \bullet};
(-1,-11)*{ \bullet};
(1,-11)*{ \bullet};
(-5,-11)*{ \bullet};
(5,-11)*{ \bullet};
(-7,-11)*{ \bullet};
(7,-11)*{ \bullet};
(-9,-11)*{ \bullet};
(9,-11)*{ \bullet};
(-11,-11)*{ \bullet};
(11,-11)*{ \bullet};
(4,-12)*{ \bullet};
(-4,-12)*{ \bullet};
(2,-12)*{ \bullet};
(-2,-12)*{ \bullet};
(0,-12)*{ \bullet};
(6,-12)*{ \bullet};
(-6,-12)*{ \bullet};
(8,-12)*{ \bullet};
(-8,-12)*{ \bullet};
(10,-12)*{ \bullet};
(-10,-12)*{ \bullet};
(12,-12)*{ \bullet};
(-12,-12)*{ \bullet};
(-3,-13)*{ \bullet};
(3,-13)*{ \bullet};
(-1,-13)*{ \bullet};
(1,-13)*{ \bullet};
(-5,-13)*{ \bullet};
(5,-13)*{ \bullet};
(-7,-13)*{ \bullet};
(7,-13)*{ \bullet};
(-9,-13)*{ \bullet};
(9,-13)*{ \bullet};
(-11,-13)*{ \bullet};
(11,-13)*{ \bullet};
(-13,-13)*{ \bullet};
(13,-13)*{ \bullet};
(-12,0); (-12,-14.5) **\dir{.};
(-7,0); (-7,-14.5) **\dir{.};
(-2,0); (-2,-14.5) **\dir{.};
(3,0); (3,-14.5) **\dir{.};
(8,0); (8,-14.5) **\dir{.};
(13,0); (13,-14.5) **\dir{.};
{\ar@{->} (-14,-15)*{};(14,-15)*{}};
(-12,-15.5)*{\scriptscriptstyle -12};
(-7,-15.5)*{\scriptscriptstyle -7};
(-2,-15.5)*{\scriptscriptstyle -2};
(3,-15.5)*{\scriptscriptstyle 3};
(8,-15.5)*{\scriptscriptstyle 8};
(13,-15.5)*{\scriptscriptstyle 13};
(0,-14.5)*{\scriptscriptstyle  \text{weight}};
\endxy
}

\def\finalrapidouno{\xy 0;/r.8pc/:
{\ar@{<->} (-16,0)*{};(16,0)*{}};
(0,.3); (0,-.3) **\dir{-};
(0,-1)*{0};
\endxy}

\def\finalrapidodos{\xy 0;/r.8pc/:
{\ar@{<->} (-16,0)*{};(16,0)*{}};
(-2,1); (-2,-1) **\dir{-};
(0,.3); (0,-.3) **\dir{-};
(-2,-2)*{-2};
(0,-1)*{0};
\endxy}

\def\finalrapidotres{\xy 0;/r.8pc/:
{\ar@{<->} (-16,0)*{};(16,0)*{}};
(3,1); (3,-1) **\dir{-};
(8,1); (8,-1) **\dir{-};
(13,1); (13,-1) **\dir{-};
(-2,1); (-2,-1) **\dir{-};
(-7,1); (-7,-1) **\dir{-};
(-12,1); (-12,-1) **\dir{-};
(-2,-2)*{-2};
(3,-2)*{3};
(8,-2)*{8};
(13,-2)*{13};
(-7,-2)*{-7};
(-12,-2)*{-12};
\endxy}

\def\circulos{\xy 0;/r.8pc/:
{\ar@{<->} (-10,0)*{};(-10,0)*{}};
{\ar@{<->} (0,-10)*{};(0,10)*{}};
(0,0)*{}*\xycircle<10pt>{};
(0,0)*{}*\xycircle<20pt>{};
(0,0)*{\bullet};
\endxy}

\def\bloblala{
\begin{pspicture}(0.5,0)(1.5,1)
\psarc(.75,1){.25}{180}{360}
\psarc(.75,0){.25}{360}{180}
\psline{-}(1.5,0)(1.5,1)
\rput[a](.75,.75){$ \bullet $}
\rput[a](.75,.25){$ \bullet $}
\rput[a](1.5,.5){$ \bullet $}
\end{pspicture}
}

\def\bloblat{
\begin{pspicture}(0.5,0)(1.5,1)
\psarc(.75,1){.25}{180}{360}
\psarc(.75,0){.25}{360}{180}
\psline{-}(1.5,0)(1.5,1)
\rput[a](.75,.75){$ \bullet $}
\rput[a](1.5,.5){$ \bullet $}
\end{pspicture}
}

\def\blobtla{
\begin{pspicture}(0.5,0)(1.5,1)
\psarc(.75,1){.25}{180}{360}
\psarc(.75,0){.25}{360}{180}
\psline{-}(1.5,0)(1.5,1)
\rput[a](.75,.25){$ \bullet $}
\rput[a](1.5,.5){$ \bullet $}
\end{pspicture}
}

\def\blobsla{
\begin{pspicture}(0.5,0)(1.5,1)
\psarc(1.25,1){.25}{180}{360}
\psarc(.75,0){.25}{360}{180}
\psline{-}(1.5,0)(.5,1)
\rput[a](.75,.25){$ \bullet $}
\rput[a](1,.5){$ \bullet $}
\end{pspicture}
}

\def\bloblas{
\begin{pspicture}(0.5,0)(1.5,1)
\psarc(.75,1){.25}{180}{360}
\psarc(1.25,0){.25}{360}{180}
\psline{-}(.5,0)(1.5,1)
\rput[a](.75,.75){$ \bullet $}
\rput[a](1,.5){$ \bullet $}
\end{pspicture}
}

\def\blobmumu{
\begin{pspicture}(0.5,0)(1.5,1)
\psarc(.75,1){.25}{180}{360}
\psarc(.75,0){.25}{360}{180}
\psline{-}(1.5,0)(1.5,1)
\rput[a](.75,.75){$ \bullet $}
\rput[a](.75,.25){$ \bullet $}
\end{pspicture}
}

\def\blobmuv{
\begin{pspicture}(0.5,0)(1.5,1)
\psarc(.75,1){.25}{180}{360}
\psarc(.75,0){.25}{360}{180}
\psline{-}(1.5,0)(1.5,1)
\rput[a](.75,.75){$ \bullet $}
\end{pspicture}
}

\def\blobvmu{
\begin{pspicture}(0.5,0)(1.5,1)
\psarc(.75,1){.25}{180}{360}
\psarc(.75,0){.25}{360}{180}
\psline{-}(1.5,0)(1.5,1)
\rput[a](.75,.25){$ \bullet $}
\end{pspicture}
}

\def\blobmuu{
\begin{pspicture}(0.5,0)(1.5,1)
\psarc(.75,1){.25}{180}{360}
\psarc(1.25,0){.25}{360}{180}
\psline{-}(.5,0)(1.5,1)
\rput[a](.75,.75){$ \bullet $}
\end{pspicture}
}

\def\blobumu{
\begin{pspicture}(0.5,0)(1.5,1)
\psarc(1.25,1){.25}{180}{360}
\psarc(.75,0){.25}{360}{180}
\psline{-}(1.5,0)(.5,1)
\rput[a](.75,.25){$ \bullet $}
\end{pspicture}
}

\def\blobnunu{
\begin{pspicture}(0.5,0)(1.5,1)
\psline{-}(.5,0)(.5,1)
\psline{-}(1.5,0)(1.5,1)
\psline{-}(1,0)(1,1)
\rput[a](.5,.5){$ \bullet $}
\end{pspicture}
}

\def\blobts{
\begin{pspicture}(0.5,0)(1.5,1)
\psarc(.75,1){.25}{180}{360}
\psarc(1.25,0){.25}{360}{180}
\psline{-}(.5,0)(1.5,1)
\rput[a](1,.5){$ \bullet $}
\end{pspicture}
}

\def\blobst{
\begin{pspicture}(0.5,0)(1.5,1)
\psarc(1.25,1){.25}{180}{360}
\psarc(.75,0){.25}{360}{180}
\psline{-}(1.5,0)(.5,1)
\rput[a](1,.5){$ \bullet $}
\end{pspicture}
}

\def\blobss{
\begin{pspicture}(0.5,0)(1.5,1)
\psline{-}(.5,0)(.5,1)
\psarc(1.25,1){.25}{180}{360}
\psarc(1.25,0){.25}{360}{180}
\rput[a](.5,.5){$ \bullet $}
\end{pspicture}
}

\def\blobtt{
\begin{pspicture}(0.5,0)(1.5,1)
\psarc(.75,1){.25}{180}{360}
\psarc(.75,0){.25}{360}{180}
\psline{-}(1.5,0)(1.5,1)
\rput[a](1.5,.5){$ \bullet $}
\end{pspicture}
}

\def\blobvu{
\begin{pspicture}(0.5,0)(1.5,1)
\psarc(.75,1){.25}{180}{360}
\psarc(1.25,0){.25}{360}{180}
\psline{-}(.5,0)(1.5,1)
\end{pspicture}
}

\def\blobuv{
\begin{pspicture}(0.5,0)(1.5,1)
\psarc(1.25,1){.25}{180}{360}
\psarc(.75,0){.25}{360}{180}
\psline{-}(1.5,0)(.5,1)
\end{pspicture}
}

\def\blobuu{
\begin{pspicture}(0.5,0)(1.5,1)
\psline{-}(.5,0)(.5,1)
\psarc(1.25,1){.25}{180}{360}
\psarc(1.25,0){.25}{360}{180}
\end{pspicture}
}

\def\blobvv{
\begin{pspicture}(0.5,0)(1.5,1)
\psarc(.75,1){.25}{180}{360}
\psarc(.75,0){.25}{360}{180}
\psline{-}(1.5,0)(1.5,1)
\end{pspicture}
}

\def\blobkk{
\begin{pspicture}(0.5,0)(1.5,1)
\psline{-}(.5,0)(.5,1)
\psline{-}(1.5,0)(1.5,1)
\psline{-}(1,0)(1,1)
\end{pspicture}
}

\newcommand{\K}{K}
\newcommand{\Shape}{\operatorname{Shape}}

\newcommand{\m}{{\mathfrak m}}
\newcommand{\OO}{{\mathcal O}}
\newtheorem{teo}{Theorem}[section]
\newtheorem{lem}[teo]{Lemma}

\newtheorem{defi}[teo]{Definition}
\newtheorem{cor}[teo]{Corollary}
\newtheorem{exa}[teo]{Example}
\newtheorem{rem}[teo]{Remark}
\newcommand{\C}{ \mathbb C}
\newcommand{\Z}{ \mathbb Z}
\newcommand{\blob}[1]{ {b_n(#1)}  }
\newcommand{\blobR}[1]{ {b_n^{\OO}(#1)}  }
\newcommand{\blobK}[1]{ {b_n^{\K}(#1)}  }
\newcommand{\F}{ {\mathbb C}}

\newcommand{\R}{R}
\newcommand{\heckedos}{{\mathcal{H}}_2(m) }
\newcommand{\hecke}{{\mathcal{H}}_n(m) }
\newcommand{\heckeintegral}{{\mathcal{H}}^{\OO}_n(m) }
\newcommand{\heckerational}{{\mathcal{H}}^{\K}_n(m) }
\newcommand{\heckenmasuno}{{\mathcal{H}}_{n+1}(m) }
\newenvironment{dem}{\noindent \textit{Proof:} }{\quad \hfill $\square$}

\newcommand{\tld}[1]{   \beta_{ {#1\tau}^{\mu}   } }

\newcommand{\sub}[1]{  {  \mathfrak{#1}}}
\newcommand{\pa}[1]{ {\text{Par}(#1)}  }
\newcommand{\pacol}[1]{ {\text{Par}_2(#1)}  }
\newcommand{\bi}[1]{ \operatorname{Bip}_{1}(#1) }
\newcommand{\bs}[1]{ \boldsymbol{#1}}
\newcommand{\jm}[1]{ {{L}}_{#1}}

\newcommand{\de}[1]{ \mathfrak{d} ({ \mathfrak{#1}})  }
\newcommand{\tc}[1]{ \tau_{\sub{ #1}} }
\newcommand{\Std}{\operatorname{Std}}
\newcommand{\Span}{\operatorname{span}}

\numberwithin{equation}{section}

\titleformat{\section}[hang]{\sc}{\thesection.}{0.5cm}{\filcenter}
\titleformat{\subsection}[hang]{\it}{\thesubsection.}{0.2cm}{\filright}

\begin{document}

\title{\bf \normalsize GRADED CELLULAR BASES FOR TEMPERLEY-LIEB ALGEBRAS OF TYPE A AND B.}
\author{\sc  david plaza{\thanks{ Supported in part by Beca Doctorado Nacional 2011-CONICYT}} and steen ryom-hansen\thanks{Supported in part by FONDECYT grants 1090701 and 1121129,
 by Programa Reticulados y Simetr\'ia
and by the MathAmSud project OPECSHA 01-math-10} }
\date{}   \maketitle
\begin{abstract}
\noindent \textsc{Abstract. }
We show that the Temperley-Lieb algebra of type $A$ and the blob algebra (also known as
the Temperley-Lieb algebra of type $ B$) at roots of
unity are $ \mathbb Z$-graded algebras.
We moreover show that they are graded cellular algebras, thus
making their cell modules, or standard modules, graded modules for the algebras.

\end{abstract}

\section{Introduction.}
In this paper we study
the Temperley-Lieb algebra $Tl_n(q)$.
It was introduced around forty years ago from considerations in statistical mechanics,
but has since turned out to be related to many topics of mathematics as well, including knot theory,
operator theory, algebraic
combinatorics and algebraic Lie theory. As of today, it is an object well known to a general audience in
physics as well as mathematics and at the same time it remains at
the center of a big number of research articles being published each year in both areas.

\medskip
Our main emphasis lies on a two-parameter
generalization $b_n(q,y_e)$ of the Temperley-Lieb algebra that was introduced by
P. Martin and H. Saleur in \cite{Mat-Sal}, as a way of introducing periodicity in the physical model.
An important feature of both $Tl_n(q)$ and $ b_n(q,y_e) $ is the fact that they are diagram algebras, that is they have bases
parameterized by certain planar diagrams, such that the multiplications are given by concatenation of these
diagrams. In the case of $Tl_n(q) $ these diagrams are the socalled {\it bridges} or {\it Temperley-Lieb diagrams},
in the case of
$ b_n(q,y_e) $ the diagrams are certain marked Temperley-Lieb diagrams and for this reason $ b_n(q,y_e) $
was called the blob algebra in \cite{Mat-Sal}.

\medskip

We are interested in the non-semisimple representation theory of $Tl_n(q)$ and $ b_n(q,y_e) $,
which is the case where $ q $ is specialized at a root of unity.
The $Tl_n(q)$-case is connected via Schur-Weyl duality to the representation theory of
the quantum group associated with $ SL_2$. The $ b_n(q,y_e) $-case is more intriguing
and has received quite a lot of attention over the last decade. It has been shown to share
a surprisingly big number of properties with objects that normally arise in Lie theory. In particular,
it was shown in \cite{martin-wood1} that the decomposition numbers are given by
evaluations at $1$ of certain Kazhdan-Lusztig polynomials associated with an infinite dihedral Weyl group.

\medskip
The fact that the decomposition numbers for $ b_n(q,y_e) $ come from polynomials gives
a first indication of the existence of a $ \Z$-graded structure on $ b_n(q,y_e) $ and on its standard modules,
and indeed a main goal of our paper is to construct such a graded structure on $ b_n(q,y_e) $.

\medskip
A main input to our paper comes from the seminal work of Brundan and Kleshchev that
constructs isomorphisms between cyclotomic Hecke algebras
and Khovanov-Lauda-Rouquier (KLR) algebras (of type $A$), see \cite{brundan-klesc}.
Since the KLR algebras are $ \Z $-graded, the various Hecke
algebras become $ \Z $-graded in this way as well. On the other hand, $ b_n(q,y_e) $ is
known to be a quotient of
the Hecke algebra $ \mathcal{H}_n(q,Q) $ of type $B$, indeed it is also sometimes referred to as the Temperley-Lieb
algebra of type $ B$. But $ \mathcal{H}_n(q,Q) $
is also the cyclotomic Hecke algebra of type $G(2,1,n)$ and our basic idea is now to exploit the result from
\cite{brundan-klesc} on this quotient construction.

\medskip
A big step towards our goal is taken already in
section 3 of our paper, where we show that the ideal $\mathcal{J}_n \subset \mathcal{H}_n(q,Q) $,
defining $ b_n(q,y_e) $, is graded,
thus making $ b_n(q,y_e) $ a $ \Z$-graded algebra. This
result relies on a realization of $\mathcal{J}_n$ due to
P. Martin and D. Woodcock in \cite{martin-wood}, in terms of certain explicitly given idempotents that
turn out to be well behaved with respect to the KLR-relations.

\medskip
On the other hand, this does not immediately imply a $ \Z$-grading on the standard modules
for $ b_n(q,y_e) $ and
indeed a major part of our paper is dedicated to this point.
An important ingredient to this comes from the recent paper by Hu and Mathas, \cite{hu-mathas},
that introduces the concept of a graded cellular algebra and shows that the cyclotomic Hecke algebras
are graded cellular with respect to the $\Z$-grading given by Brundan and Kleshchev's work. We then
achieve our goal in the sections 4-6 by showing that $ b_n(q,y_e) $ is a graded cellular algebra.

\medskip
A main difficulty in applying \cite{hu-mathas},  is due to the fact that the cell structure on
$ \mathcal{H}_n(q,Q) $ considered in \cite{hu-mathas} is related to the dominance order on bipartitions, which
is known to be incompatible with the natural order for the category of $ b_n(q,y_e) $-modules,
see \cite{Steen} and \cite{Steen1}.  We overcome this problem
by showing that $ b_n(q,y_e) $ is an algebra endowed with a family of Jucys-Murphy elements, in the sense
of Mathas \cite{Mat-So}, with respect to a natural order that we introduce in section 4.
This involves delicate arguments involving the diagram basis for $ b_n(q,y_e) $.

\medskip
It should be mentioned that our results are also valid in the Temperley-Lieb algebra case where
the relevant Hecke algebra $\mathcal{H}_n(q)$ this time is of type $ A$, and even in this
case our results seem to be new. On the other hand, in the Temperley-Lieb algebra
case there is actually a simpler way to show that the ideal of $\mathcal{H}_n(q)$ defining
$Tl_n(q)$ is graded. It is based on certain properties of Murphy's standard basis that were proved by M. H\"arterich
in \cite{Martin}.

\medskip
Let us sketch the layout of the paper. In the next section we introduce the various algebras that
play a role in the paper. In the third section we show that the ideals defining the Temperley-Lieb algebra
and the blob algebra are graded, which makes these algebras graded. In the following section we recall the diagrammatic
realizations of the Temperley-Lieb algebra and the blob algebra.
We here focus mostly on the blob algebra case. We introduce two ways of parametrizing the
blob diagrams, one via standard bitableaux of one-line bipartitions, the other via walks on the Bratteli diagram.
We also introduce an order relation $ \succ $ on the blob diagrams. In the fifth section we
show that the images in $ b_n(q,y_e) $ of the Jucys-Murphy elements of $\mathcal{H}_n(q,Q)$
make the blob algebra into an algebra with a family of Jucys-Murphy elements, in the sense of Mathas.
As we explain in the beginning of that section, this is quite surprising. We rely here on both combinatorial
descriptions of the blob diagrams. In the sixth section we obtain our main results, showing that the Temperley-Lieb
algebra and the blob are both graded cellular, and in 
the last section we give two examples illustrating our results.

\medskip
It is a pleasure to thank the referees for many suggestions that helped us improve the text.

\section{Notation and setup.} \label{}
In this section we fix the notation that is used throughout the paper. We introduce the algebras
to be studied, the Temperley-Lieb algebra, the blob algebra, the corresponding Hecke
and Khovanov-Lauda-Rouquier algebras and recall the relevant results from the literature involving them.
The important diagrammatic realizations of the Temperley-Lieb algebra and the blob algebra are postponed
to section 4.

\medskip
Throughout the paper the ground field is the complex field  $\C$ although some of our results hold in greater
generality. For $q\in \C^{\times}$
and an integer $k$ we define
$$[k]=[k]_q:= q^{k-1} + q^{k-3} + \ldots + q^{-k+1} \in \F$$
the usual Gaussian integer. All our algebras are associative and unital.

\subsection{The Temperley-Lieb algebra, the blob algebra, the Hecke algebras.}
\begin{defi} \rm
Let $ q \in \C^{\times}$.
The Temperley-Lieb algebra $Tl_{n}(q)$ is the $\F$-algebra on the generators $U_1,...,U_{n-1}$ subject to the relations
\begin{align*}
U_i^{2}& =-[2] U_i &  & \mbox{if }  1 \leq i \leq n-1 \\
U_iU_{j}U_{i}&=U_i &   & \mbox{if }   |i-j| =1 \\
U_iU_{j}&=U_{j}U_i &   & \mbox{if }   | i-j | >1.
\end{align*}
\end{defi}

The main object of the paper is the blob algebra, introduced in \cite{Mat-Sal} by  P. Martin and H. Saleur
as a generalization of the Temperley-Lieb algebra.
Let $y_e$ be an invertible element of $\F$.

\begin{defi} \rm  \label{Defblob}
The blob algebra $b_n(q,y_e)$ is the $\F$-algebra on the generators $e,U_1,...,U_{n-1}$ subject to the relations
\begin{align*}
U_{i}^{2}&=-[2] U_i &  &\mbox{if }  1 \leq i \leq n-1 \\
U_iU_{j}U_{i}&=U_i &   & \mbox{if }   |i-j| =1 \\
U_iU_{j}&=U_{j}U_i &   & \mbox{if } | i-j | >1 \\
U_1eU_{1}&=y_e U_1 &   &  \\
e^{2}&=e & &  \\
U_ie&=eU_i & & \mbox{if }  2 \leq i \leq n-1.
\end{align*}
\end{defi}
Assume that $[m]\neq 0$. The parametrization of $b_n(q,y_e)$ through $  y_e = \frac{-[m-1]}{[m]} $
includes the non-semisimple cases,
see \cite[Section 2]{martin-wood1}. Under this choice of $y_e$ we denote $b_n(q,y_e)$
for the rest of the paper by
$\blob{m}$ and replace $e$ by the rescaled generator $U_0:=-[m]e$.

\medskip
The Temperley-Lieb algebra and the blob algebra
were introduced from motivations in statistical mechanics. An important feature,
that we postpone to the next section, is that they both have diagrammatic realizations
by planar diagrams.

\medskip
We next define the related Hecke algebras.

\begin{defi} \rm  \label{HeckeA}
Let $ q \in \C$ and assume that $ q \neq 0,1$.
The Hecke algebra $\mathcal{H}_n(q)$ of type $A_{n-1}$ is the
$\F$-algebra with generators $T_1,\ldots ,T_{n-1}$, subject to the relations
\begin{align}
 (T_i-q)(T_i+1)&=0  & &  \text{ for } 1\leq i \leq n-1   \\
 T_iT_{i+1}T_{i}&=T_{i+1}T_{i}T_{i+1}   &  &  \text{ for } 1\leq i \leq n-2 \\
 T_iT_j  & =T_jT_i   &  &  \text{ for }    | i-j| >1
\end{align}
\end{defi}
It follows easily from the relations that $T_r$ is an invertible element in $\mathcal{H}_n(q)$,
with $T_r^{-1} = q^{-1}(T_r -q +1 )$.
We define elements $L_1, \ldots , L_n \in \mathcal{H}_n(q) $ by $ L_1 := 1 $ and recursively
$L_{r+1} = q^{-1} T_r L_r T_r $ for all admissible $r$. They are the first examples of Jucys-Murphy elements
that play an important role in our paper.

\begin{defi} \rm \label{hecke algebra}
Let $q, \lambda_1,\lambda_2 \in \F$ and suppose that $ q\neq 0,1$.
The cyclotomic Hecke algebra ${\mathcal{H}}_n(q;\lambda_1 ,\lambda_2)$ of type $G(2,1,n)$ is
the $\F$-algebra with generators $L _1 ,\ldots,L_n ,$ $T_1 ,\ldots ,T _{n-1}$ and relations
\begin{align*}
(L _1 -\lambda_1)(L _1 -\lambda_2) &= 0, & L_rL_s&=L_sL_r,  \\
(T_r + 1)(T_r - q)  &= 0, &  T_r L_r &  = L_{r+1} (T_r - q + 1),   \\
T _s T_{s+1} T_s  &= T _{s+1} T_s T_{s+1}, & &  \\
T_r L_s   &= L_s T_r , & \mbox{if } & |r - s| > 1, \\
T_r T_s    &= T_s T_r , & \mbox{if } & s \neq r,r + 1
\end{align*}
for all admissible $r, s$.
\end{defi}
Once again, $T_r$ is invertible with $T_r^{-1} = q^{-1}(T_r -q +1 )$.
From this one gets
that
\begin{equation} \label{jucysmurphy como t}
L_{r+1} = q^{-1} T_r L_r T_r
\end{equation}
Moreover, it follows from the relations that $ f(L_1, \ldots, L_n) $ is a central element of ${\mathcal{H}}_n(q;\lambda_1 ,\lambda_2)$
for $ f(x_1, \ldots,  x_n ) $ a symmetric polynomial. These $ L_i $ are also called Jucys-Murphy elements.

\medskip
We now explain certain relations between the algebras that we have defined.

\begin{teo} \label{TLisomorphism} There are surjections $\Phi_1 $ and $\Phi_2 $ given by
$$\begin{array}{cccl cccl}
\Phi_1 :  &  \mathcal{H}_n(q^2) & \longrightarrow  & Tl_{n}(q),  &
  & T_i  & \mapsto  &   qU_i + q^2\\
\Phi_2 :  &  \mathcal{H}_n(q^2) & \longrightarrow  & Tl_{n}(q),  &
  & T_i  & \mapsto  &   -qU_i -1.
\end{array} $$
The kernel of $ \Phi_1 $ is the ideal  generated by
$$ q^{-6} T_1 T_2 T_1 - q^{-4} T_1 T_2  - q^{-4} T_2 T_1 + q^{-2} T_1   + q^{-2} T_2  -1 $$
and the kernel of $ \Phi_2 $ is the ideal generated by
$$ T_1 T_2 T_1 + T_1 T_2  + T_2 T_1 +  T_1   +  T_2  +1. $$
\end{teo}
\begin{dem}
This is well known.
\end{dem}

\medskip
There are two, not obviously equivalent, ways to generalize this Theorem to the blob algebra case.
One is given in \cite{gra ler 1}, but for our purposes it is more convenient to work with the second one, that
appears in \cite{martin-wood}.
Set $ Q := q^m$ and
define $ \hecke= {\mathcal{H}}_n(q^{2};Q, Q^{-1})$. Assume
\begin{align}
 \label{qym} q^{4}\neq 1,&  & Q \neq Q^{-1}, &  & Q \neq q^{2} Q^{-1}, & & Q^{-1} \neq q^{2} Q.
\end{align}
With the above conditions, one can define elements $ e_2^{-1}, e_2^{-2} \in \heckedos$ by the formulas
\begin{equation*}
 e_2^{-1} =\frac{  (T_1-q^{2})(L_1-Q^{-1} ) (L_2-Q^{-1} )} {(1+q^2)(Q-Q^{-1} )(Q^{-1}-q^{-2} Q)}
\end{equation*}
\begin{equation*}
e_2^{-2} =\frac{(T_1-q^{2})(L_1-Q)(L_2-Q)}   {(1+q^2)(Q^{-1}-Q)(Q-q^{-2} Q^{-1})}.
\end{equation*}

Note that $(L_1-Q)(L_2-Q)$ and $(L_{1}-Q^{-1})(L_2-Q^{-1})$ are symmetric polynomials in $L_1$ and $L_2$. Therefore, they are central elements in $\mathcal{H}_2(m)$. Using this and $ L_2 = q^{-2} T_1 L_1 T_1$, one finds that they verify the following equations
\begin{align}
\label{carac e1} (T_1+1)e_2^{-1} & =0,  &  (T_1+1)e_2^{-2}&=0,  \\
\label{carac e2} (L_1-Q)e_2^{-1} & =0,  & (L_1-Q^{-1})e_2^{-2} & =0,\\
\label{carac e3} (L_2-Qq^{-2})e_2^{-1} & =0,  & (L_2-Q^{-1}q^{-2})e_2^{-2} & =0
\end{align}
and from this it follows that $ e_2^{-1} $ and $ e_2^{-2} $ are idempotents associated with irreducible representations of
$ \heckedos $ of dimension one.
Note that $ e_2^{-1} $ and $ e_2^{-2} $ are the unique idempotents
satisfying (\ref{carac e1}) and (\ref{carac e2}).
For all $ n$ there is a canonical embedding $\hecke \hookrightarrow \heckenmasuno $. Using it repeatedly
we consider $ e_2^{-1} $ and $ e_2^{-2} $ as elements of $\hecke$ and
denote by ${\mathcal{J}}_n$ the ideal of $\hecke$ generated by them.

\begin{teo} \label{isomorphism} The map $\Phi $ given by
$$\begin{array}{cccl}
\Phi :  &  \hecke & \longrightarrow  & \blob{m}  \\
  & T_i -q^{2} & \mapsto  &   qU_i \\
  & L_1-q^{m}  & \mapsto &   (q-q^{-1})U_0
\end{array} $$
induces a $\C$-algebra isomorphism between $\hecke / \mathcal{J}_n$ and $\blob{m}$.
\end{teo}
\begin{dem}
See \cite[Proposition 4.2]{martin-wood}.
\end{dem}

\medskip
We would like to have an integral version of the last result, but want also to avoid
those choices of the parameters that correspond to the conditions (\ref{qym}).
This can for example be achieved by localizing $  \C[q, q^{-1}, Q, Q^{-1} ] $ conveniently. To be precise,
we choose for $ R $ the localization of the Laurent polynomial ring $  \C[q, q^{-1}, Q, Q^{-1} ] $
at $ S$, defined as the multiplicatively closed subset
of $  \C[q, q^{-1}, Q, Q^{-1}  ] $ generated by the polynomials
$ 1, q^4-1,   Q- Q^{-1},     Q- Q^{-1} q^{2} $ and $ Q^{-1}- Q q^{2} $.
For integers $ l$ and $ m$ we denote by
$ {\m} $ the ideal $ \langle q - e^{ 2 \pi i/ l } ,  Q - q^m \rangle $ of $ R   $.
Then we have that either $ {\m} = R $ or else
$  {\m} $ is a maximal ideal in $ R $. In the last case we
define $ \OO := R_{\m} $ and get
that $ \OO $ is a discrete valuation ring with maximal ideal $ {\m} $,
quotient field $\K:= \C(q,Q)$ and residue field $ \OO/ {\m} =  \C $ containing the $l$'th root of unity $q$.

\medskip

Throughout the paper we assume that $ \OO, \K $ and $ \C$ are chosen
as above, and furthermore, in order to simplify notation, that $ l $ is odd.
In the next subsection we
recall the $\Z$-grading on $\mathcal{H}_n(q^2)$ and
$\mathcal{H}_{n}(m)$ given by Brundan and Kleshchev in
\cite{brundan-klesc}. Note that since $ l $ is assumed to be odd,
the condition from {\it  loc. cit.}
that $q^{m}$ be a power of $q^{2}$,
or equivalently,
that the congruence
$ 2k \equiv m \mod l $ be solvable, is always fulfilled.

\medskip
Recall that the quantum characteristic of an element $ q $ of a field $ F$ is the smallest positive integer $k$
such that $1 + q + \ldots +q^{k-1}=0$, setting $k=0$ if no such integer exists.  With our choice of $ q \in \F$
the quantum characteristic is $ l$.
We set $I=\mathbb{Z} / l\mathbb{Z}$ and
refer to $ I^n$ as the residue sequences of length $ n$. Note that in order to apply \cite{brundan-klesc},
we should acually use the quantum characteristic of $ q^2$ in the definition of $ I $, but
since $ l $ is assumed to be odd, the two definitions coincide.

\medskip

\medskip
We now
define $ \blobR{m} $ as the $ \OO $-algebra on generators $e,U_1,...,U_{n-1}$ subject to the same relations as for $ \blob{m}$.
Then $ \blobR{m}  $ is free over $ \OO $ as
can be seen using the results of the appendix of \cite{blob positive}, note that they are valid over
any commutative ring. The rational blob algebra $ \blobK{m} $ is defined the same way,
and we have base change isomorphisms $ \blobR{m}  \otimes_{\OO} \C = \blob{m}$ and
$ \blobR{m}\otimes_{\OO}  K = \blobK{m}  $.
Finally we define
$\heckeintegral$ as the $ \OO$-algebra on generators $L _1 ,\ldots,L_n ,$ $T_1 ,\ldots ,T _{n-1}$
subject to the same relations as for $ \hecke$, but using parameters $ \lambda_1 = Q $ and $\lambda_2 = Q^{-1} $.
Similarly, we define $\heckerational$ and we have base change isomorphisms as above.

\begin{teo} \label{isomorphismintegral}
There is a surjection $ \Phi: \heckeintegral  \longrightarrow \blobR{m} $.
\end{teo}
\begin{dem}
The argument given in \cite[Proposition 4.2]{martin-wood} involves verification of blob algebra relations and therefore
gives a surjection $ \heckeintegral  \longrightarrow \blobR{m} $, as claimed.
\end{dem}

\subsection{The Khovanov-Lauda-Rouquier algebra.}
In the following $ \cal H $ refers to either $\mathcal{H}_n(q^2)$  or $\hecke$ (with $ q \in \F$ chosen as above).
Let $M$ be a finite dimensional
$\cal H$-module.
By \cite[Lemma 7.1.2]{kle} the eigenvalues of each
$L_r$ on $M$ are of the form $q^{2i}$ for $i \in I$. So $M$ decomposes as the direct sum
$ M =\bigoplus _{\bs{i} \in I^{n}} M_{\bs{i}}$ of its generalized weight spaces
$$  M_{\bs{i}}:= \{  v\in M \text{  } \mid \text{  } (L_r-q^{2i_r})^{k}v=0 \text{ for } r=1,\ldots , n \text{ and } k \gg 0  \}. $$
In particular, taking $M$ to be the regular left module $ \cal H $, we obtain a system
$\{ e(\bs{i})\mid  \bs{i} \in I^{n}  \}$
of mutually orthogonal idempotents in $ \cal H $ such that $e(\bs{i})M=M_{\bs{i}}$ for each
$M$ as above.

\medskip
We can now define nilpotent elements $ y_1,\ldots ,y_n \in \cal H$ via the formula
\begin{equation}{\label{nilpotent}}
y_r=\sum_{\bs{i} \in I^n} (1-q^{-2i_{r}}L_r)e(\bs{i}).
\end{equation}
For $1\leq r < n$ and $\bs{i} \in I^{n}$,
Brundan and Kleshchev define in \cite{brundan-klesc} certain formal power series,
$P_r(\bs{i}), Q_r(\bs{i}) \in \F[[y_r, y_{r+1}]]$,
such that $Q_r(\bs{i})$ has non-zero constant term, see
 \cite[(4.27) and (4.36)]{brundan-klesc} for the explicit formulas.
Since each $y_r$ is nilpotent in $\cal H$, we can consider
$P_r(\bs{i})$ and $Q_r(\bs{i})$
as elements of $ \cal H$, with $Q_r(\bs{i})$ invertible.
We then set
\begin{equation}  \label{definicion de psi}
    \psi_r=\sum_{\bs{i} \in I^n} (T_r+P_r(\bs{i}))Q_r(\textbf{i})^{-1}e(\bs{i}).
\end{equation}
The main theorem in \cite{brundan-klesc} gives a presentation of $ \cal H $ in terms of the elements
$$ \{ \psi_1,\cdots , \psi_{n-1} \} \cup  \{ y_1,\cdots , y_{n} \}    \cup
  \{ e(\textbf{i}) \mbox{         }  | \mbox{         } \textbf{i} \in I^n \}                $$
and a series of relations between them that we describe shortly. An important point of these
relations is that they are homogeneous with respect to a nontrivial $\mathbb{Z}$-grading on $ \cal H$.
To describe the $\mathbb{Z}$-grading it is convenient to introduce the matrix $(a_{ij})_{i,j\in I}$, given by
\begin{align*}
a_{ij}& =\left\{ \begin{array}{l}
2 \\
0   \\
-1   \\
\end{array}
\right.   & \begin{array}{l}
              \mbox{if } i=j \mod l  \\
              \mbox{if } i\neq j\pm 1 \mod l \\
              \mbox{if } i= j\pm 1 \mod l .
            \end{array}
\end{align*}
With this at hand, we are now able to state \cite[Main Theorem]{brundan-klesc}. The
Theorem holds in greater generality than shown here, namely for
all cyclotomic Hecke algebras, including the degenerate algebras, but for our purpose the following version is enough.

\begin{teo} \label{definkl}
The algebra $\cal H$ is isomorphic to a cyclotomic Khovanov-Lauda-Rouquier algebra of
type $ A$.
To be precise, it is isomorphic to the $ \C$-algebra generated by
$$ \{ \psi_1,\cdots , \psi_{n-1} \} \cup  \{ y_1,\cdots , y_{n} \}    \cup
  \{ e(\textbf{i}) \mbox{         }  | \mbox{         } \textbf{i} \in I^n \}                $$
subject to the following relations for $\textbf{i},\textbf{j}\in I^n$ and all admissible $r, s$
\begin{align}
\label{aaaa}       y_1e(\textbf{i})   = 0 \mbox{ if } i_1 =    & \left\{ \begin{array}{r}
\pm k \mod{l}   \\
0 \mod{l} \\
\end{array}
\right. &  \begin{array}{l}
\mbox{if  }   {\cal H } = \hecke \\
\mbox{if  }    {\cal H } =  \mathcal{H}_n(q^2)  \\
\end{array} \\
\label{bbbb}       e(\textbf{i})   = 0 \mbox{ if } i_1 \neq    & \left\{ \begin{array}{r}
\pm k \mod{l}   \\
0 \mod{l} \\
\end{array}
\right. &  \begin{array}{l}
\mbox{if  }   {\cal H } = \hecke \\
\mbox{if  }    {\cal H } =  \mathcal{H}_n(q^2)  \\
\end{array} \\
\label{kl2}e(\textbf{i})e(\textbf{j})& =\delta_{\textbf{i,j}} e(\textbf{i}), & \\
\label{kl3}\sum_{\textbf{i} \in I^n} e(\textbf{i})& =1, \\
\label{kl4}y_{r}e(\textbf{i})& =e(\textbf{i})y_r ,& \\
\label{kl5}\psi_r e(\textbf{i})& =e(s_r \textbf{i}) \psi_r, \\
\label{kl6}y_ry_s& = y_sy_r,&
&   \\
\label{kl7}\psi_ry_s& = y_s\psi_r,&
 \mbox{ if } s\neq r,r+1 \\
\label{kl8}\psi_r\psi_s& = \psi_s\psi_r,&
 \mbox{ if } |s-r|>1 \\
\label{kl9}\psi_ry_{r+1}e(\textbf{i})& =\left\{ \begin{array}{l}
(y_r\psi_r +1)e(\textbf{i}) \\
y_r\psi_r e(\textbf{i})   \\
\end{array}
\right.  &\begin{array}{l}
\mbox{if  }  i_r=i_{r+1}  \mod l  \\
\mbox{if  }   i_r \neq i_{r+1}   \mod l  \\
\end{array} \\
\label{kl10}y_{r+1}\psi_re(\textbf{i})& =\left\{ \begin{array}{l}
(\psi_ry_r +1)e(\textbf{i}) \\
\psi_ry_r e(\textbf{i})   \\
\end{array}
\right. &  \begin{array}{lcc}
\mbox{if  }  i_r=i_{r+1}  \mod l  \\
\mbox{if  }   i_r \neq i_{r+1}  \mod l  \\
\end{array}
\end{align}
\begin{align}
\label{kl11}\psi_r^{2}e(\textbf{i})& =\left\{ \begin{array}{l}
0   \\
e(\textbf{i})  \\
(y_{r+1}-y_{r})e(\textbf{i}) \\
(y_{r}-y_{r+1})e(\textbf{i})  \\
\end{array}
\right.
\begin{array}{l} \,\,\,\,\,\,\,\,\,\,\,\,\,\,\,\,\,\,
\mbox{if  }  i_r=i_{r+1} \, \, \, \, \, \, \, \, \, \, \mod l  \\ \,\,\,\,\,\,\,\,\,\,\,\,\,\,\,\,\,\,
\mbox{if  }   i_r \neq i_{r+1} \pm 1 \mod l   \\\,\,\,\,\,\,\,\,\,\,\,\,\,\,\, \,\,\,
\mbox{if  }  i_{r+1}= i_r+1  \mod l  \\\,\,\,\,\,\,\,\,\,\,\,\,\,\,\, \,\,\,
\mbox{if  }  i_{r+1}= i_r-1  \mod l \\
\end{array}
\\
\label{kl12}\psi_r\psi_{r+1}\psi_re(\textbf{i})& =\left\{ \begin{array}{l}
(\psi_{r+1}\psi_r\psi_{r +1} +1)e(\textbf{i})  \\
(\psi_{r+1}\psi_r\psi_{r +1} -1)e(\textbf{i})    \\
(\psi_{r+1}\psi_r\psi_{r +1} )e(\textbf{i})    \\
\end{array}
\right.
\begin{array}{l}
\mbox{if  }  i_{r+2}=i_r=i_{r+1}-1  \mod l  \\
\mbox{if  } i_{r+2}=i_r=i_{r+1}+1  \mod l    \\
\mbox{otherwise   }        \\
\end{array}
\end{align}
where $ s_r := (r,r+1)$ is the simple transposition acting in $ I^n $
by permutation of the coordinates $ r, r+1$ and $k\in \mathbb{Z}$ is such that $2k \equiv  m \mod l$.
The isomorphism maps each of the generators to the element of $\cal H$ that has the same name.
The conditions
$$   \begin{array}{ccccc}
       deg\mbox{  } e(\textbf{i})=0,  &  &  deg\mbox{  } y_r=2,  &  &
 deg\mbox{  }\psi_se(\textbf{i})=-a_{i_{s},i_{s+1}} \end{array}
$$
for $1 \leq r\leq n$,  $1 \leq s \leq n-1 $ and $\textbf{i} \in I^{n}$
define a unique $\mathbb{Z}$-grading on $ \cal H $ with degree function $ deg$.
\end{teo}

\medskip
Following \cite{hu-mathas} we shall refer to the $ e(\bs{i})  $ as
the KLR-idempotents.
In the following, all statements involving a grading on $ \cal H $ refer to the above Theorem.
Note that although
the elements $L_r$ and $T_r$ are not homogeneous in $ \cal H$, they
can be expressed in terms of homogeneous generators in the following way, see equations (4.42) and
(4.43) of \cite{brundan-klesc}:
\begin{align}
\label{Lsubr} L_r & =\sum_{\bs{i} \in I^n} q^{2i_{r}}(1-y_r)e(\bs{i}) \\
\label{Tsubr} T_r & =\sum_{\bs{i} \in I^n} (\psi_rQ_r(\bs{i})-P_r(\bs{i}))e(\bs{i}).
\end{align}

\subsection{Partitions and tableaux.}
We finish this introductory section by recalling some basic combinatorial notions related to the symmetric group.
Whenever we work with a partially ordered set $ X $ with order relation say $ \geq $,
we write $x_1 > x_2$, that is omit the lower line of the order symbol,
if $x_1,x_2 \in X$ satisfy $x_1 \geq x_2$ and $x_1 \neq x_2$.
Let $ n $ be a positive integer.
A(n integer) partition of $n$ is a sequence $\lambda=(\lambda_1,\lambda_2,\ldots)$ of non-negative integers such that
$|\lambda |:= \sum_i\lambda_i =n$ and $\lambda_i \geq \lambda_{i+1}$ for all $i\geq 1$. The Young diagram of
$ \lambda$ is the set
$$[\lambda] =\{ (i,j) \in  \mathbb{N}\times \mathbb{N}  \mbox{  } | \mbox{   } 1\leq j \leq \lambda_i \text{ and } i\geq 1 \}.$$
The elements of it are called nodes or entries.
It is useful to think of $ [\lambda] $ as an array of boxes in the plane, with the indices following
matrix conventions. Thus the box with label $ (i, j)$ belongs to the $i$'th row and $ j$'th column.
If $\lambda$ is a partition of $n$ we denote by $\lambda'$ the partition of $n$ obtained from $\lambda$ by interchanging
its rows and columns. A two-column partition of $n$ is a partition  $\lambda $ of $n$ such that $\lambda_i \leq 2 $ for all $i\geq 1$.
The set of all partitions of $n$ is denoted $\pa{n}$ and the set of two-column partitions of $n$ is denoted by $\pacol{n}$.
A $\lambda$-tableau is a bijection $\tau : [\lambda] \rightarrow \{1,\ldots,n \} $. We say that $\tau $ has
shape $\lambda$ and write
$\Shape(\tau)=\lambda $. We think of it as a labeling of the diagram
of $\lambda$ using elements from $ \{1, 2, \ldots, n \}$
and in this way we
can talk of the rows and columns of a tableau
as subsets of $ \{1, 2, \ldots, n \}$.
We say that $\tau $ is row (resp. column) standard if the
entries of $\tau$ increase from left (resp. top) to right (resp. bottom) in each row (resp. column).
$\tau $ is standard if it is row standard and column standard. The set of all standard $\lambda$-tableau is denoted by $\Std(\lambda)$ and the union of all $\Std(\lambda)$ is denoted $\Std(n)$.

\medskip
Assume that $\lambda , \mu \in \pa{n}$. We say that $\lambda$ dominates $\mu$ and
write $\lambda \trianglerighteq \mu$ if
$$ \sum_{i=1}^{j} \lambda_i \geq \sum_{i=1}^{j} \mu_{i}$$
for all $j\geq 1$. Then $\pa{n}$ becomes a partially ordered set via $ \trianglerighteq $. It can be extended to $\Std(n)$
as follows.
For $\sigma, \tau \in  \Std(n)$,
we say that $\sigma $ dominates $\tau$ and write $\sigma \trianglerighteq \tau $ if $Shape(\sigma_{|_k}) \trianglerighteq  \Shape(\tau_{|_k})$, for $k=1,\ldots,n$, where $\sigma_{|_k}$ and $\tau_{|_k}$ are the
tableaux obtained from $\sigma$ and $\tau$ by removing the entries greater than $k$.

\medskip
Let $\tau^{\lambda}$ be the unique standard $\lambda$-tableau such that $\tau^{\lambda} \trianglerighteq \tau$ for all $\tau \in \Std(\lambda)$. In $\tau^{\lambda}$ the numbers $ 1, 2, \ldots, n $ are filled in
increasingly along the rows from top to bottom.
The symmetric group ${\mathfrak{S}}_n$ acts on the left on
the set of $\lambda$-tableaux permuting the entries.
For $\tau \in \Std(\lambda)$, we denote by $d(\tau)$ the permutation of
${\mathfrak{S}}_n$ that satisfies $\tau = d(\tau) \tau^{\lambda}$.

\section{Grading the Temperley-Lieb algebra and the blob algebra.} \label{}
In this section we show that the Temperley-Lieb algebra $Tl_{n}(q)$ and the blob algebra $ \blob{m} $ are $\Z$-graded
algebras. We do this by proving that the kernels of the surjections given in Theorem \ref{isomorphism}
and Theorem \ref{TLisomorphism} are graded ideals. In the $Tl_{n}(q)$-case we rely on certain properties of
Murphy's standard basis that are proved in \cite{Martin}. These properties are missing in the $ \blob{m}$-case and
so our argument is somewhat different in that case.
Let $ A = \oplus_{ n \in \Z}  A_n $ be a $ \Z$-graded ring with homogeneous parts $ A_n$.
Recall that $ I \subset A $ is called a
graded (homogeneous) ideal of $ A $ if it is an ideal and if $ I = \oplus_{ n \in \Z}  I_n$
where $ I_n := A_n \cap I$. If $ I $ is a graded ideal of $ A $ then the quotient $ A/I $ becomes a $ \Z$-graded ring as well,
with homogeneous parts $ A_n / I_n $.
We need the following Theorem.

\begin{teo} \label{chevalley}
Let $ A $ be a $\mathbb{Z}$-graded algebra. Assume that $I$ is an ideal of $A$ that is generated by
homogeneous elements. Then $ I $ is graded.
\end{teo}

\begin{dem}
See \cite[Theorem 1.3]{chevalley}.
\end{dem}

\subsection{Grading $Tl_n(q)$.}
Let us briefly recall Murphy's standard basis for the Hecke algebra $ \mathcal{H}_n(q^{2})$.
For $w=s_{i_1} \ldots s_{i_{k}}$ a reduced expression of $w \in {\sub{S}}_{n}  $
we define $T_w: =T_{i_1}\ldots T_{i_k} $. Then $ \{ T_w | w \in {\sub{S}}_{n} \}$ is a basis for $ \mathcal{H}_n(q)$.
For $\lambda \in \pa{n}$ we let ${\sub{S}}_{\lambda} \leq {\sub{S}}_{n} $ denote the row stabilizer
of $\tau^{\lambda}$ under the left action of $ {\sub{S}}_{n} $ on tableaux
and define
$$x_{\lambda}:=\sum_{w\in {\sub{S}}_{\lambda}} T_{w} $$
We let $\ast$ denote the anti-automorphism of $\mathcal{H}_n(q^2)$
determined by $T_i^{\ast}=T_i$ for all $1\leq i < n$
and define for
$\sigma , \tau \in \text{Std}(\lambda)$
$$  x_{\tau \sigma} = T_{d(\tau)}^{\ast} x_{\lambda}  T_{d(\sigma)}.
$$
Then $ \{ x_{\tau \sigma} \} $, with $\tau $ and $\sigma $ running over
standard tableaux of the same shape,
is Murphy's standard basis for $\mathcal{H}_n(q^2)$, see \cite[Theorem 4.17]{Murphy1}.

\medskip
We set
$\mathcal{I}_n :=  \ker \Phi_2 $
where $ \Phi_2 : \mathcal{H}_n(q^2)  \longrightarrow   Tl_{n}(q) $
is the second surjection given in Theorem \ref{TLisomorphism}.
Then $\mathcal{I}_n $ is an ideal of $ \mathcal{H}_n(q^2)  $ and we have $  \mathcal{H}_n(q^2) / \mathcal{I}_n  = Tl_{n}(q) $.
We can now state our first Theorem.
\begin{teo} \label{first_theorem}
$\mathcal{I}_n $ is a graded ideal of $ \mathcal{H}_n(q^2)  $. Hence $ Tl_{n}(q) $ is a $ \Z$-graded
algebra, with the grading induced from the one on $ \mathcal{H}_n(q^2)  $, via Theorem \ref{definkl}.
\end{teo}
\begin{dem}
We first note that by the results of H\"arterich, \cite[Theorem 4]{Martin}, we know that
$\mathcal{I}_n$ is spanned (over $\C$!) by those $ \{ x_{\tau \sigma} \} $ for which the underlying shape
has strictly more than two columns, that is
$ \Shape(\tau), \Shape(\sigma)  \notin \pacol{n}$.
In other words, $ \{ x_{\tau \sigma} \,  | \,  \sigma, \tau \in \Std(\lambda), \lambda \in \pa{n} \setminus \pacol{n} \} $
is a basis for $\mathcal{I}_n$.

\medskip
On the other hand,
in \cite{hu-mathas} J. Hu and A. Mathas construct a basis
$\{ \psi_{\tau \sigma} \} $
for $ \mathcal{H}_n(q^2)  $,
such that each $ \psi_{\tau \sigma} $ is a homogeneous element of $ \mathcal{H}_n(q^2)  $;
here $ (\tau, \sigma ) $ is running over the same set as for the Murphy's standard basis.
They furthermore show in \cite[Lemma 5.4]{hu-mathas} that for each
pair $ (\tau, \sigma ) $ like this, there is a non-zero scalar $c\in \C$
such that
\begin{equation} \label{matrix graded cellular a murphy}
\psi_{\tau \sigma } =cx_{\tau \sigma } + \sum_{(\upsilon,\varsigma) \vartriangleright ( \tau ,\sigma  )} r_{\upsilon \varsigma} x_{\upsilon \varsigma}
\end{equation}
where $r_{\upsilon \varsigma} \in \F$ and where $ (\upsilon,\varsigma) \trianglerighteq ( \tau ,\sigma  )$
by definition means that $ \upsilon \trianglerighteq \tau  $ and
$\varsigma \trianglerighteq \sigma   $.
But this shows that also the $ \{ \psi_{\tau \sigma} \} $ such that
$ \Shape(\tau), \,  \Shape(\sigma)  \notin \pacol{n}$ are a basis for
$\mathcal{I}_n$. From this we get,
via Theorem
\ref{chevalley}, that $\mathcal{I}_n$
is a graded ideal as claimed.
\end{dem}

\begin{rem}\rm
There is a version of the Theorem involving the homomorphism $ \Phi_1 $. For this, in the proof one should replace
$ \{ \psi_{\tau \sigma} \} $ by the dual basis $ \{ \psi_{\tau \sigma}^{\prime} \} $ of \cite{hu-mathas}.
\end{rem}

\begin{rem}\rm
In spite of the importance of the Temperley-Lieb algebra in mathematics and physics,
the above graded structure has not been mentioned before in the literature, to the best of our knowledge.
For example, in the categorification of the Temperley-Lieb algebra considered in
\cite{Stroppel}, the parameter $ q $ is not a root of unity.
The same remark applies to the supergrading used in \cite{Zhang}.

\end{rem}

\subsection{Grading  $\blob{m}$.}
Let us now turn to the blob algebra.
In order to treat that case we need the following Theorem. Note that since $ l $ is odd,
 there is always $ k \in \mathbb Z $ satisfying the condition of the Theorem.
\begin{teo} \label{idempotentes homogeneous}
Let $k \in \mathbb{Z}$  such that $2k \equiv m \mod{l} $. Then, the elements
$e_2^{-1},e_2^{-2} \in \hecke $ are homogeneous of degree zero.
More precisely, they can be written as a sum of homogeneous elements of degree zero as follows
\begin{align} \label{relations blob}
e_2^{-1}=& \sum_{\bs{i}}e(\bs{i}) & e_2^{-2}=& \sum_{\bs{j}} e(\bs{j})
\end{align}
where the left sum runs over all $\bs{i} \in I^{n}$ such that $i_1=k$ and $i_2=k-1$, and the right sum runs over all $\bs{j} \in I^{n}$ such that $j_1=-k$ and $j_2=-k-1$.
\end{teo}

\begin{dem}
We only prove the result for $e_2^{-1}$, the result for $e_2^{-2}$ is proved
similarly.

In \cite[Section 4.4]{bkw graded}, Brundan, Kleshchev and Wang note that under the
embedding $ \hecke \hookrightarrow \heckenmasuno$ one has
$e(\bs{i}) \mapsto \sum_{i\in I} e(\bs{i},i)$, and so it is enough to prove the case $n=2$,
that is that $ e_2^{-1} = e(k,k-1)$ holds. Using the uniqueness statement for $ e_2^{-1} $,
in order to prove this, it is enough to
show that $ e(k,k-1) $ verifies the equations (\ref{carac e1})
and (\ref{carac e2}), since it is clearly an idempotent.

Note first that $y_1=0$ as it follows by combining the relations (\ref{aaaa}), (\ref{bbbb})
and (\ref{kl3}).
Put now $\bs{j} =(k,k-1)$. Multiplying (\ref{Lsubr})  by $ e(\bs{j}) $ for $ n= 2 $ and $ r= 1 $,
we get $ L_1 e(\bs{j})  = q^{2k} e(\bs{j}) $, or equivalently
$ L_1 e(\bs{j})  = q^{m} e(\bs{j}) $. Hence (\ref{carac e2}) holds.

To show (\ref{carac e1}) we first recall
from \cite[Lemma 4.1(c)]{hu-mathas} that in general $ e(\bs{i}) \neq 0$ iff $\bs{i} \in I^n $ is
a residue sequence coming from a standard bi-tableau of a bipartition of $ n $. Combining this fact with the standing conditions
on $ q $ given in (\ref{qym}), we deduce $e(s_1\bs{j})=0$ and hence
$  \psi_1e (\bs{j}) = 0 $ by (\ref{kl5}).
Multiplying this equation on the left by $ \psi_1 $ and using (\ref{kl11}) we obtain $ y_2e (\bs{j}) = 0 $.

Now, recall that by definition $ P_1 (\bs{j}) $ and $ Q_1 (\bs{j}) $ are power series in $ y_1 $ and $ y_2 $.
Furthermore, in this particular case we have that the constant coefficient of $ P_1 (\bs{j}) $ is $1$ and so
(\ref{Tsubr}) gives $ (T_1+1) e (\bs{j}) = 0 $ as needed.
\end{dem}

\medskip
Later on, in Remark \ref{at-this-point}, we indicate an alternative
proof of the Theorem that uses seminormal bases.

\medskip
We are now in position to establish the main objective of this section,
namely to provide a graded structure on
$\blob{m} $. In the forthcoming Section \ref{a graded cellular}, we refine this graded structure on $ \blob{m} $
to a graded cellular basis structure.

\begin{cor} \label{main}
The kernel of the surjection $ \Phi :   \hecke  \longrightarrow   \blob{m}  $ from Theorem \ref{isomorphism} is a
graded ideal.
Hence,
the algebra $\blob{m}$ has a presentation with generators
$$ \{ \psi_1,\ldots , \psi_{n-1} \}   \cup  \{ y_1,\ldots , y_{n} \}    \cup  \{ e(\textbf{i}) \mbox{         }  | \mbox{         } \textbf{i} \in I^n \}$$
subject to the same relations as in Theorem \ref{definkl}, with the additional relations
\begin{equation} \label{relations blob2}
 e(\bs{i})=0
\end{equation}
for each $\bs{i} \in I^n $ such that
$i_1=k$ and $i_2=k-1$, or  $i_1=-k$ and $i_2=-k-1$. These relations are homogeneous with respect to the
degree function defined in Definition \ref{definkl}. Therefore,  $\blob{m}$ can be provided with
the structure of a $  \mathbb{Z}$-graded algebra
such that $ \Phi $ is a homogeneous homomorphism.
\end{cor}

\begin{dem}
The result follows by a direct application of the
Theorems \ref{isomorphism}, \ref{idempotentes homogeneous} and \ref{chevalley}.
Note that the ideal generated by a sum of orthogonal idempotents coincides with
the ideal generated by the idempotents.
\end{dem}

\begin{rem} \rm We can also give an homogeneous presentation for $Tl_n(q)$ as follows. First, note that for $\lambda=(3) \in \pa{3}$ we have
\begin{equation*}
x_{\lambda} = T_1 T_2 T_1 + T_1 T_2  + T_2 T_1 +  T_1   +  T_2  +1
\end{equation*}
On the other hand, by \cite[Corollary 4.16]{hu-mathas} if $q$ is not a cubic root of unity then we have in $\mathcal{H}_3(q^{2})$ that  $x_{\lambda}=ce(0,1,2)$, where $c\in \C^{\times}$. Therefore, we obtain a homogeneous
presentation of $Tl_{n}(q)$ by imposing to the homogeneous presentation of $\mathcal{H}_n(q^{2})$
the additional relation
\begin{equation*}
e(\bs{i}) =0
\end{equation*}
for $ \bs{i}  \in I^n $ such that $ i_{1}=0, i_{2}=1 \text{ and } i_{3}=2$.
If $q$ is a cubic root of unity, again using \cite[Corollary 4.16]{hu-mathas}, we should impose
the additional relation
\begin{equation*}
e(\bs{i})y_{3} =0
\end{equation*}
where $ i_{1}=0, i_{2}=1 \text{ and } i_{3}=2$.
\end{rem}

\section{Diagrams algebras and combinatorics of tableaux.} \label{diagrams algebra}

In this section, we briefly recall the diagram bases for the Temperley-Lieb algebra and the blob algebra,
together with the well known indexation of the
Temperley-Lieb diagrams via pairs of two-column standard tableaux.
We then go on to introduce a generalization of this to the
blob algebra diagrams, via pairs of one-line standard bitableaux $\Std(\bs{n})$ of total
degree $ n$.
In section 4.2 we endow $\Std(\bs{n})$ with a partial order structure,  different
from the usual dominance order.
In section 4.3 we
describe this partial order in terms of certain
``walks'' on the Pascal triangle.  Note that although the relevance of these
walks is indicated in section (4.6) of \cite{martin-wood}, the systematic treatment of them
seems to be new.

\subsection{Diagram basis for $Tl_n(q)$.}
We first recall the diagrammatic realization of the Temperley-Lieb algebra $Tl_n(q)$, first given by L. Kauffman,
in which the basis elements are drawn as ``$(n,n)$-bridges'' or simply ``Temperley-Lieb diagrams''.
An $(n,n)$-bridge consists of $n$ vertices, also called points or nodes, on each of two
parallel edges, the  ``top'' resp. ``bottom'' lines,
that are joined pairwise by $n$ non-intersecting lines between the two lines.
Figure 1 shows two examples.

\begin{figure}[h!]
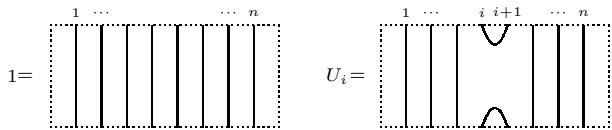

\xy 0;/r.08pc/:
(50,25)*{ \scriptscriptstyle 1};
(60,25)*{ \scriptscriptstyle \cdots };
(110,25)*{ \scriptscriptstyle \cdots };
(120,25)*{ \scriptscriptstyle n };
(50,20); (50,-20) **\dir{-};
(60,20); (60,-20) **\dir{-};
(70,20); (70,-20) **\dir{-};
(80,20); (80,-20) **\dir{-};
(90,20); (90,-20) **\dir{-};
(100,20); (100,-20) **\dir{-};
(110,20); (110,-20) **\dir{-};
(120,20); (120,-20) **\dir{-};
(40,20); (40,-20) **\dir{.};
(40,-20); (130,-20) **\dir{.};
(130,-20); (130,20) **\dir{.};
(130,20); (40,20) **\dir{.};
(25,0)*{\scriptstyle 1};
(30,0)*{\scriptstyle =};
(180,25)*{ \scriptscriptstyle 1};
(190,25)*{ \scriptscriptstyle \cdots };
(240,25)*{ \scriptscriptstyle \cdots };
(250,25)*{ \scriptscriptstyle n };
(210,25)*{ \scriptscriptstyle i };
(220,25)*{ \scriptscriptstyle i+1 };
(180,20); (180,-20) **\dir{-};
(190,20); (190,-20) **\dir{-};
(200,20); (200,-20) **\dir{-};
(210,20);(220,20)  **\crv{(215,5)};
(210,-20);(220,-20)  **\crv{(215,-5)};
(230,20); (230,-20) **\dir{-};
(240,20); (240,-20) **\dir{-};
(250,-20); (250,20) **\dir{-};
(170,20); (170,-20) **\dir{.};
(170,-20); (260,-20) **\dir{.};
(260,-20); (260,20) **\dir{.};
(260,20); (170,20) **\dir{.};
(153,0)*{ \scriptstyle U_i};
(161,0)*{\scriptstyle =};
(0,0)*{\mbox{  }};
(-45,0)*{\mbox{  }};
\endxy
  \caption{Diagrammatic generators} \label{1 y U_i tl}
\end{figure}

The set of all $(n,n)$-bridges is
denoted by $\mathbb{T}(n)$.
We define a multiplication on $ \F \mathbb{T}(n) $ by identifying the bottom of the first diagram with the top of the
second, and replacing every closed loop that may arise by a factor $-[2]$ (see Figure $2$).

\begin{figure}[h!]
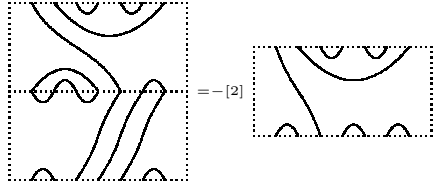
 \label{composicion temperley-lieb}
\xy 0;/r.07pc/:
(-50,0)*{\mbox{  }};
(0,0)*{\mbox{  }};
(80,80); (120,40) **\crv{(90,60) & (110,60)};
(90,80);(140,80)  **\crv{(115,50)};
(100,80);(110,80)  **\crv{(105,70)};
(120,80);(130,80)  **\crv{(125,70)};
(80,40);(110,40)  **\crv{(95,60)};
(90,40);(100,40)  **\crv{(95,50)};
(130,40);(140,40)  **\crv{(135,50)};
(70,80); (150,80) **\dir{.};
(150,80); (150,40) **\dir{.};
(150,40); (70,40) **\dir{.};
(70,40); (70,80) **\dir{.};
(70,40); (70,0) **\dir{.};
(70,0); (150,0) **\dir{.};
(150,0); (150,40) **\dir{.};
(80,40);(90,40)  **\crv{(85,30)};
(100,40);(110,40)  **\crv{(105,30)};
(80,0);(90,0)  **\crv{(85,10)};
(130,0);(140,0)  **\crv{(135,10)};
(120,40); (100,0) **\crv{(107,20) & (113,20)};
(130,40); (110,0) **\crv{(117,20) & (123,20)};
(140,40); (120,0) **\crv{(127,20) & (133,20)};
(165,40)*{\scriptscriptstyle = -[2]};
(190,60); (210,20) **\crv{(197,40) & (203,40)};
(200,60);(250,60)  **\crv{(225,30)};
(210,60);(220,60)  **\crv{(215,50)};
(230,60);(240,60)  **\crv{(235,50)};
(190,20);(200,20)  **\crv{(195,30)};
(220,20);(230,20)  **\crv{(225,30)};
(240,20);(250,20)  **\crv{(245,30)};
(180,60); (260,60) **\dir{.};
(260,60); (260,20) **\dir{.};
(260,20); (180,20) **\dir{.};
(180,20); (180,60) **\dir{.};
\endxy
  \caption{Composition in $Tl_7(q)$}
\end{figure}

With this definition $ \F \mathbb{T}(n) $ becomes a $ \C$-algebra where the identity element is the diagram
denoted by $ 1 $ in Figure \ref{1 y U_i tl}. The diagrammatic realization of the Temperley-Lieb algebra
refers to the isomorphism of $ \C $-algebras $f: Tl_n \rightarrow  \F \mathbb{T}(n) $, given by
$ f(U_i) = U_i $ where the second $ U_i $ is the diagram of Figure 1.

\medskip
Let us now recall the bijection between $(n,n)$-bridges and pairs of two-column standard tableaux
of the same shape.
Let $ \beta $ be an element of $\mathbb{T}(n)$.
We say that a line of $ \beta $ is vertical if it travels from top to bottom,
otherwise we say that it is horizontal.
Suppose now that $ \beta $
has exactly $v$ vertical lines and set $h = \frac{n-v}{2}$. The associated pair of
standard $(h+v , h)'$-tableaux $(\tau_{top}(\beta), \tau_{bot}(\beta))$ is then given by the following rules:
\begin{enumerate}
  \item $k$ is in the second column of $\tau_{top}(\beta)$  $(\tau_{bot}(\beta))$ if and only if the $k$-th point
is a right endpoint of a horizontal line in the top (bottom) edge.
  \item the numbers increase along the columns of $\tau_{top}(\beta)$ and $\tau_{bot}(\beta)$.
\end{enumerate}
For $\lambda \in \pacol{n}$ and $\sigma, \tau \in \Std(\lambda)$, we denote by $\beta_{\sigma \tau}$
the unique $(n,n)$-bridge such that $\tau_{top}(\beta_{\sigma \tau})=\sigma$ and $\tau_{bot}(\beta_{\sigma \tau})=\tau$.

\begin{exa} \rm
Let $\beta$ be the diagram to the right of Figure \ref{composicion temperley-lieb}. Then,
\begin{equation*}
\Yvcentermath1 \tau_{top}(\beta)=\young(14,26,37,5)  \qquad \qquad \tau_{bot}(\beta)=\young(12,35,47,6)
\end{equation*}
\end{exa}

\subsection{Bipartitions, bitableaux and diagrammatic realization of $\blob{m}$.}
We aim at generalizing the above results to the case of the blob algebra.
For this we first recall the concepts of bipartitions and bitableaux.
We provide them with structures of partially ordered sets, in a non-conventional way.

\medskip
A bipartition of $n$ is a pair $\bs{\lambda}=(\lambda^{(1)},\lambda^{(2)})$ of usual (integer) partitions such that
$n=|\lambda^{(1)}|+|\lambda^{(2)}|.$
By the Young diagram of a bipartition $\bs{\lambda}$ we mean the set
$$[\bs{\lambda}] =\{ (i,j,k) \in  \mathbb{N}\times
\mathbb{N} \times \{1,2\} \mbox{  } | \mbox{   } 1\leq j \leq \lambda_i^{(k)} \}.$$
Its elements are called entries or nodes.
We can visualize $[\bs{\lambda}] $ as a pair of usual Young diagrams
called the components of $ [\bs{\lambda}] $. Thus for $ d=1,2$, the $d$'th component of
$ [\bs{\lambda}] $ is
$ \{ (i,j,k) \in [\bs{\lambda}] \,  | \,  k=d\}$.
A one-line bipartition of $n$ is a
bipartition $\bs{\lambda} $ of $n$ such that $\lambda^{(k)}_i=0 $ for all $i\geq 2$ and $k=1,2$.
The set of all one-line bipartitions of $n$ is denoted $\bi{n}$.
For
$\bs{\lambda}$ a bipartition,
a $\bs{\lambda}$-bitableau is a bijection $\sub{t} : [\bs{\lambda}] \rightarrow \{ 1,\ldots,n  \}$.
We say that $\sub{t}$ has shape $\bs{\lambda} $ and write $\Shape(\sub{t})=\bs{\lambda}$.
A $\bs{\lambda}$-bitableau  is called standard if in each component its entries increase along each row and down each column.  The set of all standard $\bs{\lambda}$-bitableaux is denoted by $\Std(\bs{\lambda})$
and the union $ \bigcup_{\bs{\lambda}} \Std(\bs{\lambda})$ with $\bs{\lambda} $ running
over all bipartitions of $n$  is denoted by $\Std(\bs{n})$.

\medskip
There are several ways of endowing $\bi{n}$ with an order structure,
the most well known being dominance order, but we shall need a different order on $\bi{n}$ that
we now explain.
Let $\Lambda_n $ be the set $ \{ -n,-n+2,\ldots,n-2,n \}$.
Then the following definition makes $\Lambda_n $ into a totally ordered set with order relation $\succ$.
\begin{defi}
Suppose $\lambda, \mu \in \Lambda_n $. We then define $\mu  \succeq \lambda$
if either $| \mu |< |\lambda|$, or if
$ | \mu |= |\lambda|  $ and $ \mu \leq \lambda$.
\end{defi}
On the other hand, the map $ f $ given by
$$  f:\bi{n} \rightarrow \Lambda_n \text{,   } ((a),(b)) \rightarrow a-b $$
is a bijection and so we can define a total order $ \succeq $ on $\bi{n}$ as follows.
\begin{defi}
Suppose $\bs{\lambda},\bs{\mu} \in \bi{n}$. Then we define $\bs{\lambda} \succeq \bs{\mu}$
iff $ f(\bs{\lambda}) \succeq f(\bs{\mu})$.
\end{defi}

For $\sub{t} \in \Std(\bs{\lambda})$  let $\sub{t}_{|_k}$ be the tableau obtained from $\sub{t}$ by removing
the entries greater than $k$. We extend the order $\succeq$ to the set of all $\bs{\lambda}$-standard bitableaux as follows.

\begin{defi}\label{orden blob} \rm Suppose that $\bs{\lambda } \in \bi{n}$
and $\sub{s}, \sub{t} \in \Std(\bs{\lambda})$. We define $\sub{s} \succeq \sub{t} $ if
$\Shape(\sub{s}_{|_k}) \succeq $ $\Shape(\sub{t}_{|_k})$ for all $k=1,\ldots,n$.
\end{defi}

\begin{exa} \rm
Let $\bs{\lambda}=((6),(3))\in \bi{n}$. Define $\sub{s},\sub{t} \in \Std(\bs{\lambda})$ as follows:
\begin{equation*}
\Yvcentermath1   \sub{s}= \Bigl(\young(245689),\young(137)\Bigr)   \qquad  \qquad \sub{t}=\Bigl(\young(145679),\young(238)\Bigr)
\end{equation*}
Then, $\sub{s} \succeq \sub{t}$.
\end{exa}

Note that $\succeq$ is a partial order on $\Std(\bs{\lambda})$, but not total.
Let $\sub{t}^{\bs{\lambda}}$ be the unique standard $\bs{\lambda}$-bitableau
such that $\sub{t}^{\bs{\lambda}} \succeq \sub{t}$ for all
$\sub{t} \in \Std(\bs{\lambda})$. For  $\bs{\lambda} =(a,b)$, set $m=\min\{a,b\} $. Then
in $t^{\bs{\lambda}}$ the numbers $1,2,\ldots, n$ are located increasingly along the rows according to the following rules:
\begin{enumerate}
  \item even numbers less than or equal to $2m$ are placed in the first component.
  \item odd numbers less than $2m$ are placed in the second component.
  \item numbers greater than $2m$ are placed in the remaining boxes.
\end{enumerate}

\begin{defi} \label{secuencia blob} \rm
Suppose  that $\bs{\lambda} \in \bi{n}$ and let $\sub{t} \in \Std(\bs{\lambda})$. Define a sequence of integers inductively by the rules $\sub{t}(0)=0$ and for $1\leq j\leq n$  $$\sub{t} (j)=\sub{t}(j-1) \pm 1$$
where the $+$ ($-$) sign is used if $j$ is in the first (second) component of $\sub{t}$.
\end{defi}

Using this sequence we can now describe the order $\succeq$.

\begin{lem} \label{order with sequence}
 If $\sub{s}, \sub{t} \in \Std(\bs{\lambda})$, then $\sub{s}\succeq \sub{t}$ if and only if $|\sub{s} (j)| \leq |\sub{t}(j) |$, for all $1\leq j \leq n$, and if $|\sub{s} (j)|= |\sub{t} (j)|$ then $\sub{s} (j) \leq \sub{t} (j)$.
\end{lem}
\begin{dem}
Note that for all $\sub{t}\in \Std(\bs{\lambda}) $ and $1\leq j \leq n$, we have
$\sub{t}(j) = f(\Shape(\sub{t}|_{j}) )$. Therefore, the result is a direct consequence of Definition \ref{orden blob}.
\end{dem}

\begin{figure}[h!]
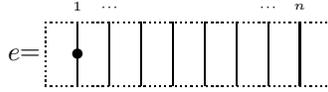

\xy 0;/r.1pc/:
(-70,0)*{\mbox{  }};
(50,25)*{ \scriptscriptstyle 1};
(60,25)*{ \scriptscriptstyle \cdots };
(110,25)*{ \scriptscriptstyle \cdots };
(120,25)*{ \scriptscriptstyle n };
(50,10)*{ \bullet};
(50,20); (50,0) **\dir{-};
(60,20); (60,0) **\dir{-};
(70,20); (70,0) **\dir{-};
(80,20); (80,0) **\dir{-};
(90,20); (90,0) **\dir{-};
(100,20); (100,0) **\dir{-};
(110,20); (110,0) **\dir{-};
(120,0); (120,20) **\dir{-};
(40,20); (40,0) **\dir{.};
(40,0); (130,0) **\dir{.};
(130,0); (130,20) **\dir{.};
(130,20); (40,20) **\dir{.};
(30,10)*{ e};
(35,10)*{ =};
(0,0)*{\mbox{  }};
\endxy
  \caption{Blob generator e.} \label{blob genera}
\end{figure}

\medskip
As is the case for the Temperley-Lieb algebra, the blob algebra has a diagrammatic realization that we now explain.
A ``blob diagram on $ n$ points'', or just a blob diagram when no confusion arises, is an $(n,n)$-bridge with possible
decorations of ``blobs'' on certain of its lines.
The blobs appear subject to the following conditions.  Each line is decorated with at most one blob; no line to the right of the
leftmost vertical line may be decorated; and to the left of it, only the outermost line in any nested formation of
loop lines can be decorated. The set of blob diagrams on $n$ points is denoted $\mathbb{B}(n)$.
Similar to the Temperley-Lieb case, there is now a
multiplication on $ \F \mathbb{B}(n) $, defined using a concatenation procedure.
This may give rise to internal loops and multiple blobs on certain lines. We then impose the rules on the multiplication
that any diagram with multiple blobs
on one or several lines is considered equal to the same diagram with a single blob on those lines,
and any internal loop is removed from the diagram
multiplying by $ y_e  = -\frac{[m-1]}{[m]}$, if the loop
is decorated, otherwise by $ -(q+q^{-1})$. The realization of
$ \blob{m} $ is
now the isomorphism $ f: \blob{m} \rightarrow \F \mathbb{B}(n) $, mapping $ U_i $ and $ e $ to the
diagrams $ U_i $ and $e$, given in Figure 1 and \ref{blob genera}.

\begin{figure}[h!]
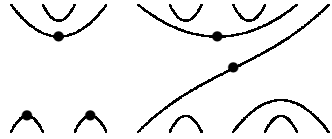

\xy 0;/r.1pc/:
(-100,0)*{\mbox{  }};
(0,40); (30,40) **\crv{(15,20)};
(10,40); (20,40) **\crv{(15,30)};
(40,40); (90,40) **\crv{(65,20)};
(50,40); (60,40) **\crv{(55,30)};
(70,40); (80,40) **\crv{(75,30)};
(0,0); (10,0) **\crv{(5,10)};
(20,0); (30,0) **\crv{(25,10)};
(50,0); (60,0) **\crv{(55,10)};
(70,0); (100,0) **\crv{(85,20)};
(80,0); (90,0) **\crv{(85,10)};
(40,0); (100,40) **\crv{(60,20) & (80,20)};
(-23,0)*{\mbox{  }};
(15,30)*{\bullet};
(65,30)*{\bullet};
(70,20)*{\bullet};
(5,5)*{\bullet};
(25,5)*{\bullet};
\endxy
  \caption{A diagram of $b_{11}(m)$} \label{blob 11}
\end{figure}

\medskip
Our next goal is to establish a bijection between
the set of blob diagrams and the set of pairs of one-line standard bitableaux of the same shape.
Let $m$ be a blob diagram. Given a horizontal line $ l $, in either edge, we put $ l = (a, b) $ where $ a $ is the left endpoint and $ b $ is the right endpoint. Let $l_1=(a_1,b_1)$ and $l_2=(a_2,b_2)$ be horizontal lines on the same edge. We say that $ l_1 $ covers $ l_2 $ if $ a_1 <a_2 <b_2 <b_1 $. We also
say that the leftmost vertical line (if any) covers all lines to the right of it. Now, we say that a node is covered if the
line to which it belongs is decorated or the line to which it belongs is covered by a decorated line. If a node is not
covered, we call it uncovered.

\begin{defi} \rm  \label{defi biy blob}
Let $m$ be a blob diagram.
Suppose that $ m $ has exactly $ v $ vertical lines and $ h = \frac{n-v}{2} $ horizontal lines on each edge.
\begin{itemize}
\item If $v\geq 0$ and the leftmost vertical line is not decorated or there is no vertical lines then
we associate to $m$ a pair of $\bs{\lambda}$-bitableaux, $\sub{t}_{top}(m)$
and $\sub{t}_{bot}(m)$, with $\bs{\lambda}= ((h+v),(h))$ by the following rules
\begin{enumerate}
\item $k$ is in the second component of $ \sub{t}_{top} (m)$ ($ \sub{t}_ {bot} (m) $) if and only if: either
$k$ is uncovered
and it is the right endpoint of a horizontal line on the top (bottom) edge,
or it is covered and it is the left endpoint of a horizontal line on the top (bottom) edge
\item the numbers increase along rows.
\end{enumerate}
  \item If $v>0$ and the leftmost vertical line is decorated then we associate to $m$ a pair of $\bs{\lambda}$-bitableaux,
$\sub{t}_{top}(m)$ and $\sub{t}_{bot}(m)$, with $\bs{\lambda}= ((h), (h+v))$ by the following rules
\begin{enumerate}
\item $k$ is in the first component of $ \sub{t}_{top} (m)$ ($ \sub{t}_ {bot} (m) $) if and only if: either
it is uncovered and it is the
left endpoint of a horizontal line on the top (bottom) edge
or it is covered and
it is the right endpoint of a horizontal line on the top (bottom) edge
  \item the numbers increase along rows.
\end{enumerate}
\end{itemize}
\end{defi}

We view these rules as a generalization of the bijection between $\mathbb{T}(n)$ and $\pacol{n}$,
with the two components of the bitableau replacing the two columns of the element of $\pacol{n}$ and with
the presence of a cover reversing the
roles of left and right.

\medskip

For $\bs{\lambda} \in \bi{n}$ and $\sub{s}, \sub{t} \in \Std(\bs{\lambda})$, we let $m_{\sub{s} \sub{t}}$
denote the unique blob diagram such that $\sub{t}_{top}(m_{\sub{s} \sub{t}})=\sub{s}$ and $\sub{t}_{bot}(m_{\sub{s} \sub{t}})=\sub{t}$.

\begin{rem}\rm \label{remark covered blob sequence}
For all $\sub{t} \in \Std(\bs{\lambda})$ and $1\leq j \leq n$, we have
\begin{description}
  \item[(\emph{i})]  If $\sub{t}(j) < 0$ then the node $j$ is covered in the top edge of $m_{\sub{tt}^{\bs{\lambda}}}$.
  \item[(\emph{ii})] If the node $j$ is covered in the top edge of $m_{\sub{tt}^{\bs{\lambda}}}$ then $\sub{t}(j)\leq 0$.
\end{description}
\end{rem}

\begin{exa} \rm   \label{ejemplo blob 11}
Let $m$ be the diagram in Figure \ref{blob 11} then
\begin{figure}[h!]
\xy 0;/r.2pc/:
(5,30); (5,25) **\dir{-};
(10,30); (10,25) **\dir{-};
(15,30); (15,25) **\dir{-};
(20,30); (20,25) **\dir{-};
(25,30); (25,25) **\dir{-};
(30,30); (30,25) **\dir{-};
(5,30); (30,30) **\dir{-};
(30,25); (5,25) **\dir{-};
(7.5,27.5)*{3};
(12.5,27.5)*{4};
(17.5,27.5)*{7};
(22.5,27.5)*{9};
(27.5,27.5)*{10};
(32.5,25)*{,};
(35,30); (35,25) **\dir{-};
(40,30); (40,25) **\dir{-};
(45,30); (45,25) **\dir{-};
(50,30); (50,25) **\dir{-};
(55,30); (55,25) **\dir{-};
(60,30); (60,25) **\dir{-};
(65,30); (65,25) **\dir{-};
(35,30); (65,30) **\dir{-};
(35,25); (65,25) **\dir{-};
(37.5,27.5)*{1};
(42.5,27.5)*{2};
(47.5,27.5)*{5};
(52.5,27.5)*{6};
(57.5,27.5)*{8};
(62.5,27.5)*{11};
(3,27.5)*{(};
(67,27.5)*{)};
(-8,27.5)*{\sub{t}_{top}(m)=};
(5,10); (5,15) **\dir{-};
(10,10); (10,15) **\dir{-};
(15,10); (15,15) **\dir{-};
(20,10); (20,15) **\dir{-};
(25,10); (25,15) **\dir{-};
(30,10); (30,15) **\dir{-};
(5,10); (30,10) **\dir{-};
(30,15); (5,15) **\dir{-};
(7.5,12.5)*{2};
(12.5,12.5)*{4};
(17.5,12.5)*{7};
(22.5,12.5)*{10};
(27.5,12.5)*{11};
(32.5,10)*{,};
(35,10); (35,15) **\dir{-};
(40,10); (40,15) **\dir{-};
(45,10); (45,15) **\dir{-};
(50,10); (50,15) **\dir{-};
(55,10); (55,15) **\dir{-};
(60,10); (60,15) **\dir{-};
(65,10); (65,15) **\dir{-};
(35,10); (65,10) **\dir{-};
(35,15); (65,15) **\dir{-};
(37.5,12.5)*{1};
(42.5,12.5)*{3};
(47.5,12.5)*{5};
(52.5,12.5)*{6};
(57.5,12.5)*{8};
(62.5,12.5)*{9};
(3,12.5)*{(};
(67,12.5)*{)};
(-8,12.5)*{\sub{t}_{bot}(m)=};
(-45,12.5)*{ };
\endxy
 \label{lalala}
\end{figure}
\end{exa}

\subsection{Walks on the Bratteli diagram.}
There is a canonical inclusion $ \blob{m} \subset b_{n+1}(m) $
which at the diagrammatic level is given by adding
to a diagram for $ \blob{m} $ a vertical undecorated line to the right, hence
producing a diagram for $ b_{n+1}(m) $ where
the two points labelled $ n+1$ are joined by a vertical line. In this way the union $ \bigcup_n b_{n}(m) $
becomes a tower of algebras and so it has an associated Bratteli diagram that
describes the generic induction and restriction rules, see for example \cite{Mat-Sal},
\cite{martin-wood} and \cite{martin-wood1}.

\medskip
We now explain how this Bratteli diagram provides a useful
interpretation of the order $ \succ $ on
$ \Std(\bs{\lambda}) $.
Let $ \mathbb{B}^{top}(n) $ (resp. $ \mathbb{B}^{bot}(n) $) denote the set of upper (lower) halves of blob diagrams. To be
more precise,
$ \mathbb{B}^{top}(n) $ (resp. $ \mathbb{B}^{bot}(n) $) consists of all
blob diagrams on $ n $ points with the information on the bottom (top) points of the vertical lines omitted.
Thus $ \mathbb{B}^{top}(n) $ (resp. $ \mathbb{B}^{bot}(n) $) is in bijection with $\Std(\bs{n})$ via $ m \mapsto
\sub{t}_{top}(m) $ (resp. $ m \mapsto \sub{t}_{bot}(m) $) and so
$ \mathbb{B}^{top}(n) $ and $ \mathbb{B}^{bot}(n) $ are in bijection with each other.
 On the diagrammatic level, the bijection can be
visualized as a reflection through an appropriate  horizontal axis.

\medskip
We know for from {\it loc. cit.} that the Bratteli diagram
gives an enumeration of $ \mathbb{B}^{top}(n) $ through a Pascal triangle pattern.
To be precise, for $ \lambda \in \Lambda_n  $ the Bratteli diagram associates with the
point $ ( \lambda, n) $ of the plane the set
$ \mathbb{B}^{top}(n, \lambda ) $, defined as those diagrams
from $ \mathbb{B}^{top}(n) $ that have exactly $| \lambda |$ vertical lines,
where the leftmost vertical line is decorated iff $ \lambda $
is negative. Set $ b_{n, \lambda } := | \mathbb{B}^{top}(n, \lambda) |$ with the convention that
$ \mathbb{B}^{top}(n, \lambda) := \emptyset $
if $ \lambda \notin \Lambda_n$.
Then there is a bijection between $  \mathbb{B}^{top}(n, \lambda ) $ and
$  \mathbb{B}^{top}(n-1, \lambda+1 ) \cup \mathbb{B}^{top}(n-1, \lambda -1) $, as we explain shortly. The
Pascal triangle formula $ b_{n, \lambda } = b_{n-1, \lambda+1 } + b_{n-1, \lambda-1 } $ is a consequence of
this bijection.

\begin{figure}[h!]
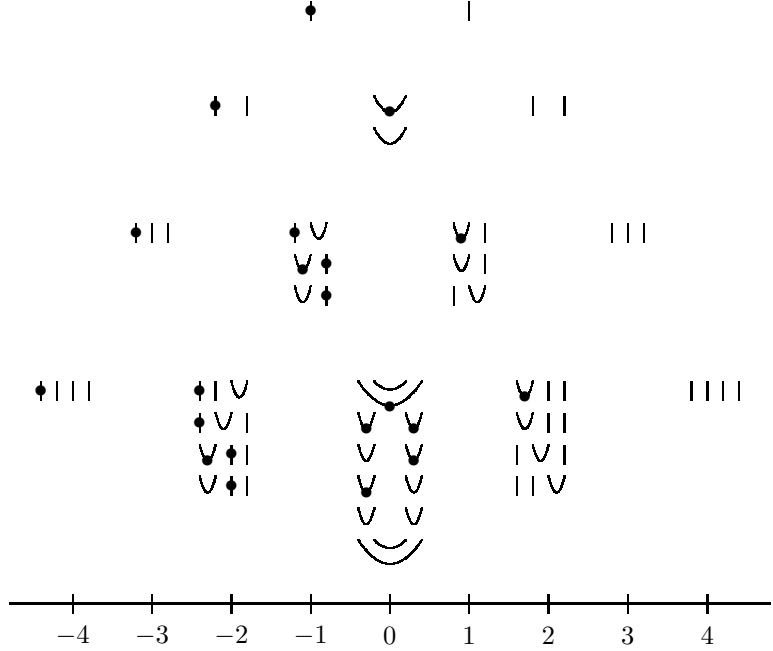

\[
\xy 0;/r1pc/:
(-12,0); (12,0) **\dir{-};
(-10,-0.3); (-10,0.3) **\dir{-};
(-5,-0.3); (-5,0.3) **\dir{-};
(0,-0.3); (0,0.3) **\dir{-};
(5,-0.3); (5,0.3) **\dir{-};
(10,-0.3); (10,0.3) **\dir{-};
(-7.5,-0.3); (-7.5,0.3) **\dir{-};
(-2.5,-0.3); (-2.5,0.3) **\dir{-};
(7.5,-0.3); (7.5,0.3) **\dir{-};
(2.5,-0.3); (2.5,0.3) **\dir{-};
(0,-1)*{0};
(-2.5,-1)*{-1};
(-5,-1)*{-2};
(-7.5,-1)*{-3};
(-10,-1)*{-4};
(2.5,-1)*{1};
(5,-1)*{2};
(7.5,-1)*{3};
(10,-1)*{4};
(-1,2); (1,2) **\crv{(0,.5)};
(-0.5,2); (0.5,2) **\crv{(0,1.5)};
(-1,3); (-0.5,3) **\crv{(-0.75,2)};
(1,3); (0.5,3) **\crv{(0.75,2)};
(-1,4); (-0.5,4) **\crv{(-0.75,3)};
(1,4); (0.5,4) **\crv{(0.75,3)};
(-1,5); (-0.5,5) **\crv{(-0.75,4)};
(1,5); (0.5,5) **\crv{(0.75,4)};
(-1,6); (-0.5,6) **\crv{(-0.75,5)};
(1,6); (0.5,6) **\crv{(0.75,5)};
(-1,7); (1,7) **\crv{(0,5.5)};
(-0.5,7); (0.5,7) **\crv{(0,6.5)};
(-6,7); (-6,6.4) **\dir{-};
(-5.5,7); (-5.5,6.4) **\dir{-};
(-5,7); (-4.5,7) **\crv{(-4.75,6)};
(-6,6); (-6,5.4) **\dir{-};
(-4.5,6); (-4.5,5.4) **\dir{-};
(-5.5,6); (-5,6) **\crv{(-5.25,5)};
(-6,5); (-5.5,5) **\crv{(-5.75,4)};
(-5,5); (-5,4.4) **\dir{-};
(-4.5,5); (-4.5,4.4) **\dir{-};
(-6,4); (-5.5,4) **\crv{(-5.75,3)};
(-5,4); (-5,3.4) **\dir{-};
(-4.5,4); (-4.5,3.4) **\dir{-};
(-11,7); (-11,6.4) **\dir{-};
(-10.5,7); (-10.5,6.4) **\dir{-};
(-10,7); (-10,6.4) **\dir{-};
(-9.5,7); (-9.5,6.4) **\dir{-};
(4,7); (4.5,7) **\crv{(4.25,6)};
(5,7); (5,6.4) **\dir{-};
(5.5,7); (5.5,6.4) **\dir{-};
(4,6); (4.5,6) **\crv{(4.25,5)};
(5,6); (5,5.4) **\dir{-};
(5.5,6); (5.5,5.4) **\dir{-};
(4,5); (4,4.4) **\dir{-};
(4.5,5); (5,5) **\crv{(4.75,4)};
(5.5,5); (5.5,4.4) **\dir{-};
(4,4); (4,3.4) **\dir{-};
(5,4); (5.5,4) **\crv{(5.25,3)};
(4.5,4); (4.5,3.4) **\dir{-};
(11,7); (11,6.4) **\dir{-};
(10.5,7); (10.5,6.4) **\dir{-};
(10,7); (10,6.4) **\dir{-};
(9.5,7); (9.5,6.4) **\dir{-};
(-8,12); (-8,11.4) **\dir{-};
(-7.5,12); (-7.5,11.4) **\dir{-};
(-7,12); (-7,11.4) **\dir{-};
(-3,12); (-3,11.4) **\dir{-};
(-2.5,12); (-2,12) **\crv{(-2.25,11)};
(-3,11); (-2.5,11) **\crv{(-2.75,10)};
(-2,11); (-2,10.4) **\dir{-};
(-3,10); (-2.5,10) **\crv{(-2.75,9)};
(-2,10); (-2,9.4) **\dir{-};
(3,12); (3,11.4) **\dir{-};
(2.5,12); (2,12) **\crv{(2.25,11)};
(2,11); (2.5,11) **\crv{(2.25,10)};
(3,11); (3,10.4) **\dir{-};
(3,10); (2.5,10) **\crv{(2.75,9)};
(2,10); (2,9.4) **\dir{-};
(8,12); (8,11.4) **\dir{-};
(7.5,12); (7.5,11.4) **\dir{-};
(7,12); (7,11.4) **\dir{-};
(-5.5,16); (-5.5,15.4) **\dir{-};
(-4.5,16); (-4.5,15.4) **\dir{-};
(-0.5,16); (0.5,16) **\crv{(0,15)};
(-0.5,15); (0.5,15) **\crv{(0,14)};
(5.5,16); (5.5,15.4) **\dir{-};
(4.5,16); (4.5,15.4) **\dir{-};
(-2.5,19); (-2.5,18.4) **\dir{-};
(2.5,19); (2.5,18.4) **\dir{-};
(-2.5,18.7)*{\bullet};
(-5.5,15.7)*{\bullet};
(-8,11.7)*{\bullet};
(-6,6.7)*{\bullet};
(-6,5.7)*{\bullet};
(-5,4.7)*{\bullet};
(-5,3.7)*{\bullet};
(-3,11.7)*{\bullet};
(-2,10.7)*{\bullet};
(-2,9.7)*{\bullet};
(-11,6.7)*{\bullet};
(0,15.5)*{\bullet};
(-2.75,10.5)*{\bullet};
(-5.75,4.5)*{\bullet};
(-0.75,5.5)*{\bullet};
(0.75,5.5)*{\bullet};
(0.75,4.5)*{\bullet};
(-0.75,3.5)*{\bullet};
(2.25,11.5)*{\bullet};
(4.25,6.5)*{\bullet};
(0,6.2)*{\bullet};
\endxy
\]
\caption{Bratteli diagram} \label{Bratteli diagram}
\end{figure}

\medskip
For $ \lambda \in \Lambda_n \setminus \{ 0 \} $ define $ \lambda^+ \in \Lambda_{n+1} $ by
$ \lambda^+ := \lambda \pm 1 $ where the sign is positive iff $ \lambda > 0 $. Similarly, for
$ \lambda \in \Lambda_n \setminus \{0\} $ define
$ \lambda^- := \lambda \pm 1 $ where the sign is positive iff $ \lambda < 0$. Finally, if
$ \lambda = 0 \in \Lambda_n $ define $ \lambda^+ := 1 $ and $ \lambda^- := -1 $. With these
definitions we have for any $ \lambda
\in \Lambda_n $ that $ \lambda^- \succ \lambda^+ $ in $ \Lambda_{n+1}$.
In other words, the map $ \lambda \mapsto \lambda^- $ moves $ \lambda $ closer to the central axis of
the Bratteli diagram consisting of the points $ \{ (0, k) , k = 0,1 \ldots \} $, whereas
$ \lambda \mapsto \lambda^+ $ takes $ \lambda $ away from the central axis.

\medskip
The above mentioned
bijection is now induced by injective maps
\begin{equation} \label{bijection_par}
f^+_{n, \lambda}: \mathbb{B}^{top}(n-1, \lambda )  \rightarrow \mathbb{B}^{top}(n, \lambda^+ ),  \, \,
f^-_{n, \lambda}: \mathbb{B}^{top}(n-1, \lambda )  \rightarrow \mathbb{B}^{top}(n, \lambda^{-} )
\end{equation}
that can be described concretely as follows. If
$ m \in \mathbb{B}^{top}(n-1, \lambda )  $ then $ f^+_{n, \lambda} $ adds an undecorated vertical line on the right hand
side of $m$. If $ \lambda \neq  0 $ then
$ f^-_{n, \lambda} $ joins the rightmost vertical line of $m$ with the new $n$-th point of the (top) edge
whereas $ f^-_{n, 0} $ adds a decorated vertical line on the right hand
side of $m$.
Finally, by convention $ f^+_{1,0} $ (resp. $ f^-_{1,0} $) maps the empty diagram to the unique
diagram of $ \mathbb{B}^{top}(1, 1 )$ (resp. $ \mathbb{B}^{top}(1,- 1 )$).

\medskip
For us the main point of this construction is
that any element of $ m \in \mathbb{B}^{top}(n) $ can be written uniquely as
\begin{equation} \label{sign-sequence}
m=    f^{\sigma_{n}}_{n, \lambda_{n}}  \ldots f^{\sigma_1}_{1,0} \, \emptyset
\mbox{ where } \sigma_k \in \{+,- \} \mbox{ for } k =1, \ldots, n.
\end{equation}
In other words, the sequence of signs $ \{ \sigma_k \}_{k=1, \ldots, n} $ uniquely determines $ m$
and hence $ \mathbb{B}^{top}(n) $ is in bijection with walks on the Bratteli diagram, starting with
the empty partition in position $(0,0)$ and
at the $k$'th step, where the walk is situated in $ (k, \lambda_{k} )$, going inwards or outwards
according to the value of $ \sigma_k$.
We denote by $ W(m) $ the walk associated with $ m \in  \mathbb{B}^{top}(n) $.

\medskip
Let us now return to the order $ \succeq $ on $ \Std(\bs{\lambda}) $ introduced above.
Suppose that $\sub{s} \in \Std(\bs{\lambda})$ for $\bs{\lambda} \in \bi{n}$.
Then $\sub{s}  $ also gives rise
to a walk, denoted $ w(\sub{s})$, on the points of the Bratteli diagram. It starts in $ (0,0) $ and
for $ k =0,  1, \ldots, n-1$ goes from $ (k,j) $ to $ (k+1,j-1) $ if $ k+1 $ is
located in the second component of $\sub{s} $ and to $ (k+1,j+1) $ if $ k+1 $ is located in the first component of
$\sub{s}  $. In other words, at the $k$'th step the walk $ w(\sub{s}) $ is situated in $ (k, \sub{s}(k)) $ where
$ \{ \sub{s}(k) \,  |\, k = 0, 1, \ldots, n \}  $ is the sequence of integers associated with $\sub{s} $
as in Definition \ref{secuencia blob}.
With this walk realization of the bitableaux, we can visualize the order $ \succeq$. Indeed, let
$\sub{s}, \sub{t} \in \Std(\bs{\lambda}) $.
Then $\sub{s} \succeq \sub{t} $ iff at each step of the two walks $ w(\sub{s}) $ is
either strictly closer than $ w(\sub{t}) $ to the central vertical axis of the Bratteli diagram or
they are at the same distance from the central axis and $ w(\sub{s}) $ is located (weakly) to the left of $ w(\sub{t}) $.

\medskip
Let us now explain the relationship between the two walks.
We denote by $ s $ the bijection $ \mathbb{B}^{top}(n) \rightarrow \Std(\bs{n}), \,  m \mapsto
\sub{t}_{top}(m) $, mentioned above.

\begin{lem} \label{lemma-sequence-sign}
Let $ m \in \mathbb{B}^{top}(n) $. Then we have $ W(m) = w(s(m))$.
 \end{lem}
\begin{dem}
This is a consequence of Remark \ref{remark covered blob sequence} and the definitions.
\end{dem}

\medskip
There is a natural surjective map $\pi: \mathbb{B}(n) \rightarrow \mathbb{T}(n)$, which sends a blob
diagram $m$ to the $(n,n)$-bridge obtained by deleting all decorations in $m$.
On the other hand,
$\mathbb{T}(n)$ is in bijection with pairs of two-column standard tableaux of the same shape and
$\mathbb{B}(n)$ is in bijection
with pairs of one-line standard bitableaux of the same shape by Definition \ref{defi biy blob},
and so our next goal is to describe the above map $ \pi $
in terms of one-line bitableaux and
two-column tableaux. For this we make a couple of definitions.

\begin{defi} \rm \label{tipo tl}
Suppose that $\bs{\lambda}=((a),(b))\in \bi{n}$
and let $\sub{t} \in \Std(\bs{\lambda})$. Set $\mu_{1}=\max \{a,b\}$ and $\mu_{2}=\min \{a,b\}$. Let $\mu$
be the two-column partition of $n$ given by $\mu= (\mu_{1},\mu_{2})'$.
Then we define $\tc{t}$ as the unique $\mu$-standard tableau that satisfies
\begin{center}
  $k$ \emph{is in the second column of} $\tc{t}$ \emph{if and only if} $| \sub{t}(k) | <| \sub{t}(k-1) |. $
  \end{center}
\end{defi}

We claim that $\tc{t}$ defined in this way is a standard tableau. For this we use that
a node $k$ of the blob diagram given by $m_{\sub{st}}$
is a right endpoint in the top (resp. bottom) edge if and
only if $ |\sub{s}(k)| <  |\sub{s}(k-1)|$ (resp. $ |\sub{t}(k)| <  |\sub{t}(k-1)|$), as
can easily be seen by analysing Definition \ref{defi biy blob}.
In other words,  $\tc{s} $ and $ \tc{t}$ can be described as the unique two-column tableaux that
satisfy $\pi(m_{\sub{st}})  = \beta_{\tc{s} \tc{t}}$, where
$ \pi $ is the map defined above,
and our claim follows.

\medskip
For $\sub{s} \in \Std(\bs{\lambda})$ we let $ | w(\sub{s}) |$ denote the walk on the Bratteli diagram
that at the $ k$'th step is located in the point $ (k, | \sub{s}(k)| ) $. The two components of its associated bitableau
are then the conjugates of the columns of $\tc{s}$, as follows from the above.

\begin{defi} \rm \label{tilde}
For $\bs{\lambda}\in \bi{n}$ and $\sub{s}, \sub{t} \in \Std(\bs{\lambda})$,
we write $\sub{s} \sim \sub{t}$
if $\tc{s} = \tc{t}$. Thus $\sub{s} \sim \sub{t}$ if and only if $|\sub{s}(k)|= |\sub{t}(k)|$ for all $1\leq k \leq n$.
\end{defi}

We give a couple of Lemmas related to these definitions.
\begin{lem} \label{orden con diagrmas}
Suppose that $\bs{\lambda} \in \bi{n}$ and let $\sub{s},\sub{t} \in \Std(\bs{\lambda})$. Then, $\tc{s} \trianglerighteq \tc{t}$
if and only if $|\sub{s} (k)| \leq | \sub{t}(k)| $ for all $1\leq k \leq n$. In particular, if $\sub{s} \succeq \sub{t}$ then $\tc{s} \trianglerighteq \tc{t}$.
\end{lem}
\begin{dem}
Notice that
$$  \Shape(\tc{s}|_{k})=\left(  \frac{k+|\sub{s}(k)|}{2} ,   \frac{k-|\sub{s}(k)|}{2} \right)'  $$ $$\Shape(\tc{t}|_{k})=\left(  \frac{k+|\sub{t}(k)|}{2} ,   \frac{k-|\sub{t}(k)|}{2} \right)'  $$
for all $ 1\leq k \leq n $. Using the property of the usual dominance order that
$\mu \trianglerighteq \nu \Longleftrightarrow \nu' \trianglerighteq \mu'$ we deduce
that $\tc{s}\trianglerighteq \tc{t}$ if and only if $|\sub{s}(k) |  \leq | \sub{t}(k)|$ for all $ 1\leq k \leq n $,
which is the first claim of the Lemma. The second claim
follows now from
Lemma \ref{order with sequence}.
\end{dem}

\medskip
Using the natural embedding
$\iota: \mathbb{T}(n) \rightarrow \mathbb{B}(n)$ we obtain a walk description of the elements of $ \mathbb{T}(n) $ as
well. Under this description, $ \mathbb{T}(n) $ corresponds to the walks on the Bratteli diagram for $ \blob{m}$
that always stay in the positive half of the Bratteli diagram, including the central vertical axis.

\medskip
The left action of ${\mathfrak{S}}_n$ on tableaux generalizes to an action of ${\mathfrak{S}}_n$ on
bitableaux. Using it we have the following Lemma.
\begin{lem} \label{orden tau}
Suppose that $\bs{\lambda} \in \bi{n}$ and let $\sub{s},\sub{t} \in \Std(\bs{\lambda})$.
Suppose moreover that $ \sub{s} \succ \sub{t} $, that $ s_k \sub{s} = \sub{t} $
for some simple transposition $s_k=(k,k+1) \in {\mathfrak{S}}_n$ and that $ \sub{s} \nsim \sub{t} $.
Then $s_k \tc{s} =  \tc{t}  $ and $ \tc{s}  \vartriangleright \tc{t}  $.
\end{lem}
\begin{dem}
Note first that by the assumptions we have $ \sub{s}(j ) = \sub{t}(j ) $ for all $ j \not= k $.
Let us first assume that $ \sub{s}(k+1 )  \geq 1 $.
Then $ \sub{s}(k ) \geq 0 $ since
$ \sub{s}(j )   $ changes by $ \pm 1 $ when $ j $ is increased by $1$. But similarly $ \sub{t}(k ) \geq 0 $
and
then we must have $ \sub{t}(k ) = \sub{s}(k ) +2 $ since $ \sub{s} \succ \sub{t} $.
Since $ k $ and $ k+1 $ are located in different components in $ \sub{s} $ and in $ \sub{t} $, this gives us the equalities
$$ \sub{s}(k-1 ) = \sub{t}(k-1 ) = \sub{s}(k ) +1 = \sub{t}(k )-1 = \sub{s}(k+1 ) = \sub{t}(k+1 ) $$ from which we get by Definition
\ref{tipo tl} that $ k $ (resp. $k+1$)  is located in second (resp. first) column of $ \tc{s} $ whereas
$ k $ (resp. $k+1$)  is placed in first (resp. second) column of $ \tc{t} $. Since
$ j $ is located in the same column of $ \tc{s} $ and $ \tc{t}$ for $ j \neq k, k+1$ we
now conclude that
$s_k \tc{s} =  \tc{t}  $ and $ \tc{s}  \vartriangleright \tc{t}  $, as needed.

\medskip
The case $ \sub{s}(k+1)  \leq -1 $ is treated similarly and so the only remaining case is $ \sub{s}(k+1 )  = 0 $.
Then $ \sub{t}(k+1)  = \sub{s}(k-1 )=\sub{t}(k-1 ) = 0 $. Moreover since $ \sub{s} \succ \sub{t} $ we
have $ \sub{s}(k )  =-1 $ and $ \sub{t}(k )  =1 $. But this implies that $ \sub{s} \sim \sub{t} $, finishing the proof.
\end{dem}

\begin{defi} \rm \label{langle-rangle}
Suppose that $\bs{\lambda} \in \bi{n}$ and that $\sub{s}, \sub{t} \in \Std(\bs{\lambda})$.
Then we say that ``$\sub{s}$ has a hook at position $ k $'' if
$ \sub{s}(k-1) =   \sub{s}(k+1)  =    \sub{s}(k) \pm 1  $ where $ 1 \leq k \leq n-1 $.
Moreover we say that ``$ \sub{t}  $ is obtained from $ \sub{s}  $ by making a hook at position
$ k $ smaller'' if
$ \sub{s}(j) =   \sub{t}(j) $ for $ j \neq k $, $ \sub{s}(k) =   \sub{t}(k) \pm 2 $ and
$ \sub{s} \succ  \sub{t} $.
\end{defi}
The last condition can also be written as $s_k \sub{s} =  \sub{t}  $ and $ \sub{s} \succ  \sub{t} $.
Geometrically, if $ \sub{t}  $ is obtained from $ \sub{s}  $ by making a hook at position
$ k $ smaller then $ \sub{t}  $ is obtained from $ \sub{s}  $ by either replacing a configuration of three consecutive points
in $ w( \sub{t})  $ forming a ``$ \langle $''
by a configuration ``$ \rangle  $'' at these three points, or reversely, depending on which side of the Bratteli diagram
the configuration is located.

\begin{lem} \label{reduced exp}
For $\sub{t} \in \Std(\bs{\lambda}) $ we define $\de{t}$ as the element of ${\mathfrak{S}}_n$ that satisfies
$\sub{t} = \de{t} \sub{t}^{\bs{\lambda}}$.  Then $ \de{t} $ can be written as product of
simple transpositions $\de{t} = s_{i_k} s_{i_{k-1}}  \ldots s_{i_1}  $ such that
$ s_{i_{j}} \ldots s_{i_1} \sub{t}^{\bs{\lambda}} $ is standard and such that
$ s_{i_j} s_{i_{j-1}} \ldots s_{i_1} \sub{t}^{\bs{\lambda}} \prec
s_{i_{j-1}} \ldots s_{i_1} \sub{t}^{\bs{\lambda}} $  for all $1\leq j \leq k$.
\end{lem}
\begin{dem}
This can be seen via the walk realization of $ \Std(\bs{\lambda}) $. Indeed the walk $ w(\sub{t}^{\bs{\lambda}})$
first zigzags on and off the central vertical line of the Bratteli diagram,
using the sign $ - $ an even number of times, and then finishes using the sign $+$ repeatedly, if $\bs{\lambda}$ is located in the positive half, or using once the sign $-$ followed by the sign $ + $ repeatedly, if $\bs{\lambda}$ is located in the negative half. Figure $\ref{ejemplo walk mayor}$ shows an example of such walks.

\begin{figure}[h!]
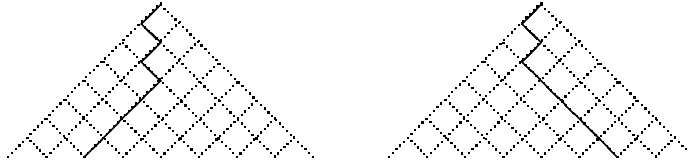

\[
\xy 0;/r.6pc/:
(0,0); (-8,-8) **\dir{.};
(0,0); (8,-8) **\dir{.};
(1,-1); (-6,-8) **\dir{.};
(2,-2); (-4,-8) **\dir{.};
(3,-3); (-2,-8) **\dir{.};
(4,-4); (0,-8) **\dir{.};
(5,-5); (2,-8) **\dir{.};
(6,-6); (4,-8) **\dir{.};
(7,-7); (6,-8) **\dir{.};
(-1,-1); (6,-8) **\dir{.};
(-2,-2); (4,-8) **\dir{.};
(-3,-3); (2,-8) **\dir{.};
(-4,-4); (0,-8) **\dir{.};
(-5,-5); (-2,-8) **\dir{.};
(-6,-6); (-4,-8) **\dir{.};
(-7,-7); (-6,-8) **\dir{.};
(0,0); (-1,-1) **\dir{-};
(0,-2); (-1,-1) **\dir{-};
(0,-2); (-1,-3) **\dir{-};
(0,-4); (-1,-3) **\dir{-};
(0,-4); (-4,-8) **\dir{-};
(0,0); (-1,-1) **\dir{-};
(0,-2); (-1,-1) **\dir{-};
(0,-2); (-1,-3) **\dir{-};
(0,-4); (-1,-3) **\dir{-};
(0,-4); (-4,-8) **\dir{-};
(0,0); (-1,-1) **\dir{-};
(0,-2); (-1,-1) **\dir{-};
(0,-2); (-1,-3) **\dir{-};
(0,-4); (-1,-3) **\dir{-};
(0,-4); (-4,-8) **\dir{-};
(20,0); (12,-8) **\dir{.};
(20,0); (28,-8) **\dir{.};
(21,-1); (14,-8) **\dir{.};
(22,-2); (16,-8) **\dir{.};
(23,-3); (18,-8) **\dir{.};
(24,-4); (20,-8) **\dir{.};
(25,-5); (22,-8) **\dir{.};
(26,-6); (24,-8) **\dir{.};
(27,-7); (26,-8) **\dir{.};
(19,-1); (26,-8) **\dir{.};
(18,-2); (24,-8) **\dir{.};
(17,-3); (22,-8) **\dir{.};
(16,-4); (20,-8) **\dir{.};
(15,-5); (18,-8) **\dir{.};
(14,-6); (16,-8) **\dir{.};
(13,-7); (14,-8) **\dir{.};
(20,0); (19,-1) **\dir{-};
(20,-2); (19,-1) **\dir{-};
(20,-2); (19,-3) **\dir{-};
(24,-8); (19,-3) **\dir{-};
(20,0); (19,-1) **\dir{-};
(20,-2); (19,-1) **\dir{-};
(20,-2); (19,-3) **\dir{-};
(24,-8); (19,-3) **\dir{-};
(20,0); (19,-1) **\dir{-};
(20,-2); (19,-1) **\dir{-};
(20,-2); (19,-3) **\dir{-};
(24,-8); (19,-3) **\dir{-};
\endxy
\]
 \caption{ Walk  $w(\sub{t}^{\bs{\lambda}}) $ for $\bs{\lambda} =((2),(6))$ (left) and $\bs{\lambda} =((6),(2))$ (right).} \label{ejemplo walk mayor}
\end{figure}

This walk can be converted into $ w(\sub{t}) $ through a series of
$ k $ transformations, say, where at each step the new walk is obtained from the previous one by making a hook
at position $j$ smaller, for some $j$.
At tableau level, each of
these transformations is given by the action of a simple transposition $ s_j $. The Lemma follows from this.
\end{dem}

\begin{exa}  \rm  We illustrate the above Lemma. Let $\bs{\lambda} = ((4),(2)) \in \bi{6}$ and
\begin{equation*}
  \Yvcentermath1 \sub{t}= \Bigl( \young(1236),\young(45)  \Bigr)
\end{equation*}
Then, we have $\de{t}=s_3s_4s_2s_3s_1$. Now, define the bitableaux $\sub{t}_0, \sub{t}_1,\sub{t}_2,\sub{t}_3,\sub{t}_4 $ and  $\sub{t}_5$ as follows:
$$
  \Yvcentermath1
  \begin{array}{ccc}
    \sub{t}_0= \Bigl( \young(2456),\young(13)  \Bigr) &  & \sub{t}_1= \Bigl( \young(1456),\young(23)  \Bigr) \\
     \sub{t}_2= \Bigl( \young(1356),\young(24)  \Bigr) &  &   \sub{t}_3= \Bigl( \young(1256),\young(34)  \Bigr) \\
    \sub{t}_4= \Bigl( \young(1246),\young(35)  \Bigr) &  & \sub{t}_5= \Bigl( \young(1236),\young(45)  \Bigr)
  \end{array}
$$
It is straightforward to check that $\sub{t}^{\bs{\lambda}}=\sub{t}_0 \succ \sub{t}_1 \succ \sub{t}_2 \succ \sub{t}_3 \succ \sub{t}_4 \succ \sub{t}_5 =\sub{t}$. The figures below show how the walk $w(\sub{t}^{\bs{\lambda}})$ is converted into $w(\sub{t})$.
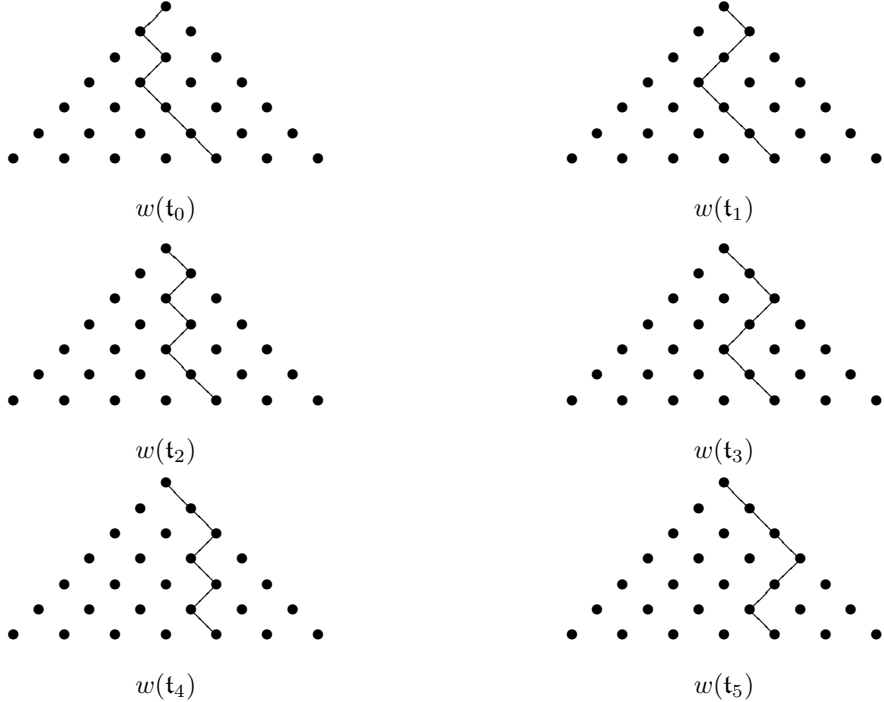
$$ \ejemploprimero  $$
$$ \ejemplosegundo  $$
$$ \ejemplotercero  $$
\end{exa}

\medskip
\begin{rem}\rm \label{remark reduced}
Although we do not need it directly, we note that $l(\de{t}) =k$
and that the expression $\de{t} = s_{i_k} s_{i_{k-1}}  \ldots s_{i_1}  $ is reduced.
\end{rem}


\section{Jucys-Murphy elements for $\blob{m}$.}
In Corollary \ref{main}, we gave a new (graded) presentation for $\blob{m}$, while in the
previous section we described the diagram basis for the blob algebra.
Unfortunately, it seems nontrivial to express the homogeneous  generators in terms of the
diagram basis of $\blob{m}$. However, it turns out that a graded cellular
basis for $\blob{m}$ can be constructed from a
precise description of the
KLR idempotents in $ \blob{m}$. Inspired by the
work of Hu and Mathas \cite{hu-mathas}, we shall obtain in the next section an expression for them
building on the results from \cite{Mat-So}. A key point for this is to make $ \blob{m} $ fit into the general setting
of an algebra with Jucys-Murphy elements. This is the main goal of this section.

\subsection{Cellular algebras.}

Before defining the concept of an algebra with Jucys-Murphy elements,
we first recall the definition of a cellular algebra, which was first given by
Graham and Lehrer in \cite{gra ler} in order to provide a common framework for a series of algebras that appear in non-semisimple representation theory.

\begin{defi} \label{cellular algebra} \rm
Let $ \R$ be an integral domain.
Suppose that $A$ is an $\R$-algebra which is free of finite rank over $\R$. Suppose that $(\Lambda, \geq)$ is a poset and that for each $\lambda \in \Lambda $ there is a finite set $T(\lambda)$ and elements $c_{\sub{st}}^{\lambda} \in A$ such that $$ \mathcal{C} =\{ c_{\sub{st}}^{\lambda} \mbox{  } | \mbox{  } \lambda\in \Lambda ; \sub{s},\sub{t} \in T(\lambda)  \}$$
is a basis of $A$. The pair $(\mathcal{C}, \Lambda)$ is a cellular basis of $A$ if
\begin{description}
  \item[(i)] The $\R$-linear map $*:A \rightarrow A$ determined by $(c_{\sub{st}}^{\lambda})^{*}=c_{\sub{ts}}^{\lambda}$, for all $\lambda \in \Lambda$ and all $\sub{s} ,\sub{t}  \in T(\lambda)$, is an algebra anti-automorphism of $A$
  \item[(ii)] If $\sub{s},\sub{t} \in T(\lambda) $, for some $\lambda \in \Lambda$, and $a\in A$ then there exist scalars $r_{\sub{u}} \in \R$ such that  $$  ac_{\sub{st}}^{\lambda} \equiv \sum_{\sub{u} \in T(\lambda)} r_{\sub{u}}   c_{\sub{ut}}^{\lambda}  \mod A^{\lambda}$$
      where $A^{\lambda}$ is the $\R$-submodule of $A$ spanned by $\{ c_{\sub{ab}}^{\mu} \mbox{  } | \mbox{  } \mu > \lambda ; \sub{a},\sub{b} \in T(\mu)  \}$
\end{description}
If $A$ has a cellular basis we say that $A$ is a cellular algebra.
\end{defi}
Note that $ r_{\sub{u}}$ depends on $\sub{u}$, $\sub{s}$ and $a$, what is important is that $ r_{\sub{u}}$ does not depend on $\sub{t}$. Now suppose that $A$ is a $\mathbb{Z}$-graded $\R$-algebra  and that each $c_{\sub{st}}^{\lambda}$ is homogeneous. If there exists a function $$   \deg: \coprod_{\lambda \in \Lambda} T(\lambda) \rightarrow \mathbb{Z}$$ such that $\deg c_{\sub{st}}^{\lambda} = \deg \sub{s} + \deg \sub{t}$, for $\lambda \in \Lambda$ and all $\sub{s},\sub{t} \in T(\lambda)$, we say that $ \mathcal{C} =\{ c_{\sub{st}}^{\lambda} \mbox{  } | \mbox{  } \lambda\in \Lambda ; \sub{s},\sub{t} \in T(\lambda)  \}$ is a graded cellular basis of $A$. If $A$ has a graded cellular basis we say that $A$ is a graded cellular algebra.

\begin{defi} \rm
Suppose that $A$ is a cellular algebra with cellular basis $ (\mathcal{C},\Lambda) $ as in the Definition \ref{cellular algebra} and fix $\lambda \in \Lambda$. Then the cell module $C^{\lambda}$ is the left $A$-module which is free as an $R$-module with basis $\{ c_{\sub{t}} \mbox{  } | \mbox{  } \sub{t} \in T(\lambda)  \}$ and where the action $A$ on $C^{\lambda}$ is given by
$$  ac_{\sub{t}}^{\lambda} = \sum_{\sub{u} \in T(\lambda)} r_{\sub{u}}   c_{\sub{u}}^{\lambda}  $$
where $r_{\sub{u}}$ is the element of $R$  that appears in the Definition \ref{cellular algebra} (ii).
\end{defi}

If in the above definition we assume that A is a graded cellular algebra, then the cell modules are graded modules.
In this case we have the direct sum decomposition
$$   C^{\lambda} =\bigoplus_{z\in \mathbb{Z}}  C_{z}^{\lambda}$$ where  $C_{z}^{\lambda} $ is the free $R$-module with basis $\{  c_{\sub{t}}^{\lambda} \mbox{   }  | \mbox{  } \sub{t}\in T(\lambda) \text{ and } \deg \sub{t} =z   \}$.

\medskip
For ${\lambda} \in \Lambda$ there is a symmetric and associative bilinear form $\langle \, \cdot  , \cdot \,\rangle_{{\lambda}}$
on $C^{{\lambda}}$ determined by
$$  c_{\sub{as}}^{\lambda}c_{\sub{tb}}^{\lambda} \equiv \langle c_{\sub{s}}^{\lambda}, c_{\sub{t}}^{\lambda} \rangle_{\lambda} c_{\sub{ab}}^{\lambda} \mod{A^{\lambda}} $$
for all $\sub{a,b,s,t} \in T(\lambda)$. Define the radical of $C^{\lambda}$ by
$$   \text{rad } C^{\lambda} = \{  x\in C^{\lambda} \, | \, \langle x,y \rangle_{\lambda} =0 \text{ for all } y\in C^{\lambda} \}        $$
It is easy to see that $\text{rad } C^{\lambda}$ is an $A$-submodule of $C^{\lambda}$ (See \cite[Proposition 2.9]{Mathas}),
and that also $\text{rad }C^{\lambda}  $ is a graded submodule of $C^{\lambda}$ (See \cite[Lemma 2.7]{hu-mathas}) . Define $D^{\lambda} = C^{\lambda} / \text{rad } C^{\lambda}$. By definition $D^{\lambda}$ is a graded $A$-module.

\medskip
In general, for a cellular algebra $ A $ we define
$\Lambda_0 = \{   \lambda \in \Lambda \mbox{  } | \mbox{  } D^{\lambda} \neq 0 \}$.
The following theorem classifies the simple $A$-modules over a field.
\begin{teo} (\cite[Theorem 3.4]{gra ler})
Suppose that R is a field. Then $\{  D^{\lambda}  \mbox{ }| \mbox{ } \lambda \in \Lambda_{0} \}$ is a complete set of pairwise non-isomorphic simple $A$-modules.
\end{teo}

\medskip
It is well known that the Temperley-Lieb algebra (\cite[Example 1.4]{gra ler}) and
the blob algebra (\cite[page 7]{gra ler 1}) are cellular, in fact the diagram bases introduced in the above section
are cellular in both cases. We now recall the various elements of these cellular structures.

\medskip
For the Temperley-Lieb algebra, according to the notation introduced in Definition
\ref{cellular algebra}, we take $\Lambda = \pacol{n}$, ordered by dominance.
We take $T(\lambda) =\Std(\lambda)$, for all $\lambda \in \Lambda$ and for $\sigma$, $\tau \in T(\lambda)$ we
 define $c_{\sigma \tau}^{\lambda} = \beta_{\sigma \tau}$.

\medskip
Similarly, for the blob algebra we take $\Lambda=\bi{n}$, ordered by  $\succeq$ introduced in the previous section.
Set $T(\bs{\lambda}) = \Std(\bs{\lambda})$, for all $\bs{\lambda} \in \bi{n}$. Given $\sub{s,t} \in T(\bs{\lambda})$
define $c_{\sub{st}}^{\bs{\lambda}} = m_{\sub{st}}$.

\medskip
For $\bs{\lambda} \in \Lambda$, let $ b_n^{\bs{\lambda}}(m)$ be the
ideal of $\blob{m}$ spanned by the set
$$\{ m_{\sub{st}}  \mbox{ } | \mbox{ } \sub{s,t} \in \Std(\bs{\mu}) \text{ ; } \bs{\mu} \succ \bs{\lambda}  \}.$$
In the cases of $\blobR{m}$ and $\blobK{m}$ we write $b_{n}^{\OO, \bs{\lambda}}(m)$ and
$b_{n}^{\K, \bs{\lambda}}(m)$
for the ideals.

\medskip
Since we are assuming that $ q+ q^{-1} \neq 0$ we get that the bilinear forms $ \langle  \cdot , \cdot  \rangle_{\lambda} $ are all
nonzero, in the Temperley-Lieb case as well as the blob algebra case. From this we get
from remark (3.10) of \cite{gra ler} that both algebras are quasi-hereditary and that the cell modules are
standard modules in the sense of quasi-hereditary algebras.

\subsection{Jucys-Murphy elements.} \label{section Jucys}
We are now in position to give the definition of an algebra with Jucys-Murphy elements.
It provides an abstract setting for carrying out much of Murphy's theory
for Young's seminormal form.
Assume
that $A$ is a cellular algebra with cellular basis $$ \mathcal{C} =\{ c_{\sub{st}}^{\lambda} \mbox{  } | \mbox{  } \lambda\in \Lambda ; \sub{s},\sub{t} \in T(\lambda)  \}$$ as in Definition \ref{cellular algebra}. Assume furthermore that each
$ T(\lambda) $ is a poset with respect to an order $ <_{\lambda} $, or just $ < $ for simplicity.
The following definition is taken from \cite{Mat-So}.

\begin{defi} \rm
A family of Jucys-Murphy (JM) elements for A is a set $\{L_1 ,\ldots , L_k\}$ of commuting
elements of A together with a set of scalars,
$$\{ c_{\sub{s}} (i) \in R \mbox{  }| \mbox{  } \sub{s} \in T (\lambda) \mbox{, } \lambda \in \Lambda \mbox{ and } 1\leq i \leq k \}$$ such that for $i = 1,\ldots , k $ we have $L_i^{*}= L_i$ and, for all $\lambda \in \Lambda $ and $\sub{s,t} \in \Lambda $,

$$   L_ic_{\sub{st}}^{\lambda} \equiv c_{s}(i)c_{\sub{st}}^{\lambda} + \sum_{\sub{v} > \sub{s} } r_{\sub{sv}} c_{\sub{vt}}^{\lambda}  \mod{  A^{\lambda} }  $$
for some $r_{\sub{sv}} \in R$ (which depends on $i$). We call $c_{\sub{s}}(i)$ the content of $\sub{s}$ at $i$.
\end{defi}

The purpose of this section is now to apply this definition to $ \blob{m}$.
By Theorem \ref{isomorphism} we have a homomorphism from $\hecke$ onto  $\blob{m}$.  Using (\ref{jucysmurphy como t}) it is easy to note that this homomorphism
maps the elements $L_{k} \in \hecke$ to
$$   (U_{k-1} +q) \ldots  (U_{1} +q)  ((q-q^{-1})U_0 +q^{m}) (U_{1} +q) \ldots (U_{k-1} +q) \in \blob{m}.$$
We shall use the same notation $L_{k}$ for this element of $\blob{m}$.
It satisfies the following commutation rules with the $ U_i $
\begin{align}
 \jm{k}U_{i}  = & U_{i}\jm{k}  &      & \text{if } k\neq i, i+1 \label{conmuttation 1} \\
  (U_{k}+q^{-1}) \jm{k+1} = & \jm{k}(U_{k}+q)  &   &  \text{for } 1\leq k <n.  \label{cunmutation 2} \\
  \jm{k+1}(U_{k}+q^{-1})  = &(U_{k}+q) \jm{k}  &   &  \text{for } 1\leq k <n.  \label{cunmutation 3}
\end{align}
It is known that the $ L_k $ are a family of JM-elements for $\hecke$ with respect to the cellular basis
used for example in \cite{hu-mathas}, in which $ <_{\lambda} $ is the dominance order
on bitableaux.
One might now hope that the set $\{ L_{1},\ldots L_{n}\}$ is also a family of JM-elements for $\blob{m}$.
That this should be the case is not at all obvious.
Indeed, the concept of a family of JM-elements depends heavily on the underlying cellular basis
and a cellular algebra may in general be endowed with several, completely different, cellular bases with different orders.
For example the conjectures of Bonnaf\'e, Geck, Iancu and Lam in \cite{BGIL}, indicate
that Lusztig's theory of cells for unequal parameters should give rise to
a cellular basis on $ \hecke $ for each choice of a weight function on the Coxeter group of type $ B$, in
dependence of
a parameter $ r$. In the setting of these conjectures,
only the asymptotic case $r> n $ corresponds to the dominance order. On the contrary, in
\cite{Steen1} it is shown that the cell structure on $ \blob{m} $ corresponds to
the other extreme case $ r = 0 $
under restriction to $ \bi{n}$ and one-line bitableaux.

\medskip
In this section we shall show that in fact
$\{ L_{1},\ldots L_{n}\}$ do form a family of JM-elements for $\blob{m}$
where the poset structure $T(\bs{\lambda}) = \Std (\bs{\lambda})$ is the one defined above.
Even more, using the surjection $ \heckeintegral   \longrightarrow \blobR{m} $ given
in Theorem \ref{isomorphismintegral}, we define elements
$\{ \jm{1},\ldots , \jm{n}\}$ of $ \blobR{m} $ using the same formula as before and we show that these form a
family of JM-elements for $ \blobR{m} $ with respect to $ \{  m_{\sub{st}} \} $, considered as elements of
$ \blobR{m} $.

\begin{defi}{\label{content}}
Suppose that $\bs{\lambda} \in \bi{n}$ and let $\sub{t} \in \Std(\bs{\lambda}) $.  Let $k$ be an integer with  $1\leq k\leq n$. Define the content of $\sub{t}$ at $k$ to be the scalar given by
$$ c_{\sub{t}}(k)= \left\{
  \begin{array}{lll}
   q^{2(c-1)} Q  & \text{ if } &  d=1  \\
   q^{2(c-1)} Q^{-1}  & \text{ if } & d=2
  \end{array}
\right.   $$
where $(1,c,d) $ is the unique node in $[\bs{\lambda}]$ such that $\sub{t}(1,c,d)=k$.
In other words, $c_{\sub{t}}(k)$ is an element of either $  \OO ,  \C(q,Q) $ or $ \C $, depending on the context.
In the $ \C $ case, we shall also write $ r_{\sub{t}}(k) := c_{\sub{t}}(k) $ and refer to it as the residue of $k$.
\end{defi}

\begin{lem}
\label{base} Suppose that $\bs{\lambda}\in \bi{n} $ and let $k$ be an integer with  $1\leq k\leq n$.
Then we have the
identity
$$\jm{k} m_{  \sub{t}^{\bs{\lambda}}  \sub{t}^{\bs{\lambda}} } \equiv c_{ \sub{t}^{\bs{\lambda}}}(k)
m_{  \sub{t}^{\bs{\lambda}}  \sub{t}^{\bs{\lambda}} }  \mod b_n^{\bs{\lambda}}(m).$$
Similar statements hold over $b_n^{\mathcal{O}}(m)$ and $b_n^{\mathcal{K}}(m)$.
\end{lem}
\begin{dem}
Using the description of $ t^{\bs{\lambda}} $ given after Definition
\ref{orden blob} together with Definition \ref{defi biy blob} we find that
the diagram corresponding to $m_{\sub{t}^{\bs{\lambda}} \sub{t}^{\bs{\lambda} }}$ is one of the diagrams that appear in
Figure \ref{holalalalal}.
But then the statement of the Lemma is
equation (28) of \cite[Lemma 7.1]{blob positive}.
Indeed, using the notation of \cite{blob positive},
equation (28) is the following one
$$
L_j \eta_t = \left\{
\begin{array}{ll}
-x^{2j  -2y +t -n} \eta_t & \mbox{if } j \geq n-| t | \mbox{ and } t > 0 \\
 -x^{2j   -t -n} \eta_t & \mbox{if } j \geq n-| t | \mbox{ and } t <  0 \\
 -x^{j -1  -2 y } \eta_t & \mbox{if } j < n-| t | \mbox{ and } j  \mbox{ is odd} \\
 -x^{j   } \eta_t & \mbox{if } j < n-| t | \mbox{ and } j  \mbox{ is even}
\end{array} \right.
$$
where we have actually corrected an error of {\it loc. cit.}
Indeed to get the correct formulas one should
subtract $2$ from all appearing exponents of $ x $, since the relation between $ L_i $ and
$L_i^{\prime} $ introduced two pages earlier of {\it loc. cit.} should be corrected the same way.
The conversion from the notation used in \cite{blob positive} to ours is now straightforward but
somewhat tedious.
For the reader's convenience we note that
$ \eta_t $ corresponds to the upper part of our
$m_{  \sub{t}^{\bs{\lambda}}  \sub{t}^{\bs{\lambda}} } $ where $ h = \frac{n-|t|}{2}$, and that $ x $ corresponds to
our $ q $ whereas $ y $ corresponds to our $  m $.
Apart from that,
the $ L_j $'s of {\it loc. cit.} have indices
belonging to $ 0, 1, \ldots, n-1 $ whereas ours have indices belonging to
$ 1, 2, \ldots, n $ and finally, since  $ L_0 $ of {\it loc. cit.} satisfies the
relation $ (L_0+x^{-2y})( L_0 +1) = 0 $, we get that our $ L_j $ corresponds to
$ -x^{y} L_{j-1} $ of {\it loc. cit.} for all relevant $ j $.
\end{dem}

\begin{figure}[h!]
\[
\xy 0;/r.11pc/:
(0,40); (10,40) **\crv{(5,25)};
(0,0); (10,0) **\crv{(5,15)};
(5,32.5)*{\bullet};
(5,7.5)*{\bullet};
(20,40); (30,40) **\crv{(25,25)};
(20,0); (30,0) **\crv{(25,15)};
(25,32.5)*{\bullet};
(25,7.5)*{\bullet};
(40,20)*{\ldots};
(50,40); (60,40) **\crv{(55,25)};
(50,0); (60,0) **\crv{(55,15)};
(55,32.5)*{\bullet};
(55,7.5)*{\bullet};
(70,0); (70,40) **\dir{-};
(80,0); (80,40) **\dir{-};
(90,20)*{\ldots};
(100,0); (100,40) **\dir{-};
(-10,0); (-10,40) **\dir{.};
(-10,0); (110,0) **\dir{.};
(110,0); (110,40) **\dir{.};
(110,40); (-10,40) **\dir{.};
(0,45)*{\scriptscriptstyle 1};
(10,45)*{\scriptscriptstyle 2};
(20,45)*{\scriptscriptstyle 3};
(30,45)*{\scriptscriptstyle 4};
(60,45)*{\scriptscriptstyle 2h};
(100,45)*{\scriptscriptstyle n};
(70,20)*{\bullet};
(50,-10)*{\scriptscriptstyle \bs{\lambda}=((h),(h+v))};
(180,-10)*{\scriptscriptstyle \bs{\lambda}=((h+v),(h)) };
(140,40); (150,40) **\crv{(145,25)};
(145,32.5)*{\bullet};
(140,0); (150,0) **\crv{(145,15)};
(145,7.5)*{\bullet};
(160,40); (170,40) **\crv{(165,25)};
(165,32.5)*{\bullet};
(160,0); (170,0) **\crv{(165,15)};
(165,7.5)*{\bullet};
(180,20)*{\ldots};
(190,40); (200,40) **\crv{(195,25)};
(195,32.5)*{\bullet};
(190,0); (200,0) **\crv{(195,15)};
(195,7.5)*{\bullet};
(210,0); (210,40) **\dir{-};
(220,0); (220,40) **\dir{-};
(230,20)*{\ldots};
(240,0); (240,40) **\dir{-};
(130,0); (130,40) **\dir{.};
(130,0); (250,0) **\dir{.};
(250,0); (250,40) **\dir{.};
(250,40); (130,40) **\dir{.};
(140,45)*{\scriptscriptstyle 1};
(150,45)*{\scriptscriptstyle 2};
(160,45)*{\scriptscriptstyle 3};
(170,45)*{\scriptscriptstyle 4};
(200,45)*{\scriptscriptstyle 2h};
(240,45)*{\scriptscriptstyle n};
\endxy
\]
\caption{$m_{\sub{t}^{\bs{\lambda}} \sub{t}^{\bs{\lambda} }}$} \label{holalalalal}
\end{figure}

\medskip
Our proof that the $ \{ \jm{k} \} $ form a family of JM-elements shall be a downwards
induction over the partial order $\succeq$ with
the preceding Lemma providing the induction basis.
To obtain the inductive step we need to understand the relationship between the action of $U_i$ and
$\succeq$ and hence we would like to have a formula for the action of $U_k $ in terms of walks on
the Bratteli diagram. In general there is no such simple formula. On the other hand, there is one situation
where the action of $ U_k $ is particularly  easy to visualize, namely that of a hook being
made smaller by the action of a simple transposition $ s_k$.


\begin{lem} \label{particularly_easy}
\label{pro 1} Suppose that $\bs{\lambda} \in \bi{n}$ and $\sub{s} , \sub{t} \in \Std(\bs{\lambda}) $.
Assume moreover that $s_k{\mathfrak{s}}={\mathfrak{t}}$ for the simple transposition $ s_k $
and that ${\mathfrak{s}}\succ {\mathfrak{t}}$ or equivalently, that
$ w(\mathfrak{t} )$ is obtained from $ w(\mathfrak{s}) $ by making a hook at position $k$ smaller.
Then the following relation holds in $\blob{m}$
$$ U_km_{\sub{s} \sub{t}^{\bs{\lambda}}} = \left\{
  \begin{array}{rll}
   m_{\sub{t} \sub{t}^{\bs{\lambda}}}  & \text{if} &  \sub{s} \nsim \sub{t}   \\
    y_em_{\sub{t} \sub{t}^{\bs{\lambda}}} & \text{if} &  \sub{s} \sim \sub{t}.
  \end{array}
\right.   $$
Similar formulas hold over $ \OO$ and $ K $.  For the Temperley-Lieb algebra we also have
\begin{equation}
U_k\beta_{\sigma \tau^{\lambda}}=\beta_{\tau\tau^{\lambda}}
\end{equation}
where $\lambda \in \pacol{n}$, $\sigma, \tau \in \Std (\lambda)$,  $s_{k}\sigma=\tau$ and $\sigma \unrhd \tau$.
\end{lem}
\begin{dem}
Let us first consider the case $\sub{s} \nsim \sub{t}$. Since
$ w(\mathfrak{t} )$ is obtained from $ w(\mathfrak{s}) $ by making a hook at position $k$ smaller
we have that the sign sequences for $ w(\sub{s}) $ and $ w(\sub{t}) $ are the same,
except at the positions $ k $ and $ k+1 $ where $ w(\sub{s}) $ has signs $ -, +$ whereas
$ w(\sub{t}) $ has signs $ +, -$.
Using the definition of
the maps $ f^{\sigma}_{n, j}$, the claim on the action of $U_k $
is a now a direct consequence. The case
$\sub{s} \sim \sub{t}$ is treated similarly whereas
the Temperley-Lieb case follows by deleting the decorations
in the blob diagrams and using the result for the blob algebra. We remark that in this
case the result in the Lemma is obtained in \cite[Lemma 8 (i)]{Martin}.
\end{dem}

\medskip
The next three Lemmas are preparations for Lemma \ref{ayuda}.

\begin{lem} \label{key lemma TL-case}
Suppose that $\lambda \in \pacol{n}$ and $\sigma, \tau, u \in \Std(\lambda)$.
Suppose moreover that $u \triangleright  \sigma \triangleright \tau $ and that
$ s_k \sigma = \tau $ for some $ k $.
Let $v \in \Std(\lambda)$ be chosen such that
$   U_k \beta_{u t^{\lambda} } = r \beta_{v t^{\lambda} }   \mod{ Tl^{\lambda}} $
for some scalar $ r \in \F$ (such $v$ always exists by the diagrammatical realization of the Temperley-Lieb
algebra and its cell modules). Then, if $ r $ is nonzero we have that $ v \triangleright \tau$.
\end{lem}
\begin{dem}
We identify $\sigma, \tau $ and $ u $ with their walks $w(\sigma), w(\tau) $ and $w( u) $ on the Bratteli diagram for
$ Tl_n$, and also with their corresponding sign sequences.
Then the sign sequences for
$ \sigma $ and $\tau$
are the same except at the $ k $'th and $k+1$'st positions where the sequence for
$\sigma  $ has $ -, + $ whereas the sequence for $\tau $ has $ +, -$.
On the other hand, for $u $ all four possibilities of signs may occur at these positions, apriori,
and so we proceed by a case by case analysis.

\medskip
The first case to analyse is the case where the signs for $ u $ are $ +,- $ at these positions. In this case we get
$v = u $ (and $ r:= -[2]$), and so the claim of the Lemma follows from the assumptions. The next case is the
one where the signs are $ +,+ $ at positions $k$ and $k+1$. On the diagrammatic level we have three options for the
top edge of $\beta_{u\tau^{\lambda}}$, illustrated in Figure \ref{diagrams TL case++}.

\begin{figure}[h!]
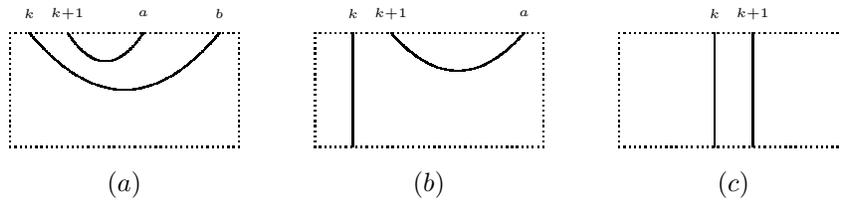

\[
\xy 0;/r.6pc/:
(0,0); (10,0) **\crv{(5,-6)};
(2,0); (6,0) **\crv{(4,-3)};
(0,0); (10,0) **\crv{(5,-6)};
(2,0); (6,0) **\crv{(4,-3)};
(0,1)*{\scriptscriptstyle k};
(2,1)*{\scriptscriptstyle k+1};
(6,1)*{\scriptscriptstyle a};
(10,1)*{\scriptscriptstyle b};
(-1,0); (11,0) **\dir{.};
(-1,0); (-1,-6) **\dir{.};
(-1,-6); (11,-6) **\dir{.};
(11,0); (11,-6) **\dir{.};
(5,-8)*{(a)};
(15,0); (27,0) **\dir{.};
(15,0); (15,-6) **\dir{.};
(15,-6); (27,-6) **\dir{.};
(27,0); (27,-6) **\dir{.};
(19,0); (26,0) **\crv{(22.5,-4)};
(17,0); (17,-6) **\dir{-};
(17,1)*{\scriptscriptstyle k};
(19,1)*{\scriptscriptstyle k+1};
(26,1)*{\scriptscriptstyle a};
(21,-8)*{(b)};
(31,0); (43,0) **\dir{.};
(31,0); (31,-6) **\dir{.};
(31,-6); (43,-6) **\dir{.};
(43,0); (43,-6) **\dir{.};
(36,0); (36,-6) **\dir{-};
(38,0); (38,-6) **\dir{-};
(36,1)*{\scriptscriptstyle k};
(38,1)*{\scriptscriptstyle k+1};
(37,-8)*{(c)};
\endxy
\]
\caption{ Top edge of $\beta_{u\tau^{\lambda}}$ } \label{diagrams TL case++}
\end{figure}

In the subcase $(a)$, the signs for $u$ at positions $k$, $k+1$, $a$ and $b$ are $+,+,-$ and $-$, respectively. For $v$ the signs in these positions are $+,-,+$ and $-$, whereas the signs for $v$ and $u$ agree at all other positions. The claim follows from this. The subcase $(b)$ is treated similarly. Finally, in the subcase $(c)$ we have $ r=0 $, contrary to the assumptions.


\medskip
The third case is the one where
the signs for $u $ are $ -,+ $ at the positions $ k, k+1 $. In that case, at the diagrammatic level, $ k $ is connected to a point $ a < k $ whereas $ k+1 $ is either
connected to $ b > k+1$ or it is the upper endpoint of a vertical line. In both cases,
we find that the sign sequence for
$ v $ is the same as the one for $u $, except at positions $ k, k+1 $ where it becomes
$ +, -$. But by the assumptions, $u $ differs from
$\sigma $ in at least one position and the result follows in this case as well. Note that this is the only case in which $u \vartriangleright v$.

\medskip
The last case is the one where the signs for $ u $ at the positions $ k, k+1 $ are $-, - $.
In this case, $ k $ is connected to $ a $ and $ k+1 $ to $ b$ and $ b < a < k < k+1$. Moreover the signs for $u$
at these positions are $ +,+,-,-$. But then the signs for $ v $ at these positions are $ +, -, +, - $ whereas
the signs for $ v $ and $ u $ agree at all other positions. The claim follows from this.
\end{dem}

\begin{lem}  \label{comparableTL}
Suppose that $\mu\in \pacol{n} $. Let $\sigma, \tau \in \text{Std}(\mu)$. Assume that $U_k\beta_{\sigma \tau^{\mu}} \equiv  \alpha \beta_{\tau \tau^{\mu}} \mod TL_{n}^{\mu}$, with $\alpha \neq 0$ and $1\leq k < n$. Then, $\tau \trianglerighteq \sigma$, or $\sigma \vartriangleright \tau $  and $s_k\sigma =\tau$.
\end{lem}

\begin{dem}
The result follows by a case by case analysis, similar to that given in the proof of the previous Lemma.
\end{dem}

\begin{lem} \label{lemma de como de descubren los nodos}
Let $\bs{\lambda} \in \bi{n}$ and $\sub{u}\in \Std(\bs{\lambda})$.
Assume that $U_km_{\sub{ut}^{\bs{\lambda}}} \equiv \alpha m_{\sub{vt}^{\bs{\lambda}}}
\mod  b_n^{\bs{\lambda}}(m)$, for $\alpha \in \C $ and $\sub{v} \in \Std(\bs{\lambda})$
(such $\sub{v} $ always exists by the diagrammatic realization of $ \blob{m}$).
Suppose moreover that $ \alpha $ is nonzero
and that the node $j$ is
covered in the top edge of $m_{\sub{ut}^{\bs{\lambda}}}$, but uncovered in the top
edge of $m_{\sub{vt}^{\bs{\lambda}}}$. Then $|\sub{u}(j)| = | \sub{v}(j)| \neq 0$ if and only if $j=k$.
Similar statements hold over $ \OO $ and $ K$.
\end{lem}
\begin{dem}
In order for the action of $U_{k}$ to transform a covered node $ j$ in the top edge of $m_{\sub{ut}^{\bs{\lambda}}}$
to an uncovered node in the top edge of $m_{\sub{vt}^{\bs{\lambda}}}$, the diagram of $m_{\sub{ut}^{\bs{\lambda}}}$
must be one of those shown in the below Figure \ref{dibujodecomosedescubrennodos}
with the position of $ j $ shown in each case.
Using this classification, the Lemma follows from Definition \ref{defi biy blob}.

\begin{figure}[h!] \label{dibujodecomosedescubrennodos}
\xy 0;/r.10pc/:
(0,0);(60,0)  **\crv{(30,-27)};
(10,0);(40,0)  **\crv{(25,-10)};
(0,5)*{\scriptscriptstyle k};
(10,5)*{\scriptscriptstyle k+1};
(40,5)*{\scriptscriptstyle a};
(30,-14)*{\bullet};
(30,-40)*{ \textstyle k\leq j <a};
(97,0);(113,0)  **\crv{(105,-15)};
(105,-8)*{\bullet};
(97,5)*{\scriptscriptstyle k};
(113,5)*{\scriptscriptstyle k+1};
(105,-40)*{k\leq j \leq k+1};
(150,0);(210,0)  **\crv{(180,-27)};
(170,0);(200,0)  **\crv{(185,-10)};
(180,-14)*{\bullet};
(200,5)*{\scriptscriptstyle k};
(210,5)*{\scriptscriptstyle k+1};
(170,5)*{\scriptscriptstyle a};
(180,-40)*{a < j \leq k+1};
(240,0); (240,-30) **\dir{-};
(240,-15)*{\bullet};
(240,5)*{\scriptscriptstyle k};
(250,0);(300,0)  **\crv{(275,-27)};
(250,5)*{\scriptscriptstyle k+1};
(300,5)*{\scriptscriptstyle a};
(265,-40)*{k \leq j < a};
\endxy
 \caption{Possibilities for $m_{\sub{u} \sub{t}^{\bs{\lambda}}}$} \label{dibujodecomosedescubrennodos}
\end{figure}
\end{dem}

\medskip
We can now finally prove the property of the order $ \succ $ that makes our induction work.
It is a generalization to the blob algebra case of Lemma \ref{key lemma TL-case}, and in fact
we shall deduce it from that Lemma.
\begin{lem} \label{ayuda}
Suppose that $\bs{\lambda} =((a),(b))\in \bi{n}$ and $\sub{s},\sub{t}, \sub{u} \in \Std(\bs{\lambda})$.
Suppose furthermore that $s_k\sub{s}=\sub{t}$ and that $\sub{u}\succ \sub{s} \succ \sub{t}$.
Let $\sub{v} \in \Std(\bs{\lambda})$ be chosen such that
$   U_k m_{\sub{ut}^{\bs{\lambda}}} = r m_{\sub{vt}^{\bs{\lambda}}}  \mod{ b_n^{\bs{\lambda}}(m) } $
for some scalar $ r \in \F$.
Then, if $r $ is nonzero we have that
$\sub{v} \succ \sub{t}$.
Similar statements are valid for $\blobR{m}$ and $\blobK{m}$.
\end{lem}
\begin{dem}
Set $\mu_1=\max \{ a,b \}, \, \mu_2= \min \{ a,b\}$
and let $\mu=(\mu_1,\mu_2)' $. Then $\mu \in \pacol{n} $ and in the
Temperley-Lieb algebra we have that
$$  U_k \tld{\tc{u}} = r_{1} \tld{\tc{v}}  \mod   Tl^{\mu}  $$
where $\tc{u},  \tc{v} $ are as in Definition \ref{tipo tl} and $ r_{1} \neq 0$; indeed
$r_{1}= -[2]$ if $r =-[2]$ or if $r =y_e$, and $r_{1}=1$ if $ r=1$.
Moreover, by Lemma \ref{orden con diagrmas} we have that $ \tc{u} \trianglerighteq \tc{s}  \trianglerighteq \tc{t} $.


\medskip
\noindent
\textbf{Case 1} ($\tc{u} \triangleright \tc{s}  \triangleright \tc{t} $). In this case we have by Lemma  \ref{orden tau} that $ s_k \tc{s} = \tc{t} $ and then Lemma \ref{key lemma TL-case} gives that $ \tc{v}  \triangleright \tc{t} $. Now by Lemma \ref{order with sequence}, in order to prove that $\sub{v} \succ \sub{t}$ it is
enough to show that
\begin{equation} \label{condicion de la secuencia}
| \sub{v}(j)| = |\sub{t}(j)| \text{ implies }   \sub{v}(j) \leq \sub{t}(j).
\end{equation}
Hence, assume that $| \sub{v}(j)| = |\sub{t}(j)|$,
but $\sub{t}(j)<0$ and $\sub{v}(j)>0$ for some $1\leq j \leq n$. We now split this case into two subcases according to Lemma \ref{comparableTL}, that is, $\tc{v} \trianglerighteq \tc{u}$ or, $\tc{u} \vartriangleright \tc{v}$ and $s_k\tc{u}=\tc{v}$.  First, we assume that $\tc{v} \trianglerighteq \tc{u}$. Then we get from $\sub{u}\succ \sub{s} \succ \sub{t}$ and Lemma \ref{order with sequence} that $ \sub{u}(j) = \sub{s}(j) = \sub{t}(j) $ and so
we get $\sub{u}(j)<0$, $\sub{v}(j) > 0$ and $| \sub{u}(j)| = |\sub{v}(j)|$. From this we conclude
via Lemma \ref{lemma de como de descubren los nodos} that $ j = k $, hence
that $\sub{s}(k)= \sub{t}(k)$, which is impossible because $s_k\sub{s}=\sub{t}$.

\medskip
So we can assume that $\tc{u} \vartriangleright \tc{v}$ and $s_k\tc{u}=\tc{v}$.
Assume first that $ k \neq j$. Then $\sub{s}(j) = \sub{t}(j) $ whereas
$| \sub{u}(j) |= \sub{v}(j) $. But $\sub{u}\succ \sub{s} $ implies $  \sub{u}(j)  = \sub{s}(j)  $ and
so we have $ \sub{u}(j) = \sub{s}(j) = \sub{t}(j) $. As before this implies via
Lemma \ref{lemma de como de descubren los nodos} that $ j = k $, contradiction.
We then assume $ j= k$. Then $ w(\mathfrak{t} )$ is obtained from $ w(\mathfrak{s}) $ by
making a hook at position $ j$ bigger and so we have $ \sub{s}(j) = \sub{t}(j) + 2 \leq 0 $ and
$ \sub{s}( j \pm 1 ) = \sub{t}(j \pm 1) = \sub{t}(j)+1$.
On the other hand, $\tc{u} \vartriangleright \tc{v}$ and $s_k\tc{u}=\tc{v}$ imply that
$| \sub{u}(j) |= \sub{v}(j) -2 $ which combined with $\sub{u}\succ \sub{s} $ gives
$\sub{u}(j) = \sub{s}(j) $ and then also $\sub{u}(j \pm 1) = \sub{s}(j \pm 1) $, since
$\sub{u}( i) $ at most changes by $ \pm 1 $ at each step. But then Lemma \ref{particularly_easy}
implies that $ \sub{v}(j) = \sub{t}(j) $, contradiction.

\medskip
\noindent
\textbf{Case 2} ($\tc{u} \trianglerighteq \tc{s}  = \tc{t} $).
By the assumptions
$ \sub{t} $ is obtained from $ \sub{s} $ by making a hook at position $ k $ smaller.
Moreover,  since $ \tc{s}  = \tc{t}  $ this hook is located on the central vertical axis of the Bratteli diagram, that is
$\sub{t}(k-1) = \sub{s}(k-1) = \sub{t}(k+1) = \sub{s}(k+1) = 0 $.
But then, since $\sub{u} \succ \sub{s} $, we have
necessarily that
$\sub{u}(k-1) = \sub{u}(k+1) = 0,  \, \sub{u}(k) = -1 $ which implies via
Lemma \ref{particularly_easy} that
$  \sub{v} $ is obtained from $   \sub{u} $ by making a hook at position $k $ smaller.
Hence we get $ \sub{v} \succ  \sub{t} $ as claimed.

\medskip
\noindent
\textbf{Case 3} ($\tc{u} = \tc{s}  \trianglerighteq \tc{t} $). By the hypothesis in this case we have $\tc{v}= \tc{t}$.
Recall that at the Bratteli diagram level this implies that at each step the walks $w(\sub{t})$ and $w(\sub{v})$
are either equal or mirror images under the reflection through the central vertical axis of the Bratteli diagram.
So, in order to prove the Lemma in this case we must prove that whenever the path $w(\sub{t})$ is on the negative side of the Bratteli diagram, the path $w(\sub{v})$ is also  on the negative part. In terms of the sequence of integers the last condition is equivalent to
\begin{equation} \label{t negative then v negative}
\sub{t}(j) <0 \text{ implies } \sub{v}(j)<0
\end{equation}
for all $1 \leq j \leq n$. Suppose by contradiction that (\ref{t negative then v negative}) is not true for some $1\leq j \leq n$. Therefore, $\sub{t}(j)<0<\sub{v}(j)$ for some $1\leq j \leq n$. Using the fact that $\sub{s}(j)=\sub{t}(j)$, for all $j\neq k$, $\tc{u} = \tc{s}$ and $\tc{v}= \tc{t}$, we can conclude via Remark \ref{remark covered blob sequence} and Lemma \ref{lemma de como de descubren los nodos} that $j=k$. Hence, at step $k$ the walk $w(\sub{t})$ (resp. $w(\sub{v})$) is on the negative side of Bratteli diagram and (\ref{t negative then v negative}) is true for all $j\neq k$. This implies that $\sub{t}(k-1)=\sub{t}(k+1)=0$ and $\sub{t}(k)=-1$. But this is impossible because $\sub{s}\succ \sub{t}$ and $s_k\sub{s}=\sub{t}$. This completes the proof of the Lemma.
\end{dem}

\medskip
We are now in position to prove the triangularity property for
$\{ \jm{1},\ldots , \jm{n}\}$.
It follows from it that the set $\{ \jm{1},\ldots , \jm{n}\}$ is a family of JM-elements for the blob algebra
with respect to the order $ \succ$.

\begin{teo} \label{triangular}
Suppose that $ \bs{\lambda} \in \bi{n}$ and $s,t\in \Std(\bs{\lambda})$. Then
$$   \jm{k}m_{\sub{st}}=c_{\sub{s}}(k) m_{\sub{st}} + \sum_{\substack{ \sub{u} \in \Std(\bs{\lambda}) \\ \sub{u} \succ \sub{s}}}  a_{\sub{u}} m_{\sub{ut}} \mod{ b_n^{\bs{\lambda}}(m) } $$
for some scalars $a_{\sub{u}} $. A similar statements holds for $\blobR{m}$ and $\blobK{m}$.
\end{teo}
\begin{dem}
By the cellularity of the diagram basis,
the statement of the Lemma is independent of
$\sub{t} $. We proceed by induction on the order $\succeq$. The induction basis $\sub{s} = \sub{t}^{\bs{\lambda}}$ is provided by Lemma \ref{base}.
Assume now that $\sub{s} \neq \sub{t}^{\bs{\lambda}}$. Then we can find $ i $ and $\sub{s}^{\prime} $ such
that $ \sub{s}^{\prime}  \succ \sub{s}  $ and $ s_i \sub{s}^{\prime}  = \sub{s}  $.
By the inductive hypothesis the Theorem is valid for $\sub{s}^{\prime} $.
We first assume that
$\sub{s} \nsim \sub{s}^{\prime} $ and $ k \neq i, i+1 $. Using Lemma \ref{pro 1} and
the commutation rule (\ref{conmuttation 1}) we then get
$$
 \jm{k}m_{\sub{st}} =  \jm{k} U_i m_{\sub{s^{\prime}t}}  =  U_i \jm{k}  m_{\sub{s^{\prime}t}} =
c_{\sub{s^{\prime}}}(k)  m_{\sub{st}} + \sum_{\substack{ \sub{u} \in \Std(\bs{\lambda}) \\ \sub{u} \succ \sub{s^{\prime}}}}
a_{\sub{u}} U_i m_{\sub{ut}} \mod{ b_n^{\bs{\lambda}}(m)}.
$$
On the other hand, by the previous Lemma the sum is a linear combination of elements of the
form $ m_{\sub{ut}} $ where $ \sub{u} \succ \sub{s} $ and since $ c_{\sub{s}}(k) =  c_{\sub{s^{\prime}}}(k)   $
we are done in this case.

\medskip
If $\sub{s} \sim \sub{s}^{\prime} $ and $ k \neq i, i+1 $ we find similarly
$$
 \jm{k}m_{\sub{st}} =  y_e^{-1} \jm{k} U_i m_{\sub{s^{\prime}t}}  = y_e^{-1}   U_i \jm{k}  m_{\sub{s^{\prime}t}} =
c_{\sub{s^{\prime}}}(k)  m_{\sub{st}} + \sum_{\substack{ \sub{u} \in \Std(\bs{\lambda}) \\ \sub{u} \succ \sub{s^{\prime}}}}
y_e^{-1} a_{\sub{u}} U_i m_{\sub{ut}} \mod{ b_n^{\bs{\lambda}}(m)}
$$
and may conclude the same way as before.
We next treat the case $\sub{s} \nsim \sub{s}^{\prime} $ and $ i= k $ where we find, using
the commutation rule  (\ref{cunmutation 2}) that
$$
 \jm{k}m_{\sub{st}} =   \jm{k} U_k m_{\sub{s^{\prime}t}}  = \jm{k} (U_k+q -q) m_{\sub{s^{\prime}t}}  =
(U_k+q^{-1}) \jm{k+1} m_{\sub{s^{\prime}t}}    -q  \jm{k} m_{\sub{s^{\prime}t}} . $$
By the inductive hypothesis, $ \jm{k} m_{\sub{s^{\prime}t}} $ and $ \jm{k+1} m_{\sub{s^{\prime}t}} $ are
linear combination of elements of the form $ m_{\sub{ut}} $ where $ \sub{u} \succ \sub{s} $
and hence we find, using the inductive hypothesis and Lemma \ref{pro 1} once more, that $\jm{k}m_{\sub{st}} $ is equal to
$$ U_k\jm{k+1} m_{\sub{s^{\prime}t}} =
c_{\sub{s^{\prime}}}(k+1)  m_{\sub{st}} +
\sum_{\substack{ \sub{u} \in \Std(\bs{\lambda}) \\ \sub{u} \succ \sub{s^{\prime}}}}
a_{\sub{u}} U_k m_{\sub{ut}} \mod{ b_n^{\bs{\lambda}}(m)}.
$$
But $ c_{\sub{s}}(k)  =  c_{\sub{s^{\prime}}}(k+1)  $ and we may conclude this case using the previous Lemma as before.
The remaining cases are treated similarly.
\end{dem}

\section{A graded cellular basis of $\blob{m}$.} \label{a graded cellular}
In this section we obtain our main results showing that
$\blob{m}$ is a graded cellular algebra.
Our methods are inspired by the ones used by Hu and Mathas in \cite[Section 4 and 5]{hu-mathas},
who construct a graded cellular basis $ \{ \psi_{st}\} $ for the cyclotomic Hecke algebra,
in terms of the Khovanov-Lauda-Rouquier generators. But
unfortunately is not possible to use their results directly. In fact, the homomorphism
$ \Phi :   \hecke  \longrightarrow   \blob{m}  $ may easily map linearly independent elements to linearly
dependent elements. Moreover,
due to the incompatibility between the dominance order used
for $ \{ \psi_{st}\} $ and the order $ \succ $ for $ \blob{m} $,
we do not know how to find a basis for $ \ker \Phi $ consisting of elements from $ \{ \psi_{st}\} $,
and so in general it seems
intractable to determine which are the subsets of $ \{ \psi_{st}\} $ that stay independent
under $ \Phi$.

\medskip
Our solution to this problem is indirect. It is
based on an alternative realization of the KLR-idempotents
$e(\bs{i})$ which is possible in the setting of an algebra with JM-elements,
see Lemma $4.2$ of \cite{Mat-So}.
It also plays a key role in
\cite{hu-mathas} in the setting of cyclotomic Hecke algebras.
To explain it we first setup the relevant notation.

\medskip
We fix $ {\OO} $ and $ \m $ as above.
Recall that $ \K= \C(q,Q) $ and $  \blobK{m} =  b_n^{\OO} (m) \otimes_{\OO} \K     $.
Over $ \K $ the contents from Definition {\ref{content}} trivially verify the
separation criterion of \cite{Mat-So} and so $ \blobK{m} $ is semisimple.
Hence we can apply \cite{Mat-So} to the algebras
$\blob{m}$, $ \blobR{m}$ and $ \blobK{m} $. We repeat the necessary definitions
in our setting.



\begin{defi} \rm
Suppose that $\bs{\lambda} \in \bi{n}$ and $\sub{s}$, $\sub{t} \in \Std(\bs{\lambda})$. Then we define
$$ F_{\sub{t}} := \prod_{k=1}^{n}  \prod_{\substack{   \sub{s} \in \Std(n) \\ c^{}_{\sub{s}}(k) \neq c^{}_{\sub{s}}(k) }}
\frac{\jm{k}-c^{}_{\sub{s}}(k)}{c^{}_{\sub{t}}(k)-c^{}_{\sub{s}}(k)} \in \blobK{m} $$
and set $f_{\sub{st}} = F_{\sub{s}} m_{\sub{st}}F_{\sub{t}}$.
\end{defi}

We extend the order $\succeq$ to pairs of bitableaux
of the same shape by declaring
$(\sub{u} , \sub{v}) \succeq (\sub{s} , \sub{t})$ if $  \sub{u},\sub{v}  \in \Std(\bs{\lambda})$
and $ \sub{s},\sub{t} \in \Std(\bs{\mu})$, and if either $\bs{\mu} \succeq \bs{\lambda}$ or $\bs{\mu} = \bs{\lambda}$
and $\sub{u} \succeq \sub{s}$ and $\sub{v} \succeq \sub{t}$. Then we get that
\begin{equation} \label{triang seminormal}
f_{\sub{st}} =  m_{\sub{st}} + \sum_{  (\sub{u,v}) \succ (\sub{s,t})} r_{\sub{uv}} m_{\sub{uv}}
\end{equation}
for some $r_{\sub{uv}}\in \K $ and hence
$$   \{ f_{\sub{st}} \mbox{  } | \mbox{   } \sub{s,t} \in \Std(\bs{\lambda}) , \bs{\lambda} \in \bi{n}    \}$$
is a basis for $\blobK{m} $, the seminormal basis. Moreover, by
\cite[Theorem 3.7]{Mat-So}, for
$ \sub{t} \in \Std(\bs{\lambda}) $ there exists a non-zero scalar $\gamma_{\sub{t}} \in \K$
such that
\begin{equation}{\label{gamma}}
f_{\sub{tt}}f_{\sub{tt}}=\gamma_{\sub{t}}f_{\sub{tt}}
\end{equation}

\medskip
Let $ \thickapprox $ be the equivalence relation on $\Std(\bs{n})$ given by $ \sub{s} \thickapprox \sub{t}  $
if $ r_{\sub{s}}(k) = r_{\sub{s}}(k) $ for $ k =1, 2, \ldots, n$. The equivalence classes for
$ \thickapprox $
are parametrized by residue sequences $ I^n$ of length $n$; for $\bs{i} \in I^{n}$ we denote by
$ \Std(\bs{i}) $ the corresponding class. Any tableau $ \sub{s} $ gives rise to a residue sequence that is
denoted $ \bs{i}^{\sub{s}} $. Then we have $ \sub{s} \in \Std( \bs{i}^{\sub{s}} ) $ but in general
$ \Std(\bs{i}) $ may be empty, of course.
For each $\bs{i} \in I^{n}$ we define idempotents
$e^b(\bs{i}) \in \blobK{m}  $ by
$$ e^b(\bs{i}):= \sum_{\sub{s} \in \Std(\bs{i})}   \frac{1}{\gamma_{\sub{s}}} f_{\sub{ss}} . $$
Then it follows from \cite{Mat-So} that actually $e^b(\bs{i}) \in b_n^{\mathcal{O}}$
and so we may reduce $e^b(\bs{i}) $ modulo $ \m$ to obtain idempotents of $ \blob{m}$ that we denote the
same way $e^b(\bs{i}) $.

\medskip
As already mentioned above, the above construction can also be
carried out for the cyclotomic Hecke $ \hecke$, where it gives rise
to idempotents that we denote $  e^{\cal H}(\bs{i})$.
The following Lemma is the key Proposition 4.8 of \cite{hu-mathas}.
\begin{lem}  \label{Proposition 4.8}
For $ \bs{i} = (i_1, i_2, \ldots, i_n ) \in I^n $ let
$$  \hecke(\bs{i})):= \{  v\in \hecke \,\mid
\, (L_r-q^{2i_r})^{k}v=0 \text{ for } r=1,\ldots , n \text{ and } k \gg 0  \} $$
be the generalized weight space for the action of $ L_i \in \hecke $.
Then we have $ \hecke (\bs{i}) = e^{\cal H}(\bs{i})  \hecke$.
In other words, $ e^{\cal H}(\bs{i}) $ is equal to the KLR-idempotent $ e(\bs{i})$.
\end{lem}

We have a similar result for $\blob{m}$.

\begin{lem}  \label{blob-idem}
For $ \bs{i} = (i_1, i_2, \ldots, i_n ) \in I^n $ let
$$  \blob{m}(\bs{i}):= \{  v\in \blob{m} \,\mid
\, (L_r-q^{2i_r})^{k}v=0 \text{ for } r=1,\ldots , n \text{ and } k \gg 0  \} $$
be the generalized weight space for the action of $ L_i \in \blob{m} $.
Then we have $ \blob{m} (\bs{i}) = e^b(\bs{i})  \blob{m}$.
\end{lem}
\begin{dem}
The proof of Proposition 4.8 of \cite{hu-mathas} carries over.
\end{dem}

\begin{lem}  \label{compatibilty-idem}
Let $ \Phi :   \hecke  \longrightarrow   \blob{m}  $ be the homomorphism in  Theorem \ref{isomorphism} and let $ \bs{i} \in I^n $.
Then $ \Phi(e(\bs{i}))=  e^b(\bs{i}) $.
In particular, $ e^b(\bs{i}) $ is a homogeneous element of $ \blob{m} $ of degree $ 0$.
\end{lem}
\begin{dem}
Since $ \Phi $ is surjective and maps the JM-elements of $ \hecke$ to the JM-elements of $ \blob{m}$,
we have $ \Phi(\hecke(\bs{i}))= \blob{m} (\bs{i}) $.
But then $$ e^b(\bs{i}) \blob{m} = \blob{m} (\bs{i})  =
 \Phi(\hecke(\bs{i}))= \Phi( e(\bs{i}) \hecke ) = \Phi( e(\bs{i}))\blob{m} .$$
Moreover, $ \Phi( e(\bs{i}))$ lies in the subalgebra of $ \blob{m} $ generated by the JM-elements
since $ e(\bs{i}) $ has the corresponding property, and so
$ \Phi(e(\bs{i}))=  e^b(\bs{i}) $ as claimed. On the other hand, by Corollary
\ref{main} we know that $\Phi $ is homogeneous and so
the second claim holds as well.
\end{dem}

\begin{rem}{\label{at-this-point}}\rm
At this point we may remark that in the case $ \heckedos $,
the separation criterion of \cite{Mat-So}
corresponds exactly to our standing conditions (\ref{qym}) on the parameters $ Q $ and $q $.
By {\it loc. cit.} it then follows that $ \heckedos $ is semisimple under (\ref{qym})
and that the classes for the corresponding relation $ \thickapprox $ are of size one. Hence,
if $ e^b(\bs{i}) $ is nonzero
we have that
$$  e^b(\bs{i}) = \frac{1}{\gamma_{\sub{s}}} f_{\sub{ss}}$$
for a bitableau $ \sub{s} $ of total degree $2$.
Using this, we obtain an alternative proof of
Theorem \ref{idempotentes homogeneous} since $ e_2^{-1} $ and $ e_2^{-2} $
are idempotents for one dimensional representations of $ \heckedos $.
\end{rem}

\medskip
We next define elements $ \psi_i^b, y_i^b $ of $ \blob{m} $ by
$ \psi_i^b := \Phi (\psi_i ) $ and $ y_i^b := \Phi ( y_i) $. As is the case for $e^b(\bs{i}) $, the elements
$ y_i^b $ and $ \psi_i^b e^b(\bs{i}) $
are homogeneous of the same degree as their Hecke algebra counterparts.
We are now in position to give the key definition of this section.

\begin{defi}\label{bases}
Suppose that $\bs{\lambda} \in \bi{n}$ and $\sub{s} , \sub{t} \in \Std (\bs{\lambda})$.
Let $\de{s} =s_{i_1}\ldots s_{i_k}$ and $\de{t}=s_{j_1} \ldots s_{j_l}$ be reduced expressions for $\de{s}$ and $\de{t}$,
chosen as in Lemma \ref{reduced exp}. Then we
define $$  \psi_{\sub{st}}^b :=    \psi_{i_1}^b\ldots \psi_{i_k}^b   e^b(\bs{i}^{\bs{\lambda}})
 \psi_{j_l}^b\ldots \psi_{j_1}^b \in \blob{m}.     $$
\end{defi}

Note that although our $ \psi_{\sub{st}}^b $ look much like the elements
$ \psi_{\sub{st}} $ introduced in \cite{hu-mathas}, this resemblance is only formal
and in general there is no obvious connection between the two families
of elements, due to the differences between the tableaux.
Note also that in our definition there is no $y$ factor, contrary to the \cite{hu-mathas}
situation. Finally, we note that our $\psi_{\sub{st}}^b  $ can be shown to be independent of the choices
of reduced expressions as above, this is also contrary to the
situation in \cite{hu-mathas}. This independence comes from the fact that
the expressions for $ \de{s} $ and $\de{t} $  are $ iji $-avoiding, that is 
any two expressions
are related through a series of Coxeter relations
of type $ s_i s_j = s_i s_j $ for $ | i-j | > 1 $.

\medskip
Our next result is parallel to Theorem 4.14 of \cite{hu-mathas}, but has no $y$ term. This 'missing'
$y$ is the reason why there is no $ y$ factor in Definition \ref{bases}.

\begin{teo}  \label{import}
Suppose that $\bs{\lambda}=((a),(b)) \in \bi{n}$.  Then there exists a nonzero scalar $r \in \F^{ \times}$ such that
$$
e^{b}(\bs{i}^{\bs{\lambda}}) \equiv rm_{  \sub{t}^{\bs{\lambda}} \sub{t}^{\bs{\lambda}}} \mod{b_n^{\bs{\lambda}}(m)  }.
$$
\end{teo}
\begin{dem}
We begin by determining
$  \gamma_{\sub{t}^{\bs{\lambda}}}$. For this we use (\ref{triang seminormal}) and
(\ref{gamma}) and find
\begin{align*}
\gamma_{\sub{t}^{\bs{\lambda}}}f_{\sub{t}^{\bs{\lambda}} \sub{t}^{\bs{\lambda}} }   &
=  f_{\sub{t}^{\bs{\lambda}} \sub{t}^{\bs{\lambda}} } f_{\sub{t}^{\bs{\lambda}} \sub{t}^{\bs{\lambda}} }  \\
  &  \equiv m_{\sub{t}^{\bs{\lambda}} \sub{t}^{\bs{\lambda}} } m_{\sub{t}^{\bs{\lambda}} \sub{t}^{\bs{\lambda}} }
\quad  \mod{ b_{n}^{\K, \bs{\lambda}}(m) } \\
   & \equiv  (y_{e})^{c} m_{\sub{t}^{\bs{\lambda}} \sub{t}^{\bs{\lambda}} } \quad \mod{ b_{n}^{\K, \bs{\lambda}}(m) }\\
   & \equiv  (y_{e})^{c} f_{\sub{t}^{\bs{\lambda}} \sub{t}^{\bs{\lambda}} } \quad \mod{ b_{n}^{\K, \bs{\lambda}}(m) }
\end{align*}
where $c= \min \{a,b\}$ and where
$ m_{\sub{t}^{\bs{\lambda}} \sub{t}^{\bs{\lambda}} }  m_{\sub{t}^{\bs{\lambda}} \sub{t}^{\bs{\lambda}} } $
can be conveniently found via the diagrammatic realization of
$ m_{\sub{t}^{\bs{\lambda}} \sub{t}^{\bs{\lambda}} }  $ in Figure 6. From this we
deduce that $\gamma_{\sub{t}^{\bs{\lambda}}} = (y_{e})^{c} $.

\medskip
On the other hand,
for $ \sub{s} \in \Std ( \bs{i}^{\bs{\lambda}} )$ with $ \sub{s} \neq \sub{t}^{\bs{\lambda}}$, we get by
combining the description of
$ \sub{t}^{\bs{\lambda}} $ given just after Definition \ref{orden blob}  with the standing conditions on
the parameters (\ref{qym})
that $ \Shape(\sub{s})  \succ \bs{\lambda} $. But then (\ref{triang seminormal})
and the definition of $ e(\bs{i}^{\bs{\lambda}})  $ imply
\begin{equation} \label{good}
e(\bs{i}^{\bs{\lambda}}) \equiv \frac{1}{ (y_{e})^{c}  }
m_{\sub{t}^{\bs{\lambda}} \sub{t}^{\bs{\lambda}}} \mod{ b_{n}^{\K, \bs{\lambda}}(m) }.
\end{equation}
Since $ e(\bs{i}^{\bs{\lambda}})  $ and $ \frac{1}{ (y_{e})^{c}  } m_{\sub{t}^{\bs{\lambda}} \sub{t}^{\bs{\lambda}}}$
both belong to $ b_{n}^{\OO, \bs{\lambda}}(m)  $, we can now replace
$ { b_{n}^{\K, \bs{\lambda}}(m) } $ by
$ { b_{n}^{\OO, \bs{\lambda}}(m) } $ in (\ref{good}). From this the proof is obtained by reducing modulo
$ {\m}   $.
\end{dem}

\medskip
We can now prove that the elements from Definition \ref{bases} form a basis for $ \blob{m}$.

\begin{teo} \label{transicion}
Suppose that $\bs{\lambda} \in \bi{n}$ and $\sub{s}, \sub{t} \in \Std(\bs{\lambda})$.
Then there are scalars $r\in \F^{\times}$ and $r_{\sub{uv}} \in \F$ such that
$$  \psi^b_{\sub{st}}= rm_{\sub{st}} + \sum_{(\sub{u,v}) \succ (\sub{s,t})} r_{\sub{uv}} m_{\sub{uv}} . $$
Hence $\{ \psi^b_{\sub{st}}  \mid \sub{s,t} \in \Std(\bs{\lambda}) \text{ for } \bs{\lambda} \in \bi{n}  \}$ is a basis for
$\blob{m}$.
\end{teo}

\begin{dem}
For $\de{s} = s_{i_{1}}\ldots s_{i_{k}} $ a reduced expression for $\de{s}$
as above we consider first $\psi_{i_1}^b\ldots \psi_{i_k}^b   e^b(\bs{i}^{\bs{\lambda}}) $.
Using (\ref{definicion de psi}) and the commutation rules
(\ref{kl7}), (\ref{kl9}) and (\ref{kl10}) between the $ y_i $ and $ \psi_j $, we get that it can be expressed as a linear
combination of
elements of the form $\Phi(T_{i_{j_1}} \ldots T_{i_{j_r}}  f_{j_1, \ldots, j_r}(y_1, \ldots, y_n) e(\bs{i}^{\bs{\lambda}}) )$
where $ (i_{j_1}, \ldots , i_{j_r}) $ is a subsequence of $( i_1 , \ldots , i_k)$ and where $f_{j_1, \ldots, j_r}(y_1, \ldots, y_n) $
is a polynomial in the $ y_i$.  But $ \Phi(T_i) = q U_i + q^2 $ and hence this can also be written
as a linear combination of elements of the form
$ U_{i_{j_1}} \ldots U_{i_{j_r}}  g_{j_1, \ldots, j_r}(y^b_1, \ldots, y^b_n) e^b(\bs{i}^{\bs{\lambda}}) )$
where $ (i_{j_1}, \ldots , i_{j_r}) $ is a subsequence of $( i_1 , \ldots , i_k)$ and
$g_{j_1, \ldots, j_r}(y^b_1, \ldots, y^b_n) $ is a polynomial in the $ y_i^b$. But from
(\ref{nilpotent}) and  Theorem \ref{import} this is a linear combination of elements of the form
$ U_{i_{j_1}} \ldots U_{i_{j_r}}  m_{  \sub{t}^{\bs{\lambda}} \sub{t}^{\bs{\lambda}}} \mod{b_n^{\bs{\lambda}}(m)  }. $

\medskip
Going through the above argument once more, we get that the coefficient of
$ U_{i_{1}} \ldots U_{i_{k}}  m_{  \sub{t}^{\bs{\lambda}} \sub{t}^{\bs{\lambda}}} $
in $\psi_{i_1}^b\ldots \psi_{i_k}^b   e^b(\bs{i}^{\bs{\lambda}}) $ is nonzero, in fact
it is essentially the product of the constant terms of the polynomials $ Q$ appearing in (\ref{definicion de psi}).
But Lemma \ref{pro 1} implies, by the choice of reduced expression for $\de{s} = s_{i_{1}}\ldots s_{i_{k}} $,
that $ U_{i_{1}} \ldots U_{i_{k}}  m_{  \sub{t}^{\bs{\lambda}} \sub{t}^{\bs{\lambda}}} = y_e^{l}m_{\sub{s} \sub{t}^{\bs{\lambda}}}
$ for some $l\in \mathbb{Z}\geq 0$
and then Lemma \ref{ayuda} implies that
$$ U_{i_{j_1}} \ldots U_{i_{j_r}}  m_{  \sub{t}^{\bs{\lambda}} \sub{t}^{\bs{\lambda}}}
= r m_{\sub{u} \sub{t}^{\bs{\lambda}}} \mod{b_n^{\bs{\lambda}}(m)  } $$
for some scalar $ r \in \C^{\times} $ and some $ \sub{u} $ such that $ \sub{u}\succ \sub{s}$.
Summing up,
this proves the Theorem in the case where $\sub{t} = \sub{t}^{\bs{\lambda}}$.

\medskip
To prove the general case, we first note that the same argument as above, only acting on the right instead of
on the left, proves the Theorem in the case where $ \sub{s} = \sub{t}^{\bs{\lambda}}$.
The general case then follows by multiplying the two versions together and using cellularity.
\end{dem}

\begin{rem} \rm
It follows from the Theorem that the subalgebra of $ \blob{m} $ generated by the $e^b(\bs{i}) $ and the $ \psi_i^b$
is equal to $ \blob{m} $ itself.
\end{rem}

\medskip

To establish our main theorem we must define a degree function on the set of all one-line standard bitableaux.
Let $\bs{\lambda} \in \bi{n}$ and $\sub{t} \in \Std(\bs{\lambda})$. Then we define the degree of $\sub{t} $ as
\begin{equation} \label{grado}
 \deg \sub{t} :=  \deg \psi_{\sub{tt}^{\bs{\lambda}}}.
\end{equation}
We can now prove our main result, namely to construct a graded cellular basis for $\blob{m}$.
Given our previous work, we can essentially follow the argument of \cite[Theorem 5.8]{hu-mathas}, just
making the corresponding changes in notation. We sketch the argument because this is the main theorem of the paper.
\begin{teo} \label{teo graded cellular bases}
The blob algebra $\blob{m}$ is a graded cellular algebra with graded cellular basis
$\{ \psi^b_{\sub{st}}  \mid \sub{s,t} \in \Std(\bs{\lambda}) \text{ for } \bs{\lambda} \in \bi{n}  \}$.
\end{teo}
\begin{dem}
First of all it follows from the triangularity property of Theorem \ref{transicion} that
$$\{ \psi^b_{\sub{st}}  \mid \sub{s,t} \in \Std(\bs{\lambda}) \text{ for } \bs{\lambda} \in \bi{n}  \}$$
is a cellular basis for $\blob{m}$, since
$\{ m_{\sub{st}}  \mid \sub{s,t} \in \Std(\bs{\lambda}) \text{ for } \bs{\lambda} \in \bi{n}  \}$ is it.
Moreover, by the
definitions, $  \psi^b_{\sub{st}} $ is a homogeneous elements of $\blob{m}$ of degree
$$  \deg{\psi_{\sub{st}}} = \deg{\sub{s}} + \deg{\sub{t}}. $$
Using Corollary \ref{main} one sees that there is a unique anti-automorphism
$*$ of $ \blob{m}$ that fixes the generators $ \psi_i^b,   y_j^b $ and $   e^b(\textbf{i}) $.
Then by the definition it is clear that
$ \psi_{\sub{st}}^{*}=\psi_{\sub{ts}} $ and so the anti-automorphism induced by the basis $\{ \psi^b_{\sub{st}}\} $
coincides with $*$.  The Theorem is proved.
\end{dem}

\medskip

By Theorem \ref{transicion} the cell modules induced by the graded cellular bases $\{ \psi^b_{\sub{st}}\}$
agree with the cell modules induced by the diagram bases $\{  m_{\sub{st}} \}$, that is
the standard modules for $\blob{m}$. Therefore, Theorem \ref{teo graded cellular bases} gives
us our main goal, to grade the standard modules for $\blob{m}$.

\begin{rem} \rm
The existence of a graded cellular basis
for the blob algebra allows one to define graded decomposition numbers.  Recently, the first author has
succeeded in calculating these graded decomposition numbers (see \cite{David}).
\end{rem}

\medskip
For completeness, we give the analogous Theorem for the Temperley-Lieb algebra.
This proof relies here on Theorem \ref{first_theorem} and the
compatibility of Murphy's standard basis with the diagram basis,
as proved in \cite{Martin}, and could have been given earlier in the paper.
Let $ \Phi_2 : \mathcal{H}_n(q^2)  \longrightarrow   Tl_{n}(q) $ as in Theorem \ref{TLisomorphism}. Define $\psi^{Tl}_{\sub{st}} := \Phi_2(\psi_{\sub{st}})$ for
$ \sub{s}, \sub{t} \in \Std(n) $, where $ \psi_{\sub{st} } $ is an element of the graded cellular basis for $  \mathcal{H}_n(q^2) $
introduced by Hu and Mathas \cite[Definition 5.1]{hu-mathas} and $ \Shape(\sub{s})= \Shape(\sub{t}) \in \pacol{n} $.
\begin{teo} \label{teo graded cellular bases tl}
The Temperley-Lieb algebra $Tl_{n}(q)$ is a graded cellular algebra with graded cellular basis
$\{ \psi^{Tl}_{\sub{st}} \} $ and degree function defined as above.
\end{teo}
\begin{dem}
According to \cite[Theorem 9]{Martin}, the diagram basis for $Tl_{n}(q)$ is upper triangularly
related to the Murphy's standard basis, with respect to the dominance order. But
$ \psi_{\sub{st}} $ is also upper triangularly related to the Murphy's standard basis with respect to the dominance order,
as already mentioned above, and hence the Theorem follows.
\end{dem}

\section{Examples.}
In this last section we illustrate our results on two examples.

\begin{exa} \rm
Our first example is $Tl_3(q)$, with $q$ chosen to be a primitive cubic root of unity, that is
$ l = 3$. This is a non-semisimple algebra and so we expect the grading to be nontrivial.
We determine the graded cellular basis $\psi^{Tl}_{\sub{st}} $
for $Tl_3(q)$, in terms of the diagrams.
Define first
\[
\xy 0;/r.2pc/:
(0,0)*{  \sub{s} = };
(5,0); (15,0) **\dir{-};
(5,5); (15,5) **\dir{-};
(5,-5); (5,5) **\dir{-};
(10,-5); (10,5) **\dir{-};
(5,-5); (10,-5) **\dir{-};
(15,5); (15,0) **\dir{-};
(7.5,2.5)*{1};
(12.5,2.5)*{2};
(7.5,-2.5)*{3};
(40,0)*{  \sub{t} = };
(45,0); (55,0) **\dir{-};
(45,5); (55,5) **\dir{-};
(45,-5); (45,5) **\dir{-};
(50,-5); (50,5) **\dir{-};
(45,-5); (50,-5) **\dir{-};
(55,5); (55,0) **\dir{-};
(47.5,2.5)*{1};
(52.5,2.5)*{3};
(47.5,-2.5)*{2};
\endxy \]
Then $ \sub{s} $ and $ \sub{t} $ are the only standard tableaux of shape $ (2,1) $.
The only other possible shape in $\pacol{3}$ is $ \lambda = (1,1,1) $ whose only standard tableau we denote by $ \sub{t}^{\lambda} $.
Hence we get that $Tl_3(q)$ has dimension five with homogeneous basis consisting of the elements
$$\psi^{Tl}_{\sub{ss}}, \psi^{Tl}_{\sub{st}} , \psi^{Tl}_{\sub{ts}}, \psi^{Tl}_{\sub{tt}}, \psi^{Tl}_{\sub{t^{\lambda}
\sub{t^{\lambda}}}}. $$

The residue sequences for $\sub{t}^{{\lambda}}$, $ \sub{s}  $ and $ \sub{t}  $ are $\bs{i}^{{\lambda}}=(0,2,1)$, $ \bs{i}^{\sub{s}}  = (0,1,2)$
and  $ \bs{i}^{\sub{t}}  = (0,2,1) $
and the degrees are $\deg(\sub{t}^{{\lambda}})=0$, $\deg(\sub{s} )=0$ and $\deg(\sub{t} )=1$. (See \cite[(3.8) and Definition 4.7]{hu-mathas}). Therefore,
using the orthogonality of the KLR-idempotents, we have
\begin{equation} \label{psiporpsi igual 0}
\psi^{Tl}_{\sub{st} } \psi^{Tl}_{\sub{ss} }= \psi^{Tl}_{\sub{ss} } \psi^{Tl}_{ \sub{ts} } =0,
\end{equation}
see \cite[Lemma 5.2]{hu-mathas}. We also have
\begin{equation} \label{psicuadrade igual psi}
\psi^{Tl}_{\sub{ss} } \psi^{Tl}_{\sub{ss} }= \psi^{Tl}_{\sub{ss} } = e(\bs{i}^{\sub{s}})
\end{equation}

Now by the triangular expansion property mentioned in the proof of the above Theorem \ref{teo graded cellular bases tl}
there exists $c\in \C^{\times}$ such that
\begin{center}
$ \psi^{Tl}_{\sub{ss}}=c $ $   \begin{array}{c}
    \blobvv
     \end{array} $
\end{center}

By (\ref{psicuadrade igual psi}) and the relation $U_{1}^{2}=-[2]U_1$ it is straightforward to check that $c= -1/[2]$. Now, using the triangular expansion property once again, there are scalars $c_1,c_2 \in \F $ with $c_1\neq 0$ such that

\begin{center}
$\psi^{Tl}_{\sub{st}}=c_1  $ $   \begin{array}{c}
    \blobvu
     \end{array} $ $+$ $c_2$ $   \begin{array}{c}
    \blobvv
     \end{array} $
\end{center}

Multiplying this equality on the right by $\psi^{Tl}_{\sub{ss} } = -\left( \frac{1}{[2]} \right)U_1$,
and using equation (\ref{psiporpsi igual 0}), we get
that $c_1=[2]c_2$. Hence the element

\begin{center}
$A:=[2]  $ $   \begin{array}{c}
    \blobvu
     \end{array} $ $+$  $   \begin{array}{c}
    \blobvv
     \end{array} $
\end{center}
is a scalar multiple of $\psi^{Tl}_{\sub{st} }$ and 
homogeneous of degree $1$. We are not able to 
determine explicitly the value of the scalar 
relating $\psi^{Tl}_{\sub{st} }$ and $A$.\\
 
 Similarly we obtain that the element
\begin{center}
$B:=[2]  $ $   \begin{array}{c}
    \blobuv
     \end{array} $ $+$  $   \begin{array}{c}
    \blobvv
     \end{array} $
\end{center}
is a scalar multiple of $\psi^{Tl}_{\sub{ts} }$ and homogeneous of degree $1$.

\medskip
Now, it is straightforward to check
$\psi^{Tl}_{\sub{ts} }
\psi^{Tl}_{\sub{st} } = \psi^{Tl}_{\sub{tt} }$. From this we obtain that $\psi^{Tl}_{\sub{tt} }$ is a scalar multiple of

\begin{center}
$C:=$ $   \begin{array}{c}
    \blobuu
     \end{array} $ $+$ $[2]$ $   \begin{array}{c}
    \blobvu
     \end{array} $ $+$ $[2]$ $   \begin{array}{c}
    \blobuv
     \end{array} $ $+$ $   \begin{array}{c}
    \blobvv
     \end{array} $
\end{center}
which is a 
homogeneous element of degree $2$. The last 
basis element, $\psi^{Tl}_{\sub{t^{\lambda}\sub{t^{\lambda}}}}$, can now 
be determined by expanding it 
in the diagram basis. On the other hand, by (\ref{kl3}) we have
\begin{equation*}
1=e(\bs{i}^{\sub{s}})+e(\bs{i}^{\lambda})=\psi^{Tl}_{\sub{ss}}+ \psi^{Tl}_{\sub{t^{\lambda}\sub{t^{\lambda}}}}
\end{equation*}
since $e(\bs{i}^{\sub{s}})$ and $e(\bs{i}^{\lambda})$ are the only non-zero KLR-idempotents. Therefore,
\begin{center}
$\psi^{Tl}_{\sub{t^{\lambda}\sub{t^{\lambda}}}} = $ $\begin{array}{c}
    \blobkk
     \end{array}$ $+$ $\frac{1}{[2]}$ $\begin{array}{c}
\blobvv
\end{array}$
\end{center}
All in all, the set $\{ \psi^{Tl}_{\sub{ss}}, A, B, C,\psi^{Tl}_{\sub{t^{\lambda}\sub{t^{\lambda}}}}\}$ is a graded cellular basis for $Tl_3(q)$. In particular,
$Tl_3(q)$ is a positively graded algebra and
$ {\cal F}_1 := \Span_{\C} \{ A, B, C \} $ and $ {\cal F}_2 := \Span_{\C} \{  C \} $
are ideals in $Tl_3(q)$. In general, $Tl_n(q)$ is not positively graded.
\end{exa}
\begin{exa} \rm
We now describe the graded cellular basis $\{\psi_{st}^{b}\}$ for  $b_3=b_3(q,y_e)$ in terms of blob diagrams, with $q$ a primitive quintic root of unity and $y_e=-\frac{[1]}{[2]}$, so in this case $l=5$ and $m=2$. First, we list all elements in $\text{Std} (3)$, with their respective residues sequences  and  degrees.

\[
\begin{tabular}{||c|c|c|c||}
  \hline
   \text{Bi-partitions} & \text{Bitableaux}   & \text{Res. Sequence} & \text{Degree}  \\ \hline
    &
\xy 0;/r0.60pc/:
(2,1.5)*{ };
(-2.3,.5)*{ \sub{t}^{\bs{\lambda}=}};
(-.5,.6)*{(};
(0,0); (1,0) **\dir{-};
(0,0); (0,1) **\dir{-};
(1,0); (1,1) **\dir{-};
(0,1); (1,1) **\dir{-};
(.5,.5)*{ \scriptstyle 2};
(1.5,0)*{,};
(2,0); (4,0) **\dir{-};
(2,1); (4,1) **\dir{-};
(2,0); (2,1) **\dir{-};
(3,0); (3,1) **\dir{-};
(4,0); (4,1) **\dir{-};
(4.5,.5)*{)};
(2.5,.5)*{ \scriptstyle 1};
(3.5,.5)*{ \scriptstyle 3};
\endxy
   & $\bs{i}^{\bs{\lambda}}=(4,1,0)$ & $0$  \\ \cline{2-4}
   $\bs{\lambda}=((1),(2))$  &
\xy 0;/r0.60pc/:
(2,1.5)*{ };
(-2,.5)*{ \sub{s}=};
(-.5,.5)*{(};
(0,0); (1,0) **\dir{-};
(0,0); (0,1) **\dir{-};
(1,0); (1,1) **\dir{-};
(0,1); (1,1) **\dir{-};
(.5,.5)*{ \scriptstyle 3};
(1.5,0)*{,};
(2,0); (4,0) **\dir{-};
(2,1); (4,1) **\dir{-};
(2,0); (2,1) **\dir{-};
(3,0); (3,1) **\dir{-};
(4,0); (4,1) **\dir{-};
(4.5,.5)*{)};
(2.5,.5)*{ \scriptstyle 1};
(3.5,.5)*{ \scriptstyle 2};
\endxy
 & $\bs{i}^{\sub{s}}=(4,0,1)$ & $ 1$ \\  \cline{2-4}
      &
 \xy 0;/r0.60pc/:
(2,1.5)*{ };
(-2,.5)*{ \sub{t}= };
(-.5,.5)*{(};
(0,0); (1,0) **\dir{-};
(0,0); (0,1) **\dir{-};
(1,0); (1,1) **\dir{-};
(0,1); (1,1) **\dir{-};
(.5,.5)*{ \scriptstyle 1};
(1.5,0)*{,};
(2,0); (4,0) **\dir{-};
(2,1); (4,1) **\dir{-};
(2,0); (2,1) **\dir{-};
(3,0); (3,1) **\dir{-};
(4,0); (4,1) **\dir{-};
(4.5,0.5)*{)};
(2.5,.5)*{ \scriptstyle 2};
(3.5,.5)*{ \scriptstyle 3};
\endxy   & $\bs{i}^{\sub{t}}=(1,4,0)$ & $0$ \\ \hline
       &  \xy 0;/r0.60pc/:
(2,1.5)*{ };
(-2.5,.5)*{ \sub{t}^{\bs{\mu}}= };
(-.5,.5)*{(};
(0,0); (2,0) **\dir{-};
(0,0); (0,1) **\dir{-};
(0,1); (2,1) **\dir{-};
(1,0); (1,1) **\dir{-};
(2,0); (2,1) **\dir{-};
(.5,.5)*{ \scriptstyle 2};
(1.5,.5)*{ \scriptstyle 3};
(2.5,0)*{,};
(3,0); (4,0) **\dir{-};
(3,1); (4,1) **\dir{-};
(3,0); (3,1) **\dir{-};
(4,0); (4,1) **\dir{-};
(4.5,0.5)*{)};
(3.5,.5)*{ \scriptstyle 1};
\endxy & $\bs{i}^{\bs{\mu}}= (4,1,2)$ & $0$ \\ \cline{2-4}
     $\bs{\mu}=((2),(1))$ &  \xy 0;/r0.60pc/:
(2,1.5)*{ };
(-2,.5)*{ \sub{v}= };
(-.5,.5)*{(};
(0,0); (2,0) **\dir{-};
(0,0); (0,1) **\dir{-};
(0,1); (2,1) **\dir{-};
(1,0); (1,1) **\dir{-};
(2,0); (2,1) **\dir{-};
(.5,.5)*{ \scriptstyle 1};
(1.5,.5)*{ \scriptstyle 3};
(2.5,0)*{,};
(3,0); (4,0) **\dir{-};
(3,1); (4,1) **\dir{-};
(3,0); (3,1) **\dir{-};
(4,0); (4,1) **\dir{-};
(4.5,0.5)*{)};
(3.5,.5)*{ \scriptstyle 2};
\endxy & $\bs{i}^{\sub{v}}=(1,4,2)$   & $0$ \\  \cline{2-4}
         & \xy 0;/r0.60pc/:
(2,1.5)*{ };
(-2,.5)*{ \sub{u}= };
(-.5,.5)*{(};
(0,0); (2,0) **\dir{-};
(0,0); (0,1) **\dir{-};
(0,1); (2,1) **\dir{-};
(1,0); (1,1) **\dir{-};
(2,0); (2,1) **\dir{-};
(.5,.5)*{ \scriptstyle 1};
(1.5,.5)*{ \scriptstyle 2};
(2.5,0)*{,};
(3,0); (4,0) **\dir{-};
(3,1); (4,1) **\dir{-};
(3,0); (3,1) **\dir{-};
(4,0); (4,1) **\dir{-};
(4.5,0.5)*{)};
(3.5,.5)*{ \scriptstyle 3};
\endxy & $\bs{i}^{\sub{u}}= (1,2,4)$ & $0$ \\ \hline
     $\bs{\nu}=((0),(3))$   & \xy 0;/r0.60pc/:
(2,1.5)*{ };
(-2.5,.5)*{ \sub{t}^{\bs{\nu}}= };
(-.5,.5)*{(};
(0,.5)*{\emptyset};
(.5,0)*{,};
(1,0); (4,0) **\dir{-};
(1,0); (1,1) **\dir{-};
(1,1); (4,1) **\dir{-};
(2,0); (2,1) **\dir{-};
(3,0); (3,1) **\dir{-};
(4,0); (4,1) **\dir{-};
(4.5,.5)*{)};
(1.5,.5)*{ \scriptstyle 1};
(2.5,.5)*{ \scriptstyle 2};
(3.5,.5)*{ \scriptstyle 3};
\endxy  & $\bs{i}^{\bs{\nu}}= (4,0,1) $ &  $0$ \\ \hline
      $\bs{\kappa}=((3),(0))$ & \xy 0;/r0.60pc/:
(2,1.5)*{ };
(-1.5,.5)*{ \sub{t}^{\bs{\kappa}}= };
(0.5,.5)*{(};
(1,0); (4,0) **\dir{-};
(1,0); (1,1) **\dir{-};
(1,1); (4,1) **\dir{-};
(2,0); (2,1) **\dir{-};
(3,0); (3,1) **\dir{-};
(4,0); (4,1) **\dir{-};
(4.5,0)*{,};
(5,.5)*{\emptyset};
(5.5,.5)*{)};
(1.5,.5)*{ \scriptstyle 1};
(2.5,.5)*{ \scriptstyle 2};
(3.5,.5)*{ \scriptstyle 3};
\endxy   & $\bs{i}^{\bs{\kappa}}= (1,2,3) $ &  $0$ \\
  \hline
\end{tabular}
\]

\medskip \medskip \medskip
\noindent
We need the following Lemma.
\begin{lem} \label{e is homogeneous}
Let $k \in \mathbb{Z}$  such that $2k \equiv m \mod{l} $. Then the element $e \in \blob{m} $ is
homogeneous of degree zero. More precisely, it can be written as a sum of homogeneous elements of degree zero as follows.
\begin{equation}
e=\sum_{\substack{
   \bs{i}\in I^n \\
   i_{1}=-k
  }} e(\bs{i})
\end{equation}
Furthermore, for all $\sub{s},\sub{t} \in \text{Std}(n)$ we have
$$  e\psi_{\sub{st}}^{b}  = \left\{
      \begin{array}{cl}
       \psi_{\sub{st}}^{b} , & \hbox{if $1$ is located in the second component of $\sub{s}$,} \\
        0, & \hbox{otherwise.}
      \end{array}
    \right.
 $$
$$  \psi_{\sub{st}}^{b} e = \left\{
      \begin{array}{cl}
       \psi_{\sub{st}}^{b} , & \hbox{if $1$ is located in the second component of $\sub{t}$,} \\
        0, & \hbox{otherwise.}
      \end{array}
    \right.
 $$
\end{lem}

\begin{dem}
The first claim follows by combining  Theorem \ref{isomorphism}, (\ref{Lsubr}) and $U_0=-[m]e$. The second is a direct consequence of the definition of $\psi_{\sub{st}}^{b}$, the orthogonality of  KLR-idempotents and the first claim.
\end{dem}

\medskip
Using  the triangularity property given in Lemma \ref{transicion}, the orthogonality of the KLR-idempotents and the previous Lemma \ref{e is homogeneous}, we can now give a description of the graded cellular bases $\{\psi_{st}^{b}\}$ of $b_3$ in terms of the diagrammatic basis. We omit the details for brevity, since $\text{dim}_{\C} (b_3) $ is quite big, of
dimension 20. The scalars $r_{\sub{ab}}$ appearing in this diagrammatic description correspond to the non-zero scalar $r$ appearing in Theorem \ref{transicion}. For brevity we omit some of the elements of the 
basis $\{\psi_{\sub{st}}\}$, but one can obtain the diagrammatic expression of the elements not 
enlisted by multiplying two of the enlisted elements. For example, $\psi_{\sub{ts}}$ is 
not enlisted but $\psi_{\sub{ts}}=\psi_{\sub{t}\sub{t}^{\bs{\lambda}}}\psi_{\sub{t}^{\bs{\lambda}}\sub{s}}$, and 
the elements $\psi_{\sub{t}\sub{t}^{\bs{\lambda}}}$ and $\psi_{\sub{t}^{\bs{\lambda}}\sub{s}}$ are enlisted. 
We finally remark that, just like in the Temperley-Lieb algebra case, in general the blob algebra is 
not positively graded.\\\\

   $\psi_{\sub{t}^{\bs{\lambda}}\sub{t}^{\bs{\lambda}}}^{b}= \frac{1}{y_e}$ $   \begin{array}{c}
    \bloblala
     \end{array} $   

\vspace{.5cm}
 
 $\psi_{\sub{t}^{\bs{\lambda}}\sub{t}}^{b}= r_{\sub{t}^{\bs{\lambda}}\sub{t}} $  $  \left( \begin{array}{c}
    \bloblat 
     \end{array} \right. $ $-$  $  \left. \begin{array}{c}
    \bloblala
     \end{array} \right) $

\vspace{.5cm}

 $\psi_{\sub{t}\sub{t}^{\bs{\lambda}}}^{b}= r_{\sub{t}\sub{t}^{\bs{\lambda}}} $  $  \left( \begin{array}{c}
    \blobtla
     \end{array} \right. $ $-$  $  \left. \begin{array}{c}
    \bloblala
     \end{array} \right) $

\vspace{.5cm}

$\psi_{\sub{s}\sub{t}^{\bs{\lambda}}}^{b}= r_{\sub{s}\sub{t}^{\bs{\lambda}}} $  $  \left( \begin{array}{c}
    \blobsla
     \end{array} \right. $ $-$ $ \frac{1}{y_e}$  $  \left. \begin{array}{c}
    \bloblala
     \end{array} \right) $

\vspace{.5cm}

$\psi_{\sub{t}^{\bs{\lambda}}\sub{s}}^{b}= r_{\sub{t}^{\bs{\lambda}}\sub{s}} $  $  \left( \begin{array}{c}
    \bloblas
     \end{array} \right. $ $-$ $ \frac{1}{y_e}$  $  \left. \begin{array}{c}
    \bloblala
     \end{array} \right) $

\vspace{.5cm}

$\psi_{\sub{t}^{\bs{\mu}}\sub{t}^{\bs{\mu}}}^{b}= \frac{1}{y_e} $  $  \left( \begin{array}{c}
    \blobmumu
     \end{array} \right. $ $-$  $  \left. \begin{array}{c}
    \bloblala
     \end{array} \right) $

\vspace{.5cm}

$\psi_{\sub{t}^{\bs{\mu}}\sub{v}}^{b}= r_{\sub{t}^{\bs{\mu}}\sub{v}} $  $  \left( \begin{array}{c}
    \blobmuv
     \end{array} \right. $ $-$  $ \begin{array}{c}
    \blobmumu
     \end{array} $ $-$  $ \begin{array}{c}
    \bloblat
     \end{array} $ $+$ $  \left. \begin{array}{c}
    \bloblala
     \end{array} \right) $

\vspace{.5cm}

$\psi_{\sub{v} \sub{t}^{\bs{\mu}}}^{b}= r_{\sub{v} \sub{t}^{\bs{\mu}}} $  $  \left( \begin{array}{c}
    \blobvmu
     \end{array} \right. $ $-$  $ \begin{array}{c}
    \blobmumu
     \end{array} $ $-$  $ \begin{array}{c}
    \blobtla
     \end{array} $ $+$ $  \left. \begin{array}{c}
    \bloblala
     \end{array} \right) $

\vspace{.5cm}

$\psi_{ \sub{t}^{\bs{\mu}}\sub{u}}^{b}= r_{ \sub{t}^{\bs{\mu}}\sub{u}} $  $  \left( \begin{array}{c}
    \blobmuu
     \end{array} \right. $ $-$  $ \begin{array}{c}
    \blobmuv
     \end{array} $ $+$
$ \begin{array}{c}
    \blobmumu
     \end{array} $ \\
\hspace*{2.1cm}
$-$ 
 $ \begin{array}{c}
    \bloblas
     \end{array} $ $+$ 
$ \begin{array}{c}
    \bloblat
     \end{array} $
$-$ $  \left. \begin{array}{c}
    \bloblala
     \end{array} \right) $

\vspace{.5cm}

$\psi_{\sub{u} \sub{t}^{\bs{\mu}}}^{b}= r_{\sub{u} \sub{t}^{\bs{\mu}}} $  $  \left( \begin{array}{c}
    \blobumu
     \end{array} \right. $ $-$  $ \begin{array}{c}
    \blobvmu
     \end{array} $ $+$
$ \begin{array}{c}
    \blobmumu
     \end{array} $ \\
\hspace*{2.1cm}$-$
 $ \begin{array}{c}
    \blobsla
     \end{array} $ $+$
$ \begin{array}{c}
    \blobtla
     \end{array} $
$-$ $  \left. \begin{array}{c}
    \bloblala
     \end{array} \right) $

\vspace{.5cm}

$\psi_{\sub{t}^{\bs{\nu}}\sub{t}^{\bs{\nu}}}^{b}=  $  $  \begin{array}{c}
    \blobnunu
     \end{array} $ $-$ $\frac{1}{y_e}$ $   \begin{array}{c}
    \blobmumu
     \end{array}  $

\vspace{.5cm}

$\psi_{\sub{t}^{\bs{\kappa}}\sub{t}^{\bs{\kappa}}}^{b}=  $  $  \begin{array}{c}
    \blobkk
     \end{array}  $ $-$  $ \begin{array}{c}
    \blobnunu
     \end{array} $ $+$ $\left(  \frac{1}{1+[2]} \right)$ 
$\left( \begin{array}{c}
    \blobts
     \end{array} \right. $ $-$
$ \begin{array}{c}
    \bloblas
     \end{array} $
\\
\hspace*{2.5cm} $-$
$\begin{array}{c}
    \blobss
     \end{array}  $
$-$ $   \begin{array}{c}
    \blobvu
     \end{array}  $
$+$ $   \begin{array}{c}
    \blobmuu
     \end{array}  $
$-$ $  \left. \begin{array}{c}
    \blobuu
     \end{array} \right) $\\
\hspace*{1.0cm} $-\left(\frac{1}{1-[2]}\right)$
$\left( \begin{array}{c}
    \bloblala
     \end{array} \right. $
$-$ $   \begin{array}{c}
    \bloblat
     \end{array}  $
$-$ $   \begin{array}{c}
    \blobtla
     \end{array}  $
$+$ $   \begin{array}{c}
    \blobtt
     \end{array}  $\\
\hspace*{2.5cm} $-$
$\begin{array}{c}
    \blobst
     \end{array}  $
$+$ $   \begin{array}{c}
    \blobsla
     \end{array}  $
$+$ $   \begin{array}{c}
    \blobuv
     \end{array}  $
$-$ $  \left. \begin{array}{c}
    \blobumu
     \end{array} \right) $\\
\hspace*{1cm} $-$
$\left( \frac{[2]}{[2]-1} \right)$
$ \left( \begin{array}{c}
    \blobvmu
     \end{array} \right.  $
$-$ $   \begin{array}{c}
    \blobvv
     \end{array}  $
$+$ $   \begin{array}{c}
    \blobmuv
     \end{array}  $
$-$ $  \left. \begin{array}{c}
    \blobmumu
     \end{array} \right) $

\end{exa}

\begin{scriptsize}

\sc Instituto de Matem\'atica y F\'isica, Universidad de Talca, Chile, mail dplaza@inst-mat.utalca.cl,
steen@inst-mat.utalca.cl,
\end{scriptsize}


\begin{thebibliography}{99}

\bibitem{BGIL} C. Bonnaf\'e, M. Geck, L. Iancu and T. Lam,  On domino insertion and Kazhdan-Lusztig cells in type $B_n$.
In: Representation theory of algebraic groups and quantum groups
(Nagoya, 2006; eds. A. Gyoja et al.), p. 33--54, Progress in Math. {\bf 284}, Birkh\"auser,
2010.
\bibitem{brundan-klesc} J. Brundan and A. Kleshchev,
Blocks of cyclotomic Hecke algebras and Khovanov-Lauda algebras, Invent. Math. {\bf 178} (2009), 451-484.
\bibitem{brundan-klesc1} J. Brundan and A. Kleshchev, Graded decomposition numbers for cyclotomic Hecke algebras, Adv.
Math., {\bf 222} (2009), 188-942.
\bibitem{bkw graded} J. Brundan,  A. Kleshchev and W.  Wang,   Graded   Specht   modules,
J. Reine und Angew. Math. {\bf 655} (2011), 61-87
\bibitem{chevalley} C. Chevalley, The Construction and Study of Certain Important Algebras. Mathematical Society of Japan (1955).
\bibitem{blob positive} A. Cox, J. Graham and P. Martin (2003). The blob algebra in positive characteristic.
Journal of Algebra, {\bf 266}(2), 584 - 635.
\bibitem{cps} E. Cline, B. Parshall and L.Scott, Finite dimensional algebras and highest weight categories,
Math. Ann., {\bf 259} (1982), 153-199.
\bibitem{gra ler}  J. Graham and G. Lehrer, Cellular algebras, Inventiones Mathematicae, {\bf 123} (1996), 1-34.
\bibitem{gra ler 1}  J. Graham and G. Lehrer,  Cellular algebras and diagram algebras in representation theory,
Advanced Studies in Pure Mathematics, {\bf 40} (2004), 141-173.
\bibitem{Martin}M. H\"arterich, Murphy Bases of generalized Temperley-Lieb algebras. Arch.Math. {\bf 72} (1999), 337-345.
\bibitem{hu-mathas} J. Hu, A. Mathas, Graded cellular bases for the cyclotomic Khovanov-Lauda-Rouquier
algebras of type A, Adv. Math., {\bf 225} (2010), 598-642.
\bibitem{kle} A. Kleshchev, Linear and Projective Representations of Symmetric Groups,
 Cambridge University Press, Cambridge, 2005.
\bibitem{Mat-Sal} P. P. Martin and H. Saleur, The blob algebra and the periodic Temperley-Lieb algebra,
Lett. Math. Phys. {\bf 30} (1994), 189-206.
\bibitem{martin-wood}P. P. Martin,  D. Woodcock, Generalized blob algebras and alcove geometry, LMS Journal
of Computation and Mathematics {\bf 6}, (2003), 249-296.
\bibitem{martin-wood1}P. P. Martin, D. Woodcock, On the structure of the blob algebra, J. Algebra {\bf 225} (2000), 957-988.
\bibitem{Mathas} A. Mathas, Hecke algebras and Schur algebras of the symmetric group, Univ. Lecture Notes, 15,
Amer. Math. Soc., 1999.
\bibitem{Mat-So} A. Mathas, Seminormal forms and Gram determinants for cellular algebras, J.
Reine Angew. Math., {\bf 619} (2008), 141-173.  With an appendix by M. Soriano.
\bibitem{Murphy} G. E. Murphy, The idempotents of the symmetric group and
Nakayama's conjecture, J. Algebra {\bf 81} (1983), 258-265.
\bibitem{Murphy1} G. E. Murphy, The Representations of Hecke Algebras of type
$ A_n $, Journal of Algebra {\bf 173} (1995), 97--121.
\bibitem{David} D. Plaza, Graded decomposition numbers for the blob algebra, Journal of Algebra {\bf 394} (2013), 182-206.
\bibitem{Steen} S. Ryom-Hansen, The Ariki-Terasoma-Yamada tensor space and the
blob-algebra, J. of Algebra {\bf 324} (2010), 2658-2675.
\bibitem{Steen1} S. Ryom-Hansen, Cell structures on the blob algebra, Representation Theory {\bf 16} (2012), 540-567. 
\bibitem{Stroppel} C. Stroppel,
Categorification of the Temperley-Lieb
category, tangles, and cobordism via
projective functors,
Duke Mathematical Journal
{\bf 126}, No. 3, 2005.
\bibitem{Zhang} R. B. Zhang,
Graded representations of the Temperley-Lieb algebra, quantum supergroups, and the Jones polynomial,
J. Math. Phys. {\bf 32}, 2605 (1991).
\end{thebibliography}
\end{document}